\documentclass[preprint,12pt]{elsarticle}

\newtheorem{theorem}{Theorem}
\newtheorem{definition}{Definition}
\newtheorem{lemma}{Lemma}
\newtheorem{corollary}{Corollary}
\newdefinition{remark}{Remark}
\newproof{proof}{Proof}
\newproof{pot}{Proof of Theorem \ref{thm2}}
\newtheorem{example}{Example}

\newcommand{\be}{\begin{equation}}
\newcommand{\ee}{\end{equation}}
\newcommand{\field}[1]{\mathbb{#1}}
\newcommand{\C}{\field{C}}
\newcommand{\R}{\field{R}}
\newcommand{\N}{\field{N}}

\newcommand{\A}{\field{A}}

\usepackage{graphicx}
\usepackage{float}
\usepackage{amssymb}
\usepackage{amsmath}
\usepackage{latexsym}
\usepackage{cancel}
\usepackage{color}
\usepackage{placeins}
\usepackage{url}
\renewcommand{\restriction}{\mathord{\upharpoonright}}
\usepackage{BOONDOX-cal}

\journal{arXiv}

\begin{document}

\begin{frontmatter}



\title{The Radon transform and the Hough transform: a unifying perspective}

\author[label1]{Riccardo Aramini\corref{cor1}\fnref{corr}}
\ead{aramini@dima.unige.it}
\fntext[corr]{Corresponding author.}

\author[label2]{Fabrice Delbary}
\ead{fdelbary@uni-mainz.de}

\author[label1]{Mauro C. Beltrametti}
\ead{beltrametti@dima.unige.it}

\author[label1,label3]{Michele Piana}
\ead{piana@dima.unige.it}

\author[label3]{Anna Maria Massone}
\ead{annamaria.massone@cnr.it}

\address[label1]{Dipartimento di Matematica, Universit\`a di Genova, via Dodecaneso 35, I-16146 Genova, Italy}
\address[label2]{Institut f{\"u}r Mathematik, Johannes Gutenberg-Universit{\"a}t Mainz, Staudingerweg 9, 55128 Mainz, Germany}
\address[label3]{CNR - SPIN, Genova, via Dodecaneso 33, I-16146 Genova, Italy}

\begin{abstract}
The Radon transform is a linear integral transform that mimics the data formation process in medical imaging modalities like
X-ray Computerized Tomography and Positron Emission Tomography. The Hough transform is a pattern recognition technique, which
is mainly used to detect straight lines in digital images and which has been recently extended to the automatic recognition of algebraic plane curves. Although defined in very different ways, in numerical applications both transforms ultimately take an image as an input and provide, as an output, a function defined on a parameter space. The parameters in this space describe
a family of curves, which represent either the integration domains considered in the (generalized) Radon transform, or the
curves to be detected by means of the Hough transform. In both cases, the 2D plot of the intensity values of the output function is the so-called (Radon or Hough) \textit{sinogram.} While the Hough sinogram is produced by an algorithm whose implementation
requires that the parameter space be discretized in cells, the Radon sinogram is mathematically defined on a continuous
parameter space, which in turn may need to be discretized just for physical or numerical reasons. In this paper, by considering a more general and $n$-dimensional setting, we prove that, whether the input image is described as a set of points (possibly with different intensity values) or as a piecewise constant function, its (rescaled) Hough sinogram converges to the corresponding Radon sinogram as the discretization step in the parameter space tends to zero. We also show that this result may have a notable impact on the image reconstruction problem of inverting the Radon sinogram recorded by a medical imaging scanner, and that the description of the Hough transform problem within the framework of regularization theory for inverse problems is worth investigating.
\end{abstract}


\end{frontmatter}

\section{Introduction}

The Radon transform \cite{ra17, helgason} is an important tool in harmonic analysis with significant conceptual impacts on both group theory and applied mathematics. The classical definition of this transform considers integrals over hyperplanes with specific orientation and distance from a reference hyperplane. For this classical Radon transform many functional properties have been investigated, including the characterization of its kernel and range, the ill-posedness of the inverse problem, as well as several inversion formulas and algorithms. The Radon transform has also been extended to integration on manifolds \cite{ku06}, although in this case important functional and computational problems are still open issues.

In biomedical imaging the classical Radon transform is the well-established mathematical model for data formation in X-ray Computerized Tomography (CT) and in Positron Emission Tomography (PET) \cite{na86, nawu01}. Indeed in X-ray CT the parameter that must be represented in the image is the density of the biological tissue, but the signal recorded by the scanner (the so-called sinogram) is a set of integrals of such density along straight lines with many different orientations and at many different distances from a reference line. On the other hand, a PET sinogram is the collection of line integrals of the concentration of a tracer that is injected into the body and whose interaction with the tissue represents a clinically sound metabolic index. Therefore, all software tools for image visualization implemented in current industrial CT and PET scanners must realize, at same stage, the numerical inversion of the Radon transform.

While the Radon transform plays a crucial role in image reconstruction, the Hough transform provides an important computational technique in pattern recognition. Indeed, the Hough transform is widely used in image-processing to detect algebraic plane curves, which are zero-loci of polynomials whose coefficients depend polynomially on a set of parameters. The basic idea of this recognition procedure (just extending the usual point-line duality in projective plane) is that a point in the image space corresponds to a locus (its Hough transform) in the parameter space. In turn, the whole curve in the image space corresponds {\it by duality} to a single point given by the intersection of all Hough transforms of the points belonging to the curve. A histogram (the Hough counter) can be constructed, representing an accumulator function defined on the discretized parameter space: for each cell in the parameter space, the value of the accumulator corresponds to the number of Hough transforms passing through that cell. The position of the maximum in the Hough counter identifies the combination of parameters characterizing the curve to be detected in the image space.

The history of the Hough transform starts in $1962$ with a patent by P. V. C. Hough \cite{ho62} to detect straight tracks of subatomic particles in bubble chamber photographs. No algebraic equations are used in the Hough patent, where the transform is defined just as {\it geometric construction by hand}. A first detailed description of the computational steps needed to implement the Hough transform technique, together with a theoretical generalization (although just outlined) to arbitrary curves, can then be found in \cite{duha72}.
At the beginning of the Nineties, a monograph \cite{le92} makes the point about the Hough transform from several perspectives: theory (in particular, generalizations, extensions and variants), numerics, applications, interpretations, future developments. In \cite{le00}, the Radon transform is applied to extract parameters characterizing the shape and angularity of powder particles: this application is somehow in the spirit of the Hough transform and contributes to highlight the link between the two transforms.

Recent papers \cite{feol08, liol15} introduce new algorithms based on Hough transform voting schemes, enabling very fast and efficient recognition of specific geometric features in large images or data sets. From a more theoretical point of view, a recent research \cite{bemapi13} provides a rigorous mathematical foundation, based on algebraic-geometry arguments, for the case of algebraic plane curves of whatever degree, together with a key lemma stating equivalent conditions under which the existence and uniqueness of the intersection point of the Hough transforms is guaranteed. This framework is then applied in \cite{macapebe15} and \cite{iciap15} to provide an atlas of algebraic curves used to recognize profiles in real astronomical and biomedical images.
Finally, we complete our short overview of the Hough transform by citing \cite{much15}, an up-to-date survey of this transform, its variants and applications. In the abstract of this paper it is claimed that more than $2,500$ research papers are concerned with the Hough transform, which represents an expression of uninterrupted interest from scholars in this field during the last decades.

In $1981$, for the first time an IEEE letter \cite{de81} guesses and shows by examples how the Hough transform can be considered a particular case of the Radon transform. Although influential and constructive, this letter is somewhat heuristic and does not consider any formal definition of the Hough transform. In \cite{prilki92}, the limitations of \cite{de81} are noticed and the similarity between the two transforms is investigated by relying on a formal definition of the Hough transform. However, both \cite{de81} and \cite{prilki92} fail to present a general and sound mathematical framework for studying in depth the relationship between the two transforms. In particular, some notation drawn from the theory of distributions is occasionally adopted, without any formal assumption and specification of the conditions making this notation mathematically meaningful. In $2004$, an inspirational report \cite{vaetal04} both reviews the literature about the relationship between the two transforms and outlines a sort of (mainly mathematical) research program to properly understand their link; in particular, it points out the importance of using concepts and results from distribution theory.

The aim of this paper is to present a general framework to describe and explain the relationship between the (generalized) Radon transform and the Hough transform. Specifically, our aim is to prove that given a digital image, the corresponding Hough counter tends to become the Radon transform of the image itself as the discretization of the parameter space becomes finer and finer.

Depending on the context, the spatial extent of a pixel in a two-di\-men\-sio\-nal image may be regarded as negligible or not.
If the pixel is considered as dimensionless, a mathematical model describing it can be chosen as the Dirac delta centered at a point,
multiplied by a number representing the grey level\footnote{\label{bellissimi}We recall that grey levels are a calibrated sequence of
grey tones, represented by integers and classified into grey-scale bands, ranging from black (usually, level $0$) to white
(usually, level $255$).}
or an analogous information about the intensity of the pixel itself. On the other hand,
if the pixel is assumed to take up a small square region, it can be mathematically described by a function being constant on the square and
zero outside. Accordingly, throughout the paper we shall speak of ``discrete image'' whenever the underlying mathematical model consists
of a linear combination (with real or even complex-valued coefficients) of Dirac deltas centered at a finite number of points in $\R^2$
(or, more generally, in $\R^n$); instead, we shall speak of ``piecewise continuous image'' if the corresponding mathematical description
is given in terms of a piecewise continuous (or, in particular, piecewise constant) function. With a slight abuse of language,
we shall often identify an image with its mathematical model.

The plan of the paper is as follows. In Section \ref{RT} we recall some basic notation and definitions concerning the Radon
transform of a (piecewise continuous) function $m$ describing an image, both in its traditional formulation
(as the set of all the surface integrals of $m$ over hyperplanes in $\R^n$, see e.g. \cite[chap.~1]{helgason})
and in a distributional framework (as inspired by \cite[chap.~I]{gegr5}), whereby the integral over a hyperplane is replaced
by the action of an appropriate distribution on $m$, regarded as a test function. Such distribution is the Dirac delta of the function
describing the hyperplane in Cartesian coordinates and is supported on the hyperplane itself. In this regard, \ref{distrapp} is devoted
to a short survey of some concepts and results of distribution theory, as needed and applied throughout the paper, with particular
attention to the definition of the Dirac delta of a function and its connection with the coarea formula.

In Section \ref{ziogen} the distributional definition of the Radon transform is generalized in such a way that hyperplanes can be
replaced by a $\lambda$-parametrized family of smooth submanifolds of $\R^n$, being $\lambda\in E\subset\R^t$ a $t$-dimensional parameter.
For each $\lambda\in E$, the corresponding submanifold is the zero locus in $\R^n$ of a continuously differentiable
function $f(\cdot;\lambda)$ expressible in the $\lambda_t$-solvable form $f(x;\lambda)=\lambda_t-F(x;\lambda_1,\ldots,\lambda_{t-1})$.
Some regularity results for this specific version of the so-called ``generalized Radon transform'' (cf. \cite{ku06} and references therein)
are established, as well as a physical interpretation allowing for a further extension, i.e., the determination of the generalized Radon
transform of a Dirac delta concentrated at a point in $\R^n$ and, by linearity, of any discrete image. The latter result, together with
a short analysis of the concept of ``sinogram'' as a visual representation of the intensity values of the generalized Radon transform,
is presented in Section \ref{sinogrammi}.

The Hough transform (in the case of a discrete image) is introduced in Section \ref{HT}. Here, we first recall some basic notions
and definitions, referring mainly to \cite{bemapi13, macapebe15} for several details and applications. Then, we describe the
discretization of the parameter space and define on it the weighted Hough counter, a function of crucial importance in the implementation
of any algorithm based on the Hough transform. Next, we focus on some important consequences of the $\lambda_t$-solvability property
for the function $f(x;\lambda)$ whose zero loci (either for a fixed $x$ or for a fixed $\lambda$) are at the basis of the whole Hough
transform process. In particular, we define the rescaled Hough counter as the
ratio between the weighted Hough counter and the solvable parameter $\lambda_t$, and we introduce the concept of ``rescaled
Hough sinogram'' as a visual representation of the intensity values of the rescaled Hough counter.

Section \ref{ziodiscreto} is concerned with the case of discrete images. Its main result is the theorem stating that, as the discretization
of the parameter space $E\subset\R^t$ becomes infinitely fine, the rescaled Hough counter, obtained for a given discrete image and for a
$\lambda_t$-solvable function $f(x;\lambda)$, tends (in a distributional sense) to the generalized Radon transform of the image itself,
provided that the latter transform is computed by integrating over submanifolds that are just the zero loci of $f(\cdot;\lambda)$ in $\R^n$
for any fixed $\lambda\in E$. Section \ref{ziocontinuo} extends the analysis and results of Section
\ref{ziodiscreto} to the case of piecewise continuous images. Sections \ref{ziodiscreto} and \ref{ziocontinuo} together represent the core
of the paper, since they provide a quite general framework for describing and explaining in detail the close but not evident
relationship between the Radon transform and the Hough transform.

Section \ref{zionumerico} presents a numerical example in which a digital phantom is recovered from a very noisy Radon sinogram,
by regarding it as a Hough sinogram.
Finally, in Section \ref{zioprospettico} we point out that the numerical technique just outlined in Section \ref{zionumerico}, if properly understood and implemented, might find an interesting application to all cases (like in Positron Emission Tomography) in which the Radon sinograms are inherently affected by a high level of noise, so that the traditional (i.e., Radon-based) inversion techniques cannot provide a satisfactory reconstruction of the unknown object.

In order to make the paper as readable and self-contained as possible, we added an Appendix recalling and collecting some notation, definitions, theorems and properties that are often used throughout the paper itself.

\section{The Radon transform}\label{RT}

Let $\gamma\in\R$ and $\widehat{\omega}\in\mathbb{S}^{n-1}:=\{x\in\R^n: |x|=1\}$, with $n\in\N\setminus\{0,1\}$. Then, we define the hyperplane
$\mathcal{P}(\widehat{\omega},\gamma)$ in $\R^n$ as
\be \label{pinocchio}
\mathcal{P}(\widehat{\omega},\gamma):=\left\{x\in\R^n : \gamma-\widehat{\omega}\cdot x=0 \right\},
\ee
where, of course, $x=(x_1,\ldots,x_n)$, $\widehat{\omega}=(\widehat{\omega}_1,\ldots,\widehat{\omega}_n)$ and the dot ``$\cdot$'' between two
elements of $\R^n$ denotes the canonical scalar product in $\R^n$.

\begin{definition}\label{defradon}
Let $m:\R^n\rightarrow\C$ be a function such that $m\in L^1(\mathcal{P}(\widehat{\omega},\gamma))$ $\forall (\widehat{\omega},\gamma)\in
\mathbb{S}^{n-1}\times\R$. Then, the \textnormal{Radon transform} of $m$ is defined as the function
$(Rm):\mathbb{S}^{n-1}\times\R\rightarrow\C$ given by
\be\label{radon}
(Rm)(\widehat{\omega},\gamma):=\int_{\mathcal{P}(\widehat{\omega},\gamma)}m(x)\,d\sigma(x)
\ \ \ \ \forall (\widehat{\omega},\gamma)\in \mathbb{S}^{n-1}\times\R,
\ee
where $d\sigma(x)$ is the Euclidean element of area on $\mathcal{P}(\widehat{\omega},\gamma)$ \textnormal{\cite{helgason}}.
\end{definition}

For each $(\widehat{\omega},\gamma)\in \mathbb{S}^{n-1}\times\R$, we consider the map defined by
\be\label{gloria1}
\R^n\ni x\mapsto f(x;\widehat{\omega},\gamma) := \gamma-\widehat{\omega}\cdot x \in \R,
\ee
and assume (just for notational simplicity) that $\widehat{\omega}_n\neq 0$. Then, we have
\be\label{gloria2}
\mathcal{P}(\widehat{\omega},\gamma)=\left\{x\in\R^n : f(x;\widehat{\omega},\gamma)=0 \right\}=
\left\{x\in\R^n : x_n=\mathsf{F}(x';\widehat{\omega}',\gamma) \right\},
\ee
where the notation
$\widehat{\omega}':=(\widehat{\omega}_1,\ldots,\widehat{\omega}_{n-1})$, $x':=(x_1,\ldots,x_{n-1})$ and
$\mathsf{F}(x';\widehat{\omega}',\gamma):=(\gamma -\widehat{\omega}'\cdot x')/\widehat{\omega}_n$ has been adopted.
Accordingly, by (\ref{dini}), (\ref{intsup}) and (\ref{desigma}) in the appendix,
expression (\ref{radon}) can be explicitly rewritten as
\be \label{detorrente}
(Rm)(\widehat{\omega},\gamma):=\frac{1}{|\widehat{\omega}_n|}\int_{\R^{n-1}}m\big(x',\mathsf{F}(x';\widehat{\omega}',\gamma)\big)\,dx'.
\ee

We also recall that, by assumption, $|\widehat{\omega}|=1$: then, from (\ref{gloria1}),
we have $|\mathrm{grad}\,f(x;\widehat{\omega},\gamma)|=1$
$\forall x\in\R^n$, i.e., condition (\ref{grad1}) is fulfilled. Hence, by (\ref{defdelta}), (\ref{desigma}) and (\ref{viotti}),
for a function\footnote{The space $\mathcal{PD}_0(\R^n)$ is the vector space
$PC^0_{C}(\R^n)$ of piecewise continuous and compactly supported functions,
endowed with an appropriate topology. The corresponding space of linear and continuous functionals on $\mathcal{PD}_0(\R^n)$
will be denoted by $\mathcal{PD}'_0(\R^n)$. See Appendices A.1--A.3 for more details.}
$m\in \mathcal{PD}_0(\R^n)$
definition (\ref{detorrente}) can be equivalently restated as the action of the linear and continuous functional
$\delta\big(f(\cdot;\widehat{\omega},\gamma)\big)\in \mathcal{PD}'_0(\R^n)$ on the test function $m$, i.e.,
\be\label{radondelta}
(Rm)(\widehat{\omega},\gamma) :=\int_{\R^n}\delta\big(f(x;\widehat{\omega},\gamma)\big) \,m(x)\, dx=
\int_{\R^n}\delta(\gamma-\widehat{\omega}\cdot x) \,m(x)\, dx.
\ee

Of course, definition (\ref{radon}) is more general than definition (\ref{radondelta}), since the former does not require $m$ to be
piecewise continuous and compactly supported; anyway, the two definitions coincide whenever $m\in \mathcal{PD}_0(\R^n)$.

Interestingly, definition (\ref{radondelta}) is naturally generalized from the case $\widehat{\omega}\in\mathbb{S}^{n-1}$
to the case $\omega\in\R^n\setminus\{0\}$:
let us discuss this point in detail.
For each $a\in\R\setminus\{0\}$ and $(\widehat{\omega},\gamma)\in\mathbb{S}^{n-1}\times\R$, we define the map
$f_{a}(\cdot;\widehat{\omega},\gamma):\R^n\rightarrow \R$ as $f_a(x;\widehat{\omega},\gamma):= a \gamma - a \widehat{\omega}\cdot x$.
Thus, the natural extension of definition (\ref{radondelta}) follows by setting, for all
$(\widehat{\omega},\gamma)\in \mathbb{S}^{n-1}\times\R$,
\be\label{radondeltalfa}
(Rm)(a\widehat{\omega},a\gamma) :=\int_{\R^n}\delta\big(f_a (x;\widehat{\omega},\gamma)\big) \,m(x)\, dx=
\int_{\R^n}\delta(a\gamma - a\widehat{\omega}\cdot x) \,m(x)\, dx.
\ee

Now, it is clear that
$|\mathrm{grad}\,f_{a}(x;\widehat{\omega},\gamma)|=|a|$ $\forall x\in\R^n$. Then,
by relations (\ref{defdelta}), (\ref{radon})  (\ref{radondelta}) and (\ref{radondeltalfa}),
we have, for any $m\in \mathcal{PD}_0(\R^n)$ and $a\in\R\setminus\{0\}$,
\begin{align}
  (Rm)(a\widehat{\omega}, a\gamma) &
  =|a|^{-1} \int_{\R^n}\delta(\gamma-\widehat{\omega}\cdot x)\,m(x)\,dx=\label{halle}\\
  & = |a|^{-1}\int_{\mathcal{P}(\widehat{\omega},\gamma)} m(x)\,d\sigma(x)=|a|^{-1}(Rm)(\widehat{\omega},\gamma).\nonumber
  \end{align}
Moreover, for any $\omega\in\R^n\setminus\{0\}$, let $\widehat{\omega}=\omega/|\omega|$ be the corresponding unit vector.
Thus, by (\ref{halle}),
for all $(\omega,\gamma)\in \left(\R^{n}\setminus\{0\}\right)\times\R$
we have
\be\label{martinoli1}
(Rm)(a\omega,a\gamma)=
|a\omega|^{-1}(Rm)\left(\widehat{\omega},\frac{\gamma}{|\omega|}\right)=|a|^{-1}(Rm)\left(\omega,\gamma \right).
\ee
Of course, by (\ref{martinoli1}), the Radon transform $(Rm)$ is uniquely determined by its values on $\mathbb{S}^{n-1}\times\R$.

Since the distance of the hyperplane
$\mathcal{P}(\omega,\gamma):=\{x\in\R^n : \gamma-\omega\cdot x=0\}$
from the origin $0\in\R^n$ is
$d\left(\mathcal{P}(\omega,\gamma),0\right)=|\gamma|/|\omega|$ and
the support of $m$ is compact, we have that
 \be\label{steffani}
 \forall \bar{\gamma}\in\R\setminus\{0\}\ \ \ \ \exists\,\lim_{(|\omega|,\gamma)\rightarrow (0,\bar{\gamma})}(Rm)(\omega,\gamma)=0.
 \ee
Then, property (\ref{martinoli1}) extends by continuity to $(\omega,\gamma)\in\R^{n+1}\setminus\{0\}$, i.e.,
\be \label{martinoli2}
(Rm)(a\omega, a \gamma)= |a|^{-1} (Rm)(\omega,\gamma)
\ \ \forall (\omega,\gamma)\in \R^{n+1}\setminus\{0\},\, \forall a\in\R\setminus\{0\},
\ee
being understood that
$(Rm)(0,\gamma):=0$ $\forall \gamma\in\R\setminus\{0\}$, as suggested by (\ref{steffani}).

As we are going to prove in a more general setting\footnote{See Theorem \ref{monteverdi} and Remark \ref{ziolineare1}
in Section \ref{ziogen}.},
$(Rm)$ is a locally integrable function on
$\left(\R^n\setminus\{0\}\right)\times\R$ and then, by (\ref{steffani}), onto $\R^{n+1}\setminus\{0\}$. Moreover,
relation (\ref{martinoli2}) shows that $(Rm)$ is an even homogeneous function
of $\omega$ and $\gamma$ of degree $-1$, which implies that
the singularity of $(Rm)$ at $(\omega,\gamma)=(0,0)\in\R^{n+1}$ is integrable, since $-1> -(n+1)$ for $n\geq 1$. Then,
recalling the inclusion map\footnote{See the end of \ref{test}, in particular definition (\ref{pairint}), as well as the end of
\ref{ziopiecewise}.}
$\widetilde{\iota}_k:L^1_{\mathrm{loc}}\left(\R^{n+1}\right)
\hookrightarrow \mathcal{PD}'_k\left(\R^{n+1}\right)$ for any $k\in\N$ or $k=\infty$,
we have
 \be \label{coero}
 (Rm)\in L^1_{\mathrm{loc}}\left(\R^{n+1}\right)\ \ \ \mbox{and}\ \ \ \widetilde{\iota}_k(Rm)\in\mathcal{PD}'_k\left(\R^{n+1}\right).
 \ee
Summing up, from now on we shall adopt the following definition of the Radon transform (cf. \cite[chap.~I]{gegr5}).
\begin{definition}\label{defradondelta1}
The \textnormal{Radon transform} of $\,m\in \mathcal{PD}_0(\R^n)$ is defined as the function
$(Rm):\R^{n+1}\setminus\{0\}\rightarrow\C$ given by
\be\label{radondelta1}
(Rm)(\omega,\gamma):=
\left\lbrace
\begin{array}{ll}
\int_{\R^n}\delta(\gamma-\omega\cdot x) \,m(x)\, dx & \forall (\omega,\gamma)\in
\left(\R^{n}\setminus\{0\}\right)\times\R;\\[2mm]
0 & \forall (\omega,\gamma)\in \{0\}\times\left(\R\setminus\{0\}\right).
\end{array}
\right.
\ee
\end{definition}

\subsection{Radon transform of the characteristic function of a square}\label{zioquadrato}

As an example and for future purpose, we now want to compute the Radon transform of the characteristic function of a square. This is
the key tool to solve the problem of
computing the Radon transform of any square-wise constant image, i.e.,
any plane image formed by square pixels, and then described by a function assuming, on each pixel, a constant
value (which may represent, e.g., the grey level\footnote{Cf. footnote no. \ref{bellissimi}.} associated with the pixel itself).
Indeed, by the linearity and translation properties of the Radon transform,
this problem is reduced to that of computing the Radon transform of
a single square pixel, with side of positive length $2a$ and centre at
the origin of the image plane. Thus, we are led to compute the Radon transform of the function $m:\R^2\rightarrow\R$ defined as
\be\label{emme}
m(x)=m(x_1,x_2):=
\left\lbrace
\begin{array}{ll}
1 & \mbox{if } (x_1,x_2)\in [-a,a]\times[-a,a],\\
0 & \mbox{otherwise.}
\end{array}
\right.
\ee

By setting $\omega=(\omega_1,1)\in\R^2\setminus\{0\}$ and $\mathsf{F}(x_1;\omega,\gamma)=-\omega_1 x_1 +\gamma$, with $\gamma\in\R$,
the equation of any straight line (not parallel to the $x_2$-axis)
in the image plane can be written as $x_2=\mathsf{F}(x_1;\omega,\gamma)$, i.e.,
\be\label{retta1}
\gamma - \omega_1 x_1 - x_2=0,
\ee
so that $f(x;\omega,\gamma)=\gamma - \omega_1 x_1 - x_2$.
Note that $|\mathrm{grad}_x\,f(x;\omega,\gamma)|^2=\omega_1^2+1=1+|\partial \mathsf{F}(x_1;\omega,\gamma)/\partial x_1|^2$, which implies the fulfilment
of property (\ref{viotti}). Then, by definitions (\ref{radondelta1}), (\ref{emme}) and (\ref{defdelta}), we have
\be\label{sdoura}
(Rm)(\omega,\gamma)
=\int_{-a}^{a} m\big(x_1,\mathsf{F}(x_1;\omega,\gamma)\big)\,dx_1.
\ee
Now, the integrand function in (\ref{sdoura}) does not vanish if and only if $\mathsf{F}(x_1;\omega,\gamma)\in [-a,a]$.
Accordingly, the integral in (\ref{sdoura}) coincides with the length of the interval obtained as the intersection of $[-a,a]$ with the
interval of variability for $x_1$ obtained from the condition $\mathsf{F}(x_1;\omega,\gamma)=x_2\in [-a,a]$, i.e.,
$-a\leq -\omega_1 x_1 +\gamma\leq a$. Depending on the possible values of
$a$, $\omega_1$, $\gamma$, the length of the intersection interval varies, as well as its analytical expression as a function of these
three parameters. However, it is also possible to obtain a single algebraic expression\footnote{Cf. \cite{we15}, with the identifications
$p=-\omega_1$, $\tau=\gamma$.}, given by
\be \label{radonquadrato}
(Rm)(\omega,\gamma)=\frac{|a-a\omega_1-\gamma|+|a-a\omega_1+\gamma|-|a+a\omega_1-\gamma|-|a+a\omega_1+\gamma|}{-2\omega_1}.
\ee

\section{The generalized Radon transform}\label{ziogen}

Taking inspiration from Definition \ref{defradondelta1} in the previous section,
it is natural to make a step further, i.e., to replace hyperplanes in $\R^n$
with $(n-1)$-dimensional submanifolds in an open subset of $\R^n$, parameterized by a finite number of parameters
$\lambda_1,\ldots,\lambda_t$ varying in an open subset of $\R^t$.

Often, when making such a generalization (see, e.g., \cite{be84, ku06}),
it is assumed that these submanifolds verify several specific conditions (e.g., smoothness, homogeneity, relationship between the
dimensions $n$ and $t$, definite positivity of the Hessian matrix), so that the corresponding generalized Radon
transform is endowed with structural properties preserving or resembling those of the classical Radon transform, in particular
its link with the Fourier transform. This approach is motivated by the need of investigating the most important issues of any integral
transform, i.e., 1) its injectivity (on an appropriate function space); 2) the characterization of its range; 3) inversion formulas and
corresponding algorithms; 4) the ill-posedness of the inverse problem (e.g., the stability of the reconstruction). While for the classical
Radon transform these problems have been solved (see, e.g., \cite{bebo98, helgason, na86, nawu01}),
only partial answers are known even for the spherical Radon transform \cite{ku06}, not to mention the case of more general submanifolds.

However, the focus of this paper is on the link between the Radon and the Hough transform and, to this end, only property (\ref{coero})
is of interest. Accordingly, in the following, we shall not be concerned with points 1)--4) above, thus being allowed to consider
submanifolds that are more general than those usually considered in the literature on this subject.

\begin{definition}\label{radongen}
For $n\in \N\setminus\{0,1\}$ and $t\in\N\setminus\{0\}$, let $W$ and $E$ be non-empty
open subsets of $\,\R^n$ and $\,\R^{t}$
respectively, and let\footnote{$E'$ is understood to be empty if and only if $t=1$.}
$E':=\{\lambda'\in \R^{t-1} : \exists \lambda_t\in\R : \lambda=(\lambda',\lambda_t)\in E\}$.
Moreover, let $f:W\times E \rightarrow\R$ be a function expressible in the $\lambda_t$-solvable form, i.e., as
$f(x;\lambda):=\lambda_t-F(x;\lambda')$, being $F:W\times E'\rightarrow\R$ such that $F\in C^1\left(W\times E'\right)$, and assume that
\begin{itemize}
\item[\textnormal{(i)}] $\displaystyle\mathcal{S}(\lambda):=\{x\in W : f(x;\lambda)=0\}\neq\emptyset$ $\forall\lambda\in E$;
\item[\textnormal{(ii)}] $\mathrm{grad}_x\, f(x;\lambda)\neq 0$ $\forall\lambda\in E$, $\forall x\in\mathcal{S}(\lambda)$.
\end{itemize}
Finally, let $m\in\mathcal{PD}_0(W)$. Then, the \textnormal{generalized Radon transform} of $m$ is defined as the function
$(R_f\, m):E\rightarrow\C$ given by\footnote{Note that the assumptions on $f$ allow defining the functional
$\delta\big(f(\cdot;\lambda)\big)$ for each $\lambda\in E$: see \ref{mementodelta} for details.}
\be\label{radongeneq}
(R_f\, m)(\lambda):=\int_{W} \delta\big(f(x;\lambda)\big)m(x)\,dx\ \ \ \forall \lambda\in E.
\ee
\end{definition}

We now want to prove the analogous of property (\ref{coero}) for the generalized Radon transform: this task is (step-wise)
accomplished by the following Theorem \ref{mozart}, Corollary \ref{haydn} and Theorem \ref{monteverdi}.

\begin{theorem}\label{mozart}
Notation and assumptions as in Definition \textnormal{\ref{radongen}},
except that $\mathcal{PD}_0(W)$ is to be replaced by $\mathcal{D}_0(W)$.
Thus, $(R_f\, m)\in C^0(E)$ for all $m\in\mathcal{D}_0(W)$.
\end{theorem}

\proof
Since $f(x;\lambda):=\lambda_t-F(x;\lambda')$ and $F\in C^1\left(W\times E'\right)$ by hypothesis, we have that
$f\in C^1\left(W\times E\right)$. Now, given $\widetilde{\lambda}\in E$, consider, according to
assumption (i) of Definition \ref{radongen}, the corresponding non-empty submanifold
$\mathcal{S}\big(\widetilde{\lambda}\big)$. By condition (ii) of Definition \ref{radongen},
for any point $\widetilde{x}\in\mathcal{S}\big(\widetilde{\lambda}\big)$
there exists $i\in\{1,\ldots,n\}$ such that
\be\label{battifollo}
\frac{\partial f}{\partial x_i}\big(\widetilde{x};\widetilde{\lambda}\big)\neq 0;
\ee
just for notational simplicity, assume that $i=n$.

Next, for positive $\epsilon_1$, $\epsilon_2$ and $\epsilon_3$, define
\begin{align}
W' & :=\left\{x'\in \R^{n-1} : \exists\,x_n\in\R : (x',x_n)\in W\right\},\\
W_n & :=\left\{x_n\in \R : \exists\,x'\in\R^{n-1} : (x',x_n)\in W\right\},\\
B\left(\widetilde{x}',\epsilon_1\right) & :=\left\{x'\in W' : \big|x'-\widetilde{x}'\big|<\epsilon_1 \right\},\\
B\big(\widetilde{\lambda},\epsilon_2\big) & := \big\{\lambda\in E :
\big|\lambda-\widetilde{\lambda}\big| < \epsilon_2 \big\},\\
B\big(\widetilde{x}_n,\epsilon_3\big) & :=\left\{x_n\in W_n : \big|x_n-\widetilde{x}_n\big| < \epsilon_3\right\}.
\end{align}
Then, by condition (\ref{battifollo}) and the implicit function theorem, for $k=1,2,3$ we can take
$\epsilon_k=\epsilon_k\big(\widetilde{x};\widetilde{\lambda}\big)$ so small that a function
\be
\mathsf{F}:B\big(\widetilde{x}',\epsilon_1\big(\widetilde{x};\widetilde{\lambda}\big)\big)\times
B\big(\widetilde{\lambda},\epsilon_2\big(\widetilde{x};\widetilde{\lambda}\big)\big)\rightarrow
B\big(\widetilde{x}_n,\epsilon_3\big(\widetilde{x};\widetilde{\lambda}\big)\big)
\ee
exists, satisfying the following properties:
\be\label{priola}
\mathsf{F}\in C^1\left(\overline{B\big(\widetilde{x}',\epsilon_1\big(\widetilde{x};\widetilde{\lambda}\big)\big)}\times
\overline{B\big(\widetilde{\lambda},\epsilon_2\big(\widetilde{x};\widetilde{\lambda}\big)\big)}\right)
\ee
and
\begin{align}\label{huet}
&\big\{\big(x',\lambda,\mathsf{F}(x';\lambda)\big):(x',\lambda) \in B\big(\widetilde{x}',\epsilon_1\big(\widetilde{x};\widetilde{\lambda}\big)\big)
\times B\big(\widetilde{\lambda},\epsilon_2\big(\widetilde{x};\widetilde{\lambda}\big)\big)\big\} = \\
&\big\{\big(x',\lambda, x_n\big)\in B\big(\widetilde{x}',\epsilon_1\big(\widetilde{x};\widetilde{\lambda}\big)\big)\times
B\big(\widetilde{\lambda},\epsilon_2\big(\widetilde{x};\widetilde{\lambda}\big)\big)\times
B\big(\widetilde{x}_n,\epsilon_3\big(\widetilde{x};\widetilde{\lambda}\big)\big) : f(x;\lambda)=0\big\}. \nonumber
\end{align}
Moreover, by condition (\ref{battifollo}) and the continuous differentiability of $f$, it is not restrictive to assume that
\be \label{bagnasco}
c\big[\epsilon_k\big(\widetilde{x};\widetilde{\lambda}\big)\big]:=\inf\left\{\left|\frac{\partial f}{\partial x_n}(x',x_n;\lambda)\right|
: (x',x_n;\lambda)\in
N\big[\epsilon_k\big(\widetilde{x};\widetilde{\lambda}\big)\big]\right\}>0,
\ee
where
\be\label{nucetto}
N\big[\epsilon_k\big(\widetilde{x};\widetilde{\lambda}\big)\big]:=
B\big(\widetilde{x}',\epsilon_1\big(\widetilde{x};\widetilde{\lambda}\big)\big)\times
B\big(\widetilde{x}_n,\epsilon_3\big(\widetilde{x};\widetilde{\lambda}\big)\big)\times
B\big(\widetilde{\lambda},\epsilon_2\big(\widetilde{x};\widetilde{\lambda}\big)\big).
\ee

For the same $\widetilde{\lambda}$, we can repeat the above construction for each $\widetilde{x}\in\mathcal{S}\big(\widetilde{\lambda}\big)$.
In particular, the union $\bigcup_{\widetilde{x}\in\mathcal{S}(\widetilde{\lambda})} B\big(\widetilde{x}',
\epsilon_1\big(\widetilde{x};\widetilde{\lambda}\big)\big)\times B\big(\widetilde{x}_n,\epsilon_3\big(\widetilde{x};
\widetilde{\lambda}\big)\big)$ is an open
covering of the closed subset $\mathcal{S}\big(\widetilde{\lambda}\big)$ of $\R^n$ and \textit{a fortiori} of the compact subset
$\mathcal{S}\big(\widetilde{\lambda}\big)\cap S_m$, being $S_m$ the compact support of $m\in\mathcal{D}_0(W)$. We can then extract a finite
subcovering of $\mathcal{S}\big(\widetilde{\lambda}\big)\cap S_m$, i.e., there exist a finite set of indices $r=1,\ldots,R$, with
$R=R\big(\widetilde{\lambda},m\big)\in\N$, and a corresponding finite subset $\big\{\widetilde{x}(r)\big\}_{r=1}^{R}$ of
$\mathcal{S}\big(\widetilde{\lambda}\big)\cap S_m$ such that
\be\label{putin}
\mathcal{S}\big(\widetilde{\lambda}\big)\cap S_m \subset \bigcup_{r=1}^R B\big(\widetilde{x}'(r),
\epsilon_1\big(\widetilde{x}(r);\widetilde{\lambda}\big)\big)\times B\big(\widetilde{x}_n(r),\epsilon_3\big(\widetilde{x}(r);\widetilde{\lambda}\big)\big).
\ee
Now, an equality analogous to (\ref{huet}) holds true for all $r=1,\ldots,R$: in each subset of the form (\ref{nucetto}), i.e., in each
$N\big[\epsilon_k\big(\widetilde{x}(r);\widetilde{\lambda}\big)\big]$,
the equation $f(x;\lambda)=0$ can be equivalently rewritten as $x_n=\mathsf{F}(x';\lambda)$. This implies that, by defining
$\epsilon_{2,m}\big(\widetilde{\lambda}\big):=\min\big\{\epsilon_2\big(\widetilde{x}(r);\widetilde{\lambda}\big):r=1,\ldots,R\big\}$ and taking
$\lambda\in B\big(\widetilde{\lambda},\epsilon_{2,m}\big(\widetilde{\lambda}\big)\big)$, the same covering on the right-hand side of
(\ref{putin}) also holds for $\mathcal{S}(\lambda)\cap S_m$.

Hence, for any such $\lambda$, this covering, together with any partition of unity $\{\rho_i\}_{i=1}^R$ subordinated to it,
can be used to compute the integral\footnote{Cf. relations (\ref{defdelta}) and (\ref{desigma}).}
\be\label{eltsin}
\int_W\delta\left(f(x;\lambda)\right)m(x)\,dx := \int_{\mathcal{S}(\lambda)}\frac{m(x)}{|\mathrm{grad}_x\,f(x;\lambda)|}\,d\sigma(x).
\ee
Indeed, by covering $\mathcal{S}(\lambda)\cap S_m$ as in (\ref{putin}), the integral on the right-hand side of (\ref{eltsin}) can be
computed as the finite sum of $R$ addenda: the generic $r$-th addendum is
\be \label{pievetta}
\int_{B\left(\widetilde{x}'(r),\epsilon_1\left(\widetilde{x}(r);\widetilde{\lambda}\right)\right)}\rho_r\left(x';\mathsf{F}(x';\lambda)\right)
m\left(x';\mathsf{F}(x';\lambda)\right)
\frac{\sqrt{1+\left|\mathrm{grad}_x\,\mathsf{F}(x';\lambda)\right|^2}}{|\mathrm{grad}_x\,f(x',\mathsf{F}(x';\lambda);\lambda)|}\,dx'.
\ee

We now note two properties concerning integral (\ref{pievetta}). First, since $\rho_r$ is infinitely differentiable,
$m$ is continuous and both $f$ and $\mathsf{F}$ are
continuously differentiable on their domain of definition,
it follows that the integrand function in (\ref{pievetta}) is continuous in $x'$ and $\lambda$. Second, taking into account conditions
(\ref{priola}), (\ref{bagnasco}) and the boundedness of $\rho_r$ and $m$, there exists a constant $K(r)\in\R^+$ such that
\be\label{guitton}
\rho_r\left(x';\mathsf{F}(x';\lambda)\right)
m\left(x';\mathsf{F}(x';\lambda)\right)
\frac{\sqrt{1+\left|\mathrm{grad}_x\,\mathsf{F}(x';\lambda)\right|^2}}{|\mathrm{grad}_x\,f(x',\mathsf{F}(x';\lambda);\lambda)|}
\leq K(r)
\ee
for all $(x',\lambda)\in B\big(\widetilde{x}'(r),\epsilon_1\big(\widetilde{x}(r);\widetilde{\lambda}\big)\big)\times
B\big(\widetilde{\lambda},\epsilon_{2,m}\big(\widetilde{\lambda}\big)\big)$. Thus, integral (\ref{pievetta}) converges uniformly with respect
to $\lambda$.

The two properties mentioned above imply that the function defined by
\be\label{maritain}
\lambda \mapsto
\int_{B\left(\widetilde{x}'(r),\epsilon_1\left(\widetilde{x}(r);\widetilde{\lambda}\right)\right)}
\rho_r\left(x';\mathsf{F}(x';\lambda)\right)
m\left(x';\mathsf{F}(x';\lambda)\right)
\frac{\sqrt{1+\left|\mathrm{grad}_x\, \mathsf{F}(x';\lambda)\right|^2}}{|\mathrm{grad}_x\,f(x',\mathsf{F}(x';\lambda);\lambda)|}\,dx'
\ee
is continuous on $B\big(\widetilde{\lambda},\epsilon_{2,m}\big(\widetilde{\lambda}\big)\big)$. Since this is true for each $r=1,\ldots,R$,
it immediately follows that also the function defined by
$\lambda\mapsto (R_f\,m)(\lambda):=\int_W\delta\left(f(x;\lambda)\right)m(x)\,dx$ is continuous
on the same domain, i.e., on a neighbourhood of $\widetilde{\lambda}$.

Finally, the same argument holds for any $\widetilde{\lambda}\in E$: hence,
$(R_f\, m)$ is continuous onto $E$ itself.
This concludes the proof.~$\square$

\vspace{3mm}
The following corollary is an immediate consequence of Theorem \ref{mozart}.
\begin{corollary}\label{haydn}
Assumptions as in Theorem \textnormal{\ref{mozart}}. Then,
$(R_f\, m)\in L^1_{\mathrm{loc}}(E)$ for all $m\in \mathcal{D}_0(W)$.
\end{corollary}

Actually, Corollary \ref{haydn} is a particular case of the following theorem, allowing for the case of piecewise continuous functions.

\begin{theorem}\label{monteverdi}
Assumptions as in Definition \textnormal{\ref{radongen}}. Then, $(R_f\, m)\in L^1_{\mathrm{loc}}(E)$ for all $m\in \mathcal{PD}_0(W)$.
\end{theorem}

\proof The proof is similar to that of Theorem \ref{mozart}. In fact, it is just the same up to expression (\ref{pievetta}). Then,
inequality (\ref{guitton}), which still holds true by conditions (\ref{priola}), (\ref{bagnasco}) and the boundedness of $\rho_r$ and $m$,
implies that the function defined as
\be
(x',\lambda)\mapsto \rho_r\left(x';\mathsf{F}(x';\lambda)\right)
m\left(x';\mathsf{F}(x';\lambda)\right)
\frac{\sqrt{1+\left|\mathrm{grad}_x\,\mathsf{F}(x';\lambda)\right|^2}}{|\mathrm{grad}_x\,f(x',\mathsf{F}(x';\lambda);\lambda)|}
\ee
is an element of $L^1\left(B\left(\widetilde{x}'(r),\epsilon_1\big(\widetilde{x}(r);\widetilde{\lambda}\big)\right)\times
B\big(\widetilde{\lambda},\epsilon_{2,m}\big(\widetilde{\lambda}\big)\big)\right)$. Thus, by Fubini theorem, the function defined by
(\ref{maritain}) is an element of $L^1\left(B\big(\widetilde{\lambda},\epsilon_{2,m}\big(\widetilde{\lambda}\big)\big)\right)$.
Since this is true for each $r=1,\ldots,R$,
it immediately follows that also the function defined by
$\lambda\mapsto (R_f\,m)(\lambda):=\int_W\delta\left(f(x;\lambda)\right)m(x)\,dx$ is Lebesgue-integrable
on the same domain, i.e., on a neighbourhood of $\widetilde{\lambda}$.

Finally, the same argument holds for any $\widetilde{\lambda}\in E$: hence,
$(R_f\, m)$ is an element of $L^1_{\mathrm{loc}}(E)$.
This concludes the proof.~$\square$

\begin{remark}\label{ziolineare1}
In the case of the Radon transform considered in Definition \ref{defradondelta1}, the assumptions of Theorem \ref{monteverdi}
(i.e., of Definition \ref{radongen}) are fulfilled for $t=n+1$, $W=\R^n$, $E=\left(\R^n\setminus\{0\}\right)\times\R$,
$\lambda'=\omega\in\R^n\setminus\{0\}$, $\lambda_t=\gamma\in\R$ and $f(x;\lambda)=\lambda_t-\lambda'\cdot x$.
Accordingly, we have $(Rm)\in L^1_{\mathrm{loc}}\big(\left(\R^n\setminus\{0\}\right)\times\R\big)$ by Theorem \ref{monteverdi},
then $(Rm)\in L^1_{\mathrm{loc}}\left(\R^{n+1}\setminus\{0\}\right)$ by setting $(Rm)(0,\gamma):=0$ $\forall\gamma\in\R\setminus\{0\}$,
as suggested by limit (\ref{steffani}). Finally, the fact that $(Rm)\in L^1_{\mathrm{loc}}\left(\R^{n+1}\right)$, i.e.,
property (\ref{coero}), follows
from the weak-singularity argument explained just below relation (\ref{martinoli2}).
\end{remark}

\subsection{A physical interpretation of the generalized Radon transform}

Taking inspiration from simple physical concepts, it is possible to establish an important result concerning the generalized Radon transform.
In order to accomplish this task, we need to prove a preliminary lemma.

\begin{lemma}\label{solvableL}
For $n\in\N\setminus\{0,1\}$ and $t\in\N\setminus\{0\}$, let $W$ and $E'$ be non-empty\footnote{Actually, $E'$ is understood to be empty
if and only if $t=1$.} open subsets of $\R^n$ and $\R^{t-1}$ respectively. Moreover, let $f:W\times\left(E'\times\R\right)\rightarrow\R$
be a function of the $\lambda_t$-solvable form $f(x;\lambda):=\lambda_t-F(x;\lambda')$
$\forall (x,\lambda)\in W\times\left(E'\times\R\right)$, with
$\lambda=(\lambda',\lambda_t)$. Finally, for each $\lambda\in E'\times\R$, let
\be\label{esselunga}
\mathcal{S}(\lambda):=\{x\in W : f(x;\lambda)=0\},\ \ \ \  \mathcal{S}^+(\lambda):=\{x\in W : f(x;\lambda)\geq 0\}.
\ee
Then, the following two properties hold true:
\begin{itemize}
\item[\textnormal{(i)}] if $\lambda_t\neq\widetilde{\lambda}_t$, then $\mathcal{S}(\lambda',\lambda_t)\cap
\mathcal{S}\big(\lambda',\widetilde{\lambda}_t\big)=\emptyset$ $\forall\lambda'\in E'$;
\item[\textnormal{(ii)}] if $\lambda_t\geq \widetilde{\lambda}_t$, then $\mathcal{S}^+(\lambda',\lambda_t)\supset
\mathcal{S}^+\big(\lambda',\widetilde{\lambda}_t\big)$ $\forall\lambda'\in E'$.
\end{itemize}
\end{lemma}
\proof (i)
According to the first of definitions (\ref{esselunga})
and in view of the specific form of $f(x;\lambda)=\lambda_t-F(x;\lambda')$, it holds that
\be\label{odio}
\mathcal{S}(\lambda',\lambda_t)=\{x\in W : \lambda_t=F(x;\lambda')\},\ \
\mathcal{S}\big(\lambda',\widetilde{\lambda}_t\big)=\{x\in W : \widetilde{\lambda}_t=F(x;\lambda')\}.
\ee
Then, $x\in\mathcal{S}(\lambda)\cap \mathcal{S}\big(\widetilde{\lambda}\big)$ implies that $\lambda_t=F(x;\lambda')=\widetilde{\lambda}_t$,
which contradicts the hypothesis $\lambda_t\neq\widetilde{\lambda}_t$.
It follows that $\mathcal{S}(\lambda)\cap \mathcal{S}\big(\widetilde{\lambda}\big)=\emptyset$.

(ii) By the second of definitions (\ref{esselunga}) and the specific form of $f$, we have that if $\lambda_t\geq \widetilde{\lambda}_t$, then
$\mathcal{S}^+(\lambda',\lambda_t)=\{x\in W : F(x;\lambda)\leq \lambda_t\}$ contains the set
$\{x\in W : F(x;\lambda)\leq \widetilde{\lambda}_t\}=\mathcal{S}^+\big(\lambda',\widetilde{\lambda}_t\big).$~$\square$

\vspace{3mm}
Now, if real-valued, a test function $m\in \mathcal{D}_0(W)$ can be regarded as the density with which some finite electric charge
(or mass, if $m$ is non-negative) is continuously distributed in free space.
Accordingly, in view of Lemma \ref{solvableL}, we shall denote by $M(\lambda',\lambda_t)$ the charge contained in the region
$\mathcal{S}^+(\lambda',\lambda_t)$, which, in general, becomes larger and larger as $\lambda_t$ increases.

Interestingly, by means of the generalized Radon transform of $m$, the following theorem establishes a link between the charge density
$m(x)$ and the charge $M(\lambda',\lambda_t)$ contained in $\mathcal{S}^+(\lambda',\lambda_t)$:
in this sense, the result can be considered as a physical interpretation of the generalized Radon transform itself. However, note that,
in the following, $m$ is not required to be real-valued.

\begin{theorem} \label{tristis}
Notation and hypotheses as in
Lemma \textnormal{\ref{solvableL}}. Moreover, assume the following properties:
\begin{itemize}
\item[\textnormal{(i)}] there exists a non-empty open subset $E$ of $E'\times\R$ such that $\mathcal{S}(\lambda)\neq\emptyset$
    $\forall\lambda=(\lambda',\lambda_t)\in E$;
\item[\textnormal{(ii)}] $f\in C^1(W\times E)$ and, for all $\lambda\in E$, $\mathrm{grad}_x\, f(x;\lambda)\neq 0$
$\forall x\in\mathcal{S}(\lambda)$.
\end{itemize}
Moreover, for any $m\in\mathcal{D}_0(W)$, let $M: E \rightarrow\C$ be defined as
\be\label{iperico}
 M(\lambda',\lambda_t):=\int_{\mathcal{S}^+(\lambda',\lambda_t)}m(x)\,dx=\int_{\{x\in W : F(x;\lambda') \leq \lambda_t\}} m(x)\,dx.
\ee
 Then, it holds that
 \be \label{cipster}
 \frac{\partial M}{\partial \lambda_t}(\lambda',\lambda_t)=(R_f\,m)(\lambda',\lambda_t)
 \ \ \ \forall (\lambda',\lambda_t)\in E,
 \ee
 where $(R_f\,m)$ is the generalized Radon transform of $m$, as defined in \textnormal{(\ref{radongeneq})}.
\end{theorem}

\proof
By definition,
\be \label{ipericina}
\frac{\partial M}{\partial \lambda_t}(\lambda',\lambda_t):=\lim_{h\rightarrow 0}\frac{M(\lambda',\lambda_t+h)-M(\lambda',\lambda_t)}{h}.
\ee
Assume that $h\rightarrow 0^+$ (the proof for the case $h\rightarrow 0^-$ is analogous). Then, from (\ref{iperico}) and (\ref{ipericina}),
we have
\be\label{ask}
\frac{\partial M}{\partial \lambda_t}(\lambda',\lambda_t)
=\lim_{h\rightarrow 0^+}\frac{1}{h}\int_{\left\{x\in\R^n : \lambda_t \leq
F(x;\lambda')\leq \lambda_t+h\right\}} m(x)\,dx.
\ee
By the representation of $\mathcal{S}(\lambda)$ as in (\ref{odio})
and the coarea formula\footnote{See (\ref{coarea}), with the identifications
$A=\left\{x\in\R^n : \lambda_t \leq F(x;\lambda')\leq \lambda_t+h\right\}$, $\Psi=F(\cdot;\lambda')$, $g=m$ and $s=\bar{\lambda}_t$.},
for each $(\lambda',\lambda_t)\in E$ we can rewrite (\ref{ask}) as
\be\label{ernia}
\frac{\partial M}{\partial\lambda_t}(\lambda',\lambda_t)=\lim_{h\rightarrow 0^+}\frac{1}{h}
\int_{\lambda_t}^{\lambda_t+h}\left(
\int_{\mathcal{S}(\lambda',\bar{\lambda}_t)}\frac{m(x)}{|\mathrm{grad}_x\, F(x;\lambda')|}\, d\sigma(x)\right) d\bar{\lambda}_t.
\ee
The internal integral in (\ref{ernia}), i.e., in view of (\ref{defdelta}),
\be\label{bay}
G\left(\lambda',\bar{\lambda}_t\right):=
\int_{\mathcal{S}(\lambda',\bar{\lambda}_t)}\frac{m(x)}{|\mathrm{grad}_x\, F(x;\lambda')|}\, d\sigma(x)=
\int_W\delta\big(f(x;\lambda',\bar{\lambda}_t)\big)m(x)\,dx,
\ee
is a continuous function on $E$ by Theorem \ref{mozart}. In particular, if we set
$E_t:=\{\lambda_t\in\R : \exists \lambda'\in\R^{t-1} : (\lambda',\lambda_t)\in E\}$, then $G(\lambda',\cdot)\in C^0(E_t)$ for all
$\lambda'\in E'$. Since it is not restrictive to assume $[\lambda_t,\lambda_t+h]\subset E_t$,
by the integral mean value theorem we find from (\ref{ernia})--(\ref{bay}) that
\be\label{discale}
\exists\, \widetilde{\lambda}_t(h)\in\left[\lambda_t, \lambda_t+h\right]\ :\
\frac{\partial M}{\partial\lambda_t}(\lambda',\lambda_t)=\lim_{h\rightarrow 0^+}\frac{1}{\xcancel{h}}\cdot \xcancel{h}\,\,
G\big(\lambda',\widetilde{\lambda}_t(h)\big).
\ee
Moreover, it holds that $\lim_{h\rightarrow 0^+}\widetilde{\lambda}_t(h)=\lambda_t$. Thus, by (\ref{bay})--(\ref{discale}) and the
continuity of $G(\lambda',\cdot)$, we have
\be \label{amburgo}
\frac{\partial M}{\partial \lambda_t}(\lambda',\lambda_t) = G\left(\lambda',\lambda_t\right)=
\int_W\delta\big(f(x;\lambda',\lambda_t)\big)m(x)\,dx.
\ee
Finally, by comparing (\ref{amburgo}) with (\ref{radongeneq}), assertion (\ref{cipster}) easily follows.~$\square$

\begin{remark}
By using the characteristic function $\Theta(f)$ introduced in (\ref{defHf1}), definition (\ref{iperico}) can be
equivalently rewritten in the form
\be\label{letta}
M(\lambda',\lambda_t):=\int_{W}\Theta\big(f(x;\lambda',\lambda_t)\big)\,m(x)\,dx.
\ee
Then, the result of Theorem \ref{tristis} can be heuristically obtained from (\ref{letta}) by formally interchanging the partial
derivative operator $\partial/\partial\lambda_t$ with the integral symbol and taking into account relations (\ref{derH1}) and
(\ref{radongeneq}), i.e.,
\begin{align}\label{lite}
\frac{\partial M}{\partial \lambda_t}(\lambda',\lambda_t) & =
\int_{W}\frac{\partial \Theta\big(\lambda_t-F(x;\lambda')\big)}{\partial \lambda_t}\,m(x)\,dx \\
& =\int_{W}\delta\big(\lambda_t-F(x;\lambda')\big)\,m(x)\,dx = (R_f\,m)(\lambda',\lambda_t). \nonumber
\end{align}
\end{remark}

\section{Generalized sinograms}\label{sinogrammi}

The classical Radon transform can be generalized so to act on distributions \cite{gegr5,helgason}. Under appropriate
assumptions, a corresponding theory could be developed for the generalized Radon transform.
However, here we are only interested in determining
the generalized Radon transform of the Dirac delta $\delta(\cdot -\widetilde{x})=\delta_{\widetilde{x}}\in\mathcal{D}'_0(W)$ centred at a point
$\widetilde{x}\in W\subset \R^n$. To this end, we shall adopt
an \textit{ad hoc} argument based on the property (\ref{cipster}) of the generalized Radon transform stated in Theorem \ref{tristis},
thus avoiding any approach concerned with distributions in general.

First, we observe that $\delta_{\widetilde{x}}$ can be regarded as an infinite charge density corresponding to a unit charge
concentrated at the point $\widetilde{x}\in W$. Then, we can imagine that such a density is the limit in $\mathcal{D}'_0(W)$
of a sequence of (feasible) charge densities $m_i$ as $i\rightarrow\infty$: this is made precise by the following Lemma \ref{safena}.
Finally, we can compute the Radon transform of
$\delta_{\widetilde{x}}$ as the limit in $\mathcal{D}'_1\left(E\right)$ of the generalized
Radon transforms $(R_f\, m_i)$ as $i\rightarrow\infty$: this is formalized by the subsequent Theorem \ref{daflon}.

\begin{lemma}\label{safena}
Let $W$, $\widetilde{x}\in W$ and $K\subset W$ be a non-empty open subset of $\R^n$,
a given point of $W$ and a compact subset of $W$ containing a neighbourhood $U_{\widetilde{x}}$ of $\widetilde{x}$, respectively.
Then, there exist (infinitely many) sequences of functions $\left\{m_i\right\}_{i\in\N}\subset\mathcal{D}_0(W)$ such
that\footnote{See (\ref{kant})--(\ref{pairint}) for details about the convergence in $\mathcal{D}'_k(W)$ and the inclusion map
$\iota_k:L^1_{\mathrm{loc}}(W)\rightarrow \mathcal{D}'_k(W)$.}
$\mathrm{supp}\,m_i\subset K$ $\forall i\in \N$ and $\iota_0(m_i)\rightarrow \delta_{\widetilde{x}}$ in $\mathcal{D}'_0(W)$ as $i\rightarrow\infty$, i.e.,
\be\label{reni}
\lim_{i\rightarrow\infty}\langle \iota_0(m_i),\phi\rangle=\langle\delta_{\widetilde{x}},\phi\rangle=\phi(\widetilde{x})
\ \ \forall\phi\in\mathcal{D}_0(W).
\ee
\end{lemma}

\proof
The proof can be obtained from an easy adaptation of standard results that can be found, e.g., in \cite[pp. 43--44]{zem87}.~$\square$

\begin{theorem}\label{daflon}
 Let $\{m_i\}_{i\in \N}\subset\mathcal{D}_0(W)$ be a sequence of functions as in Lemma \textnormal{\ref{safena}}, and let
 $f:W\times E\rightarrow\R$ be as in Theorem \textnormal{\ref{tristis}}.
 Then, it holds that
 $\iota_1\left(R_f\,m_i\right)\rightarrow \delta\big(f(\widetilde{x};\cdot)\big)$
 in $\mathcal{D}'_1\left(E\right)$ as $i\rightarrow\infty$, i.e.,
recalling \textnormal{(\ref{kant})--(\ref{pairint})} and \textnormal{(\ref{defdelta})},
 \be \label{diuretico}
 \lim_{i\rightarrow\infty}\left\langle \iota_1\left(R_f\, m_i\right), \psi\right\rangle=
 \left\langle\delta\big(f(\widetilde{x};\cdot)\big),\psi\right\rangle\ \ \forall\psi\in\mathcal{D}_1\left(E\right),
 \ee
 where $\delta\big(f(\widetilde{x};\cdot)\big)$ is the Dirac delta of the function $f(\widetilde{x};\cdot):E\rightarrow\R$.
 \end{theorem}

\proof
From Corollary \ref{haydn} and relations (\ref{pairint}), (\ref{derdistr}), (\ref{cipster}), (\ref{letta}), we have
\begin{align} \label{aiuto1}
  \langle \iota_1\left(R_f\,m_i\right),\psi \rangle & =\left\langle\iota_1\left(\frac{\partial M_i}{\partial \lambda_t}\right),
  \psi\right\rangle=
  -\left\langle\iota_0(M_i),\frac{\partial\psi}{\partial \lambda_t}\right\rangle\\
  & = - \int_{E}\left[\int_{W} m_i(x)\, \Theta\big(f(x;\lambda)\big)\,dx\right]
  \frac{\partial\psi}{\partial \lambda_t}(\lambda)\,d\lambda\nonumber\\
  & = - \int_{W}m_i(x)\left[\int_{E}\frac{\partial\psi}{\partial \lambda_t}(\lambda)\,
  \Theta\big(f(x;\lambda)\big)\,d\lambda\right] dx,\nonumber
\end{align}
where the last equality follows from Fubini theorem. Now, let us consider the internal integral, i.e., the function $\Phi:W\rightarrow\C$
defined by
\be\label{defPhi}
W \ni x\mapsto \Phi(x):= \int_{E}\frac{\partial\psi}{\partial \lambda_t}(\lambda)
\,  \Theta\big(f(x;\lambda)\big)\,d\lambda\in\C.
\ee
By (\ref{defHf1}) and the compactness of $\mathrm{supp}\,\psi$, the integration domain in (\ref{defPhi}) can be restricted to the
intersection $I(\psi,x)$ of $\mathrm{supp}\,\psi$ with the set
$\mathcal{S}^+(x):=\left\{\lambda\in E: f(x;\lambda)\geq 0\right\}$, i.e.,
\be\label{defPhihalf}
\Phi(x)= \int_{I(\psi,x)}\frac{\partial\psi}{\partial \lambda_t}(\lambda) \,d\lambda,\ \ \ \ \mbox{with}\
I(\psi,x):=\mathrm{supp}\,\psi\cap \mathcal{S}^+(x).
\ee
We note that $I(\psi,x)$ is compact for any $\psi\in\mathcal{D}_1\left(E\right)$ and any $x\in W$. Moreover, we can prove that
$\Phi$ is continuous on $W$, i.e.,
\be \label{pergolesi}
\lim_{x\rightarrow x^\ast}\left|\Phi(x)-\Phi(x^\ast)\right|=0\ \ \forall x^\ast\in W.
\ee
To prove limit (\ref{pergolesi}), we first set
$M_t:=\max_{\lambda\in E}\left|\frac{\partial\psi}{\partial \lambda_t}(\lambda)\right|$; then, we respectively denote by
$\Delta$ and $\mathcal{L}^{t}$ the symmetric difference between two sets and the Lebesgue measure on $\R^{t}$. Accordingly,
from (\ref{defPhihalf}) we have
\be\label{durante}
\left|\Phi(x)-\Phi(x^\ast)\right|=\left|\int_{I(\psi,x)\,\Delta\, I(\psi,x^\ast)}
\frac{\partial\psi}{\partial \lambda_t}(\lambda) \,d\lambda \right|\leq M_t \,
\mathcal{L}^{t}\big(I(\psi,x)\,\Delta\,I(\psi,x^\ast)\big).
\ee
Now, it is easy to realize that
\be\label{leo}
\lim_{x\rightarrow x^\ast}\mathcal{L}^{t}\big(I(\psi,x)\,\Delta\,I(\psi,x^\ast)\big)=0\ \ \ \forall x^\ast\in W,
\ \forall\psi\in\mathcal{D}_1\left(E\right).
\ee
Hence, limit (\ref{pergolesi}) readily follows from relations (\ref{durante})--(\ref{leo}).

In general, the function $\Phi$ is not compactly supported, but $m_i$ is, with $\mathrm{supp}\,m_i\subset K$ as in Lemma \ref{safena}.
Then, let $A$ be an open and bounded subset of $W$ such that $K\subset A$. By Urysohn lemma, there exists a continuous
function $u:W\rightarrow [0,1]$ such that $u(x)=1$ for $x\in K$ and $u(x)=0$ for $x\in W\setminus A$. It follows that the function
mapping
$x$ into $\widetilde{\Phi}(x):=u(x)\Phi(x)$ is both continuous and compactly supported, i.e., $\widetilde{\Phi}\in \mathcal{D}_0(W)$; moreover,
it clearly holds that $\big\{x\in W : \widetilde{\Phi}(x)=\Phi(x)\big\}\supset K\supset \mathrm{supp}\,m_i\cup U_{\widetilde{x}}$.
Then, recalling (\ref{pairint}) and (\ref{defPhi}),
we can rewrite the last equality in (\ref{aiuto1}) as
\be\label{vivaldi}
\left\langle \iota_1\left(R_f\, m_i\right),\psi\right\rangle=- \int_{W} m_i(x)\widetilde{\Phi}(x)\,dx=
- \left\langle\iota_0(m_i),\widetilde{\Phi}\right\rangle.
\ee
From (\ref{vivaldi}) and property (\ref{reni}), which is satisfied by assumption, we find
\be \label{strozzi}
\lim_{i\rightarrow\infty}\left\langle\iota_1\left(R_f\,m_i\right),\psi\right\rangle =
- \lim_{i\rightarrow\infty}\left\langle\iota_0(m_i),\widetilde{\Phi}\right\rangle= - \widetilde{\Phi}(\widetilde{x}).
\ee
Finally, by (\ref{pairint}), (\ref{derdistr}), (\ref{derH1}), (\ref{defPhi}) and since $u(\widetilde{x})=1$, we have
\begin{align} 
-\widetilde{\Phi}(\widetilde{x}) & = - u(\widetilde{x}) \int_{E}\frac{\partial\psi}{\partial \lambda_t}(\lambda)\,
  \Theta\big(f(\widetilde{x};\lambda)\big)\,d\lambda=
  -\left\langle \iota_0 \left[\Theta\big(f(\widetilde{x};\cdot)\big)\right], \frac{\partial\psi}{\partial\lambda_t}\right\rangle\nonumber \\
  & = \left\langle \frac{\partial\,\iota_0\left[\Theta\big(f(\widetilde{x};\cdot)\big)\right]}{\partial\lambda_t},\psi \right\rangle=
  \left\langle\delta\big(f(\widetilde{x};\cdot)\big),\psi \right\rangle.\label{albinoni}
  \end{align}
Then, relation (\ref{diuretico}) is obtained from an immediate comparison between (\ref{strozzi}) and (\ref{albinoni}).~$\square$

\vspace{3mm}
Summing up, from Theorem \ref{daflon} it follows that the appropriate definition of the generalized Radon transform of the Dirac delta
$\delta_{\widetilde{x}}\in\mathcal{D}'_{0}(W)$ is
\be\label{radeltagen}
\left(R_f\,\delta_{\widetilde{x}}\right)(\cdot):=\delta\big(f(\widetilde{x};\cdot)\big)\in \mathcal{D}_1'(E).
\ee

The specific form of (\ref{radeltagen}) for the Radon transform considered in Definition \ref{defradondelta1} deserves a short discussion.
In view of the usual identifications
made in Remark \ref{ziolineare1}, we can rewrite (\ref{radeltagen}) as
\be\label{radelta}
(R\delta_{\widetilde{x}})(\omega,\gamma):=\delta(\gamma - \omega\cdot \widetilde{x}),
\ee
where $\delta(\gamma - \omega\cdot \widetilde{x})$ is the Dirac delta $\delta(f(\widetilde{x};\cdot))\in\mathcal{D}_1'\left(\R^{n+1}\right)$
of the function\footnote{Analogously to Remark \ref{ziolineare1}, the hypotheses of Theorem \ref{tristis} are satisfied for
$E=\left(\R^n\setminus\{0\}\right)\times\R$, then its thesis (\ref{cipster}) does not hold true, in principle, for all $\lambda\in\R^{n+1}$.
However, $\mathcal{L}^{n+1}\left(\R^{n+1}\setminus E\right)=0$ and, as observed in Remark \ref{ziolineare1}, $(Rm_i)\in
L^1_{\mathrm{loc}}(\R^{n+1})$, so that the first equality in (\ref{aiuto1}) is valid
on $\R^{n+1}$ (i.e., $R_f\,m_i$ and $\partial M_i/\partial\lambda_t$ are equal
as elements of $L^1_{\mathrm{loc}}(\R^{n+1})$). As a result, in Theorem \ref{daflon} we can set $E=\R^{n+1}$ and then regard here
$\delta(f(\widetilde{x};\cdot))$ as an element of $\mathcal{D}_1'\left(\R^{n+1}\right)$.}
mapping $(\omega,\gamma)$ into $f(\widetilde{x};\omega,\gamma):=\gamma-\omega\cdot \widetilde{x}$.

\begin{remark}\label{def-sinogramma}
Definition (\ref{radelta}) is the mathematical justification of the name ``sinogram'' given to the two-dimensional representation
of the intensity values of the Radon transform of an
image in X-ray Computerized Tomography (CT).
Indeed, for $n=2$, we can model a single point (a dimensionless pixel) $P$
in the image as a Dirac
delta\footnote{We denote by $x(P)$ the $n$ coordinates of $P$, i.e., $x(P):=\left(x_1(P),\ldots,x_n(P)\right)\in\R^n$.}
$\delta_{x(P)}$, whose Radon transform is given by $\delta(\gamma - \omega\cdot x(P))$, according to (\ref{radelta}).
Now, the support of $\delta(\gamma - \omega\cdot x(P))$ in the parameter space is the plane 
$\mathcal{P}(x(P))=\{(\omega_1,\omega_2,\gamma)\in\R^3\,:\,\gamma - \omega\cdot x(P)=0\}$. By intersecting such plane with the cylinder $\mathbb{S}^1\times\R$,
which amounts to expressing $\omega=(\omega_1,\omega_2)\in\mathbb{S}^1$ as $\omega=(\cos\vartheta,\sin\vartheta)$ for $\vartheta\in [0,2\pi)$,
we find a set of points described by the equation $\gamma= x_1(P)\cos\vartheta + x_2(P)\sin\vartheta$,
which is a sinusoidal curve in the $(\vartheta,\gamma)$-plane.
\end{remark}

Moreover, an image of greater complexity or an
object $m$ can be modelled as a set of a
finite number $\nu$ of dimensionless pixels $P_1,\ldots,P_{\nu}$ having respective grey levels $\mu_1,\ldots,\mu_\nu$,
which correspond (in X-ray CT) to
the values of the linear attenuation coefficient of the object at those points.
This amounts to taking $m$ as
\be\label{vinci1}
m(x)=\sum_{j=1}^\nu \mu_j\,\delta\big(x-x(P_j)\big),\ \mbox{with}\
\mu_j\in \R , \  x(P_j)\in \R^2 \ \ \forall j=1,\ldots,\nu.
\ee
From (\ref{radelta}), (\ref{vinci1}) and the linearity of the Radon transform, we can then compute the Radon transform of $m(x)$ as
\be\label{trip1}
(Rm)(\gamma,\omega)=
\sum_{j=1}^{\nu} \mu_j\,\delta\big(\gamma-\omega\cdot x(P_j)\big),
\ee
thus obtaining, as its support, a superposition of $\nu$ sinusoidal curves in the $(\vartheta,\gamma)$-plane, which is just how a
sinogram appears. Of course, the same result (\ref{trip1}) also holds for $x(P_j)\in\R^n$, with $n>2$.

Even though no sinusoidal curve is involved in the general case, by analogy we shall call ``(generalized) sinogram'' any visual
representation of the intensity values of the generalized Radon transform of (\ref{vinci1}), which is, according to (\ref{radeltagen}),
\be\label{bilzerian}
\left(R_f\,m\right)(\lambda)=\sum_{j=1}^{\nu} \mu_j\,\delta\big(f(x(P_j);\lambda)\big).
\ee
By analogy, we shall speak of (generalized) sinogram also in the case of a piecewise continuous image, i.e., to indicate any
visual representation of the intensity values of the generalized Radon transform of $m\in\mathcal{PD}_0(W)$, as given by (\ref{radongeneq}).

\section{The Hough transform}\label{HT}

The Hough transform is a pattern recognition technique
for the automated detection of curves in images.
We refer, e.g., to  \cite{ba81, bemapi13, bero12, de81, duha72, ho62, le92, le93, macapebe15, prilki92, vaetal04}
for background material and complete details.
Here, we limit ourselves to recalling that the problem solved by this technique can be formulated in short as follows.
Given an image whose points are contained in an open subset $W$ of $\R^n$,
a set of points $\left\{P_j\right\}_{j=1}^{\nu}$ of
interest in the image itself
and a $\lambda$-parametrized family of functions $f_\lambda: W \rightarrow \R$,
find, among all possible values of the parameters $\lambda=(\lambda_1,\ldots,\lambda_t)\in E\subset\R^t$, the values
$\bar{\lambda}=(\bar{\lambda}_1,\ldots,\bar{\lambda}_t)$ for which
the corresponding zero locus of $f_{\bar{\lambda}}$, i.e., $\mathcal{S}(\bar{\lambda})=\{x\in W : f_{\bar{\lambda}}(x)=0\}$
(typically, a curve for $n=2$), best fits the set of points $\left\{P_j\right\}_{j=1}^{\nu}$.

In this section, we recall few basic definitions and discuss some concepts enlightening a new approach.

\subsection{The general setting}\label{SU}
The Cartesian product of copies of $\R$ (or $\A^1(\R)$, the affine space)
considered below is equipped with the Euclidean topology. First, let us fix some notation and preliminaries:

\begin{itemize}
\item[(i)] $x:=(x_1,\ldots,x_n)$, orthogonal Cartesian coordinates in the {\em image space} $\A^n(\R)$ (with $n\geq 2$),
also denoted by $\A_x^n(\R)$ and often identified with $\R^n$ itself;
\item[(ii)] $\lambda:=(\lambda_1,\ldots,\lambda_t)$, orthogonal Cartesian coordinates in the {\em parameter space} $\A^t(\R)$ (with $t\geq 1$),
also denoted by $\A_{\lambda}^t(\R)$ and often identified with $\R^t$ itself;
\item[(iii)] $W$ and $E$, non-empty open subsets of points in $\A^n(\R)$ and $\A^t(\R)$, respectively.
For notational simplicity and homogeneity with respect to the previous sections,
we shall also indicate by $x$ or $\lambda$ a point in $W$ or $E$, which amounts (by a slight abuse of language) to identifying
$W$ or $E$ with their coordinate representation in $\R^n$ or in $\R^t$, respectively;
\item[(iv)] $f:W\times E \to \R$, a function such that,
for each $\lambda\in E$, the map $f_\lambda:=f(\cdot;\lambda):W\rightarrow\R$ defined by $x\mapsto f(x;\lambda)$
satisfies the following conditions: (a) ${\mathcal S}(\lambda):=\{x\in W : f_{ \lambda}(x)= 0\}\neq \emptyset$; (b)
$f_{\lambda}\in C^1(W)$; (c) $(\mathrm{grad}_x\, f_{\lambda})(x)\neq 0$ $\forall x\in\mathcal{S}(\lambda)$.
\end{itemize}
Note that, as shortly explained between equalities (\ref{ulisse}) and (\ref{dini}) in \ref{mementodelta},
the previous conditions (a)--(c)
imply that ${\mathcal S}(\lambda)$ is a smooth, closed, orientable and $(n-1)$-dimensional
submanifold of $W\subset\R^n$, for each $\lambda\in E$.

We now propose the following definition of Hough transform.

\begin{definition}\label{def1}
Let $f:W\times E \to \R$ be a function satisfying conditions \textnormal{(a)--(c)} above, let $P\in W$ be a point in the image
space having coordinates $x=(x_1(P),\ldots, x_n(P))$, and let $f_x:E\rightarrow\R$ be the map
defined by $\lambda\mapsto f(x;\lambda)$. Then we say that the zero locus of $f_x$, defined as
${\mathcal H}(x):=\left\lbrace\lambda\in E : f_x(\lambda)=0\right\rbrace$, is the
{\em Hough transform of the point $P$ with respect to the function $f$ and to the coordinate system $x=(x_1,\ldots,x_n)$}.
If no confusion will arise, we simply say that $\mathcal{H}(x)$
is the {\em Hough transform of $x$}.
\end{definition}

Summarizing, the function $f:W\times E\rightarrow\R$ introduced above, when evaluated either at a fixed point $\lambda\in E$ of
the parameter space or at a fixed point $P \in W$ of the image space, defines, respectively,
\be\label{pesto}
{\mathcal{S}}(\lambda) =\left\lbrace x\in W : f_{\lambda}(x)=0\right\rbrace; \ \ \
\mathcal{H}(x) =\left\lbrace \lambda\in E : f_x(\lambda)=0\right\rbrace.
\ee

Clearly, for each $(x,\lambda)\in W\times E$, the {\em duality condition} (already understood in the algebraic case in \cite{bemapi13})
\begin{equation}\label{Dual}
x\in {\mathcal S}(\lambda) \Longleftrightarrow 0=f_{\lambda}(x)=f(x;\lambda)=f_x(\lambda)=0
\Longleftrightarrow\lambda \in {\mathcal H}(x)
\end{equation}
holds true, allowing us to conclude that {\em the Hough transform ${\mathcal H}(x)$ of a point $x\in W$ contains a point $\lambda\in E$
if and only if ${\mathcal S}(\lambda)$ passes through $x$.}

Note that, in general, the Hough transform operator $\mathcal{H}:W\rightarrow E$ mapping $x$ to $\mathcal{H}(x)$ is not injective, since
we may have ${\mathcal H}(x)={\mathcal H}(x')$ for different points $x,x'\in W$ (see example \ref{CS} at the end of this section).

An issue naturally arising from the previous setting is that of investigating the geometrical properties of the Hough transform
$\mathcal{H}(x)$.
Since $\mathcal{H}(x)$ is the zero locus of the function $f_{x}:E\rightarrow\R$, such properties will
depend on corresponding properties of $f_x$.
Here we limit ourselves to shortly recall the following.

A plain situation occurs if $f_x\in C^1(E)$ and $(\mathrm{grad}_\lambda\, f_x)(\lambda)\neq 0$ for each $\lambda\in\mathcal{H}(x)\neq\emptyset$:
in this case, $\mathcal{H}(x)$ is a smooth, closed, orientable and $(t-1)$-dimensional submanifold of $E\subset\R^t$
(thus paralleling the
properties of $\mathcal{S}(\lambda)$ in $W\subset\R^n$; cf. item (iv) above).
However, the condition $(\mathrm{grad}_\lambda\, f_x)(\lambda)\neq 0$ may hold only for some $(x,\lambda)\in W\times E$: in this case,
${\mathcal H}(x)$ is a $(t-1)$-dimensional submanifold locally around $\lambda$.

In general, ${\mathcal H}(x)$ may be empty, or may contain (or even may be made up of) irreducible components of dimension not greater than
$t-2$ (for instance, a single point). On the other hand, $\mathcal{H}(x)$ may be equal to the whole $E\subset\R^t$: this happens if $x$ is a
\textit{base point} of the family $\{{\mathcal S}(\lambda)\}_{\lambda\in E}$, i.e., a point belonging to ${\mathcal S}(\lambda)$ for all
$\lambda\in E$. Indeed, in this case the duality condition (\ref{Dual}) implies that ${\mathcal H}(x)=E$ (see example \ref{CS} again).

\subsection{The weighted Hough counter}\label{WHC}

We briefly describe here, with slight modifications with respect to \cite{bemapi13}, the basic steps of the algorithm leading to
the construction of the weighted Hough counter, which is a key tool of the pattern recognition technique based on the Hough transform,
as implemented in \cite{bemapi13}.

First, let $f:W\times E\rightarrow\R$ be a function satisfying conditions (a)--(c) stated in item (iv) of Subsection \ref{SU}.
Then, consider the following steps.
\begin{enumerate}
\item[I.] \textit{Discretization of the parameter space.} Identify a suitable (and bounded) investigation domain $\mathcal{T}\subset E$
in the parameter space $\R^t$. Next, choose an \textit{initialization point} $\lambda^\ast=(\lambda_1^\ast,\ldots,\lambda_t^\ast)$ in
$\mathcal{T}$ and, for each $k=1,\ldots,t$, a \textit{sampling distance} $d_k$ with respect to the component $\lambda_k$. Then, set
\begin{equation}\label{IN}
\lambda_{k,\mathcal{n}_k}:=\lambda_k^{*}\pm \mathcal{n}_k d_k, \;\; k=1,\ldots,t, \;\; \mathcal{n}_k=0,\ldots,\mathcal{N}_k-1,
\end{equation}
where $\mathcal{N}_k$ is half the number of considered samples for such component, and $\mathcal{n}_k$ the index labelling the sample.
Moreover, denote by
\begin{equation}\label{cells}
\hspace{-1cm}
C(\mathcal{n}):=\left\{\lambda=(\lambda_1,\ldots,\lambda_t) \in \mathcal{T} : \lambda_k\in
\left[\lambda_{k,\mathcal{n}_k}-\frac{d_k}{2},\,\lambda_{k,\mathcal{n}_k}+\frac{d_k}{2}\right)\ \forall k=1,\ldots,t\right\}
\end{equation}
the rectangular cell with centre in the \textit{sampling point}
$\lambda_{\mathcal{n}}:=(\lambda_{1,\mathcal{n}_1},\ldots,\lambda_{t, \mathcal{n}_t})$ of the discretized region
$\mathcal{T}$, where $\mathcal{n}\in\N^t$ denotes the multi-index $(\mathcal{n}_1,\ldots, \mathcal{n}_t)$ labelling the cells.
Finally, denote by $C(\lambda)$ the cell containing the point $\lambda\in \mathcal{T}$.
We point out that the discretization is defined by relation (\ref{IN}), that is, by the choice of the initialization point
$\lambda^\ast\in\mathcal{T}$ and the \textit{discretization step} represented by the multi-index $d:=(d_1,\ldots,d_t)\in\R_{+}^t$.
In the following, a discretization will be also denoted by $\{\lambda^\ast,d\}$.

\item[II.] \textit{Definition of the Hough transform kernel.} Let $P$ be a point in a subset $W\subset \R^n$ of the image space,
with coordinates $x=(x_1(P),\ldots,x_n(P))$, and let $\mathcal{H}(x)$ be the Hough transform of $P$ with respect to the function
$f:W\times E\rightarrow\R$ and the coordinate system $x=(x_1,\ldots,x_n)$ (cf. Definition \ref{def1}). Then consider the following map,
depending on both the function $f$ and the discretization (\ref{IN}), defined by
\begin{equation}\label{K1}
p(x, \lambda;\lambda^\ast,d):=
\begin{cases}
1 \; \;\; \mbox{if} \;\; \mathcal{H}(x)\cap C(\lambda)\neq \emptyset , \\
0\; \;\; \mbox{otherwise}.
\end{cases}
\end{equation}
For a given discretization $\{\lambda^\ast,d\}$, the map
\be
p(\cdot,\cdot\, ;\lambda^\ast,d) : W\times \mathcal{T} \to \{0,1\}
\ee
is also called \textit{Hough transform kernel (with respect to the function $f$)}.

From a numerical viewpoint, the problem of computing $p(x,\lambda;\lambda^\ast,d)$, i.e., establishing whether the Hough transform
$\mathcal{H}(x)$ intersects a cell $C(\mathcal{n})$ or not, is not so easy as it might appear at first sight: see, e.g., \cite{tobe14}
for a discussion of this problem in the algebraic case. Here we shall not deal with such an issue, since we are going to make an assumption
on the analytic form of $f$ (i.e., $\lambda_k$-solvability), whereby $p(x,\lambda;\lambda^\ast,d)$ can be properly
redefined and easily computed (see Subsection \ref{solF}).

\item [III.] \textit{Introduction of weights and construction of the Hough accumulator.}
For any given set of points of interest in the image space, say $P_{j}$,
$j=1,\ldots,\nu$, denote by $\mu_j$ the grey level\footnote{Cf. footnote no. \ref{bellissimi}.}
associated with $P_{j}$. For all $j=1,\ldots,\nu$,
let $x(P_{j})\in\R^n$ denote the $n$ coordinates of the point $P_{j}$. Accordingly, the mathematical description of this set of points,
regarded as a discrete image, can be given in terms of a linear combination of Dirac deltas centred at $x(P_{j})$, whose respective
coefficients are the weights $\mu_j$, just as in equality (\ref{vinci1}).
We then introduce the \textit{weighted Hough counter (with respect to the function $f$)}, also called \textit{weighted Hough accumulator},
as the map $H(\cdot;\lambda^\ast,d): \mathcal{T} \to \N$ defined by\footnote{\label{zioambiguo}Some authors define the Hough transform itself as
the right-hand side of (\ref{WHT}): see, e.g., \cite{prilki92}, eq. (7).}
\begin{equation}\label{WHT}
H(\lambda;\lambda^\ast,d):=\sum_{j=1}^\nu \mu_j\, p\big(x(P_j), \lambda;\lambda^\ast,d\big).
\end{equation}
\end{enumerate}

For sake of completeness, we just recall that the set of points $\{P_j\}_{j=1}^{\nu}$ can be often selected by processing the image
through an appropriate edge-detection algorithm. Then, the grey level of all the points $P_j$ is usually set to $1$ and the values
$\bar{\lambda}=(\bar{\lambda}_1,\ldots,\bar{\lambda}_t)\in\mathcal{T}$ for which the manifold
$\mathcal{S}(\bar{\lambda})$ best fits the points $P_j$ can be found as those maximizing the Hough counter, i.e., those representing
the centre of the rectangular cell with the maximum number of intersections with all the Hough transforms of the points $P_j$.
In this framework, the so-called Hough regularity, i.e., the property whereby $\mathcal{S}(\lambda)=\mathcal{S}(\lambda')$
implies $\lambda=\lambda'$, plays an important role. However, here we shall not deal with any pattern recognition technique;
we again refer to \cite{bemapi13, macapebe15}
for details, examples and discussion of some numerical issues.

\subsection{The Hough transform kernel and $\lambda_t$-solvability}\label{solF}

Motivated from the framework of Section \ref{ziogen} (see, in particular, Definition \ref{radongen})
and in order to avoid possible pathologies,
like those highlighted in the final part of Subsection \ref{SU} above,
we now focus on a specific form of the function $f:W\times E\rightarrow\R$ and on some properties following from it.

To this end, set $\lambda:=(\lambda',\lambda_k)$, with $k\in\{1,\ldots,t\}$, and define
$E':=\left\{\lambda'\in \R^{t-1} : \exists \lambda_k\in \R : \lambda:=(\lambda',\lambda_k)\in E\right\}$,
being $E'=\emptyset$ if and only if $t=1$. Up to renaming the variables $\lambda_i$'s, we can always assume that $k=t$.
Consider now a function $f:W\times E \rightarrow \R$ satisfying conditions (a)--(c) stated in item (iv) of Subsection \ref{SU}.
Furthermore, assume the following two conditions to be true:
(d) $f$ is $\lambda_t$-\textit{solvable}, i.e., as introduced in Definition \ref{radongen}, $f(x;\lambda)$ is of the form
$f(x;\lambda_1,\ldots,\lambda_t)= \lambda_t-F(x;\lambda_1,\ldots,\lambda_{t-1})$; (e) for each $x\in W$,
the function $f_x:E\rightarrow\R$ introduced in Definition \ref{def1} is continuously differentiable, i.e., in view of (d),
$F_x\in C^1(E')$, where, for each $x\in W$, the map $F_x:E'\rightarrow\R$ is obviously defined as
$\lambda'\mapsto F(x;\lambda')$. Thus, we can prove the following lemma.

\begin{lemma}\label{HTgood}
Let $f:W\times E\rightarrow\R$ be a function satisfying conditions \textnormal{(d)--(e)} above, and let $P\in W$ be a point of coordinates
$x=(x_1(P),\ldots,x_n(P))$. Thus, if non-empty, the Hough transform $\mathcal{H}(x)=\left\lbrace \lambda\in E : f_x(\lambda)=0\right\rbrace$
of $x$ is a smooth, closed, orientable and $(t-1)$-dimensional submanifold of $E\subset\R^t$.
\end{lemma}
\proof 
By conditions (d)--(e), the gradient of $f$ with respect to $\lambda$ can be computed as
$(\mathrm{grad}_{\lambda}\,f_x)(\lambda)=\big(\frac{\partial F_x}{\partial\lambda_1}(\lambda),
\ldots,\frac{\partial F_x}{\partial\lambda_{t-1}}(\lambda),1\big)$, thus showing that it never vanishes on $E\subset\R^t$.
Then, the same remark just below item (iv) in Subsection \ref{SU} suffices to conclude the proof.~$\square$

\begin{figure}[!htbp]
\begin{center}
\includegraphics[scale=0.49]{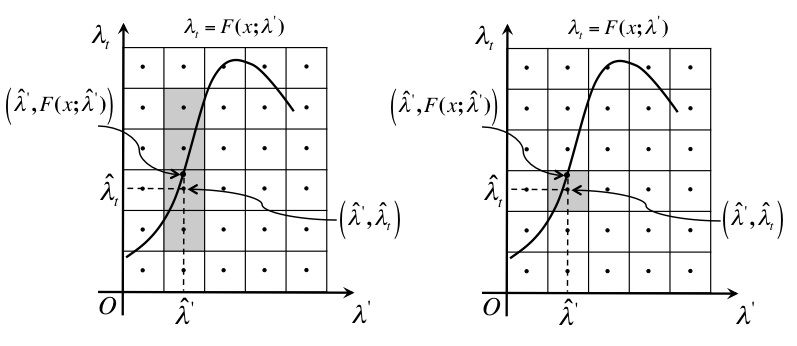}\\
(a)\ \ \ \ \ \ \ \ \ \ \ \ \ \ \ \ \ \ \ \ \ \ \ \ \ \ \ \ \ \ \ \ \ \ \ \ \ \ \ \ \ \ \ \ (b)
\caption{(a) The Hough transform $\mathcal{H}(x)$, which is the graph of the function $\lambda_t=F(x;\lambda')$, intersects four cells
in the same column parallel to the $\lambda_t$-axis; the first $t-1$ coordinates of the centres of the cells coincide, and are equal to
$\widehat{\lambda}'$. (b) The $\lambda_t$-solvability of
$\mathcal{H}(x)$ allows choosing a single cell among the four previous ones, according to a minimum distance criterion.}
\label{manualiter}
\end{center}
\end{figure}

\vspace{3mm}
As anticipated in Subsection \ref{WHC}, the $\lambda_t$-solvability of $f$ inspires an appropriate redefinition of the
Hough transform kernel $p(x,\lambda;\lambda^\ast,d)$ and an easy way to compute it. To address this issue, we first observe that
if $f$ is $\lambda_t$-solvable, then the Hough transform $\mathcal{H}(x)$ can be regarded as the graph of a function $F(x;\cdot)$ of
$\lambda'\in E'$, i.e.,
$\mathcal{H}(x)=\left\lbrace(\lambda',\lambda_t)\in E : \lambda_t=F(x;\lambda')\right\rbrace$. On the other hand, as shown in Figure
\ref{manualiter}(a),
it might happen that $\mathcal{H}(x)$ intersects two or more cells whose centres only differ by the $\lambda_t$-coordinate, i.e., cells belonging
to a column parallel to the $\lambda_t$-axis in the parameter space. Clearly, this circumstance is only due to the discretization of
the investigation domain $\mathcal{T}$ (i.e., to the fact that the length of the cells along the $\lambda'$-axes is positive),
while $\lambda_t$-solvability would rather suggest that for a certain discretized value of $\lambda'$, at most one cell
should be intersected by $\mathcal{H}(x)$, as a graph of a function of $\lambda'$. Interestingly, this drawback can be easily overcome
without refining or changing the kind of the discretization $\left\lbrace\lambda^\ast,d\right\rbrace$ in $\mathcal{T}$.
To this end, it suffices to choose,
among all the crossed cells in the same column, the one whose centre $(\widehat{\lambda}',\widehat{\lambda}_t)$ has the minimum distance
from the point $\big(\widehat{\lambda}',F(x;\widehat{\lambda}')\big)\in\mathcal{H}(x)$ (see Figure \ref{manualiter}(b)).
This amounts to redefining the Hough transform kernel, since we need to drop all the other crossed cells in the same column.

\FloatBarrier

The new definition can be made explicit, from a computational viewpoint, as follows. As explained in item I of Subsection \ref{WHC},
the sampling points $\lambda_{\mathcal{n}}=(\lambda_{1,\mathcal{n}_1},\ldots,\lambda_{t,\mathcal{n}_{t}})$ are the centres of the cells $C(\mathcal{n})$ covering the investigation domain $\mathcal{T}$.
Clearly, for each $k=1,\ldots,t$, the set of all the $k$-th coordinates $\lambda_{k,\mathcal{n}_k}$ of the sampling points induces a
corresponding discretization (with step $d_k$) of the values of
the $k$-th continuous variable $\lambda_k$.
In particular, if we set
$\mathcal{T}_k:=\left\lbrace\lambda_k\in\R : \exists (\lambda_1,\ldots,\xcancel{\lambda_k},\ldots,\lambda_t)\in\R^{t-1}
: (\lambda_1,\ldots,\lambda_k,\ldots,\lambda_t)\in \mathcal{T}\right\rbrace$,
we can define the function $c_k:\mathcal{T}_k\rightarrow\R$ mapping $\lambda_k\in\mathcal{T}_k$ to
its closest discretized value. The ambiguity arising when $\lambda_k$ is just half-way between two discretized values is removed by taking
the larger value, in agreement with (\ref{cells}).
Accordingly, in view of (\ref{IN}), (\ref{cells}), the analytic expression of the function $c_k$ is
\be\label{mappack}
c_k(\lambda_k)=\lambda^\ast_k+\left\lfloor \frac{1}{2}+\frac{\lambda_k-\lambda^\ast_k}{d_k}\right\rfloor d_k,
\ee
where $\lfloor\cdot\rfloor$ denotes, as usual, the floor function, mapping $x\in\R$ to the largest integer not greater than $x$.
Moreover, by using (\ref{mappack}) and setting $\mathcal{T}':=\left\lbrace\lambda'\in\R^{t-1} : \exists
\lambda_t\in\R : (\lambda',\lambda_t)\in \mathcal{T}\right\rbrace$,
we can also define the map $c': \mathcal{T}'\rightarrow \R^{t-1}$ as
$\lambda'=(\lambda_1,\ldots,\lambda_{t-1})\mapsto c'(\lambda')= \left(c_1(\lambda_1),\ldots,c_{t-1}(\lambda_{t-1})\right)$.
Then, for fixed $x\in W$ and $\lambda'\in \mathcal{T}'$, among all the cells intersected by $\mathcal{H}(x)$ and having centres in
$\big(c'(\lambda'),\lambda_{t,\mathcal{n}_t}\big)$, we want to select the one whose centre has coordinates
$\big(c'(\lambda'), c_t(F(x;c'(\lambda')))\big)$. This amounts to replacing the Hough transform kernel (\ref{K1}) by the following one:
\be\label{nuovokernel}
p(x,\lambda;\lambda^\ast,d):=
\begin{cases}
1 \; \;\; \mbox{if} \;\; (x,\lambda)\in C(x;\lambda), \\
0\; \;\; \mbox{otherwise,}
\end{cases}
\ee
where
$C(x;\lambda):=\left\{(x,\lambda)\in W\times \mathcal{T} : -d_t/2 \leq \lambda_t - F(x;c'(\lambda')) <d_t/2 \right\}$.
From now on throughout the paper, we shall always use definition (\ref{nuovokernel}) for the Hough transform kernel, both in the theoretical
discussion and in the numerical computation of the weighted Hough counter.

Finally, we give a definition that will prove useful in the following sections, where the link between the Hough transform
and the Radon transform will be investigated.

\begin{definition}\label{Hsinog}
Let $f:W\times E\rightarrow\R$ be a function satisfying conditions \textnormal{(a)--(e)} above, and let $\mathcal{T}\subset E$ be
the investigation domain introduced in item \textnormal{I} of Subsection \textnormal{\ref{WHC}}.
Then, the function mapping $\lambda\in\mathcal{T}$ to
$H(\lambda;\lambda^\ast,d)/d_t$ (i.e., the ratio between the weighted Hough counter defined in \textnormal{(\ref{WHT})} and the
sampling distance $d_t$ with respect to the component $\lambda_t$) will be called \textnormal{rescaled (weighted) Hough counter}.
A visual representation of the intensity values of $H(\lambda;\lambda^\ast,d)/d_t$ in the coordinate system $(\lambda_1,\ldots,\lambda_t)$
will be called \textnormal{Hough sinogram}.
\end{definition}

\subsection{The algebraic case}\label{AC}

An important class of functions satisfying the conditions imposed above is those of polynomials.
Indeed, the Hough transform is a standard pattern recognition technique initially introduced for the detection of straight lines,
circles and ellipses. Some foundational results, based on algebraic geometry arguments, strongly support an extension of this method
to the automated recognition of special plane algebraic curves in images, to detect profiles of interest of various shapes.
We refer to \cite{bemapi13, macapebe15, bero12} for complete details, examples, and further developments.

Let $\alpha:=(\alpha_1,\ldots,\alpha_n)$ be the multi-index characterizing monomials
$x_1^{\alpha_1}\cdots x_n^{\alpha_n}$ of degree $|\alpha |:=\sum_{i=1}^n\alpha_i$.
Then, consider a $\lambda$-parametrized family of irreducible polynomials in the variable $x$,
of a given degree $d$ independent of $ \lambda$, that is,
\be\label{polinomi}
f_{\lambda}(x):=\sum_{|\alpha|=0}^dx_1^{\alpha_1}\cdots x_n^{\alpha_n}g_\alpha(\lambda)\in \R[x_1,\ldots,x_n],
\ee
where $g_\alpha(\lambda)$ is a polynomial expression in $\lambda_1,\ldots,\lambda_t$,
such that for each $\lambda\in E\subset\R^t$ there exists $\alpha=\alpha(\lambda)$ with $|\alpha|=d$
and $g_\alpha(\lambda)\neq 0$.

Of course, whenever we take $W$ as an open subset of $\R^n$ that is disjoint from the set of singular points of
$\mathcal{S}(\lambda)$,
the polynomial $f_{\lambda}(x)$ satisfies conditions (b) and (c) stated in item (iv) of Subsection \ref{SU} above,
so that $\mathcal{S}(\lambda)$ is $(n-1)$-dimensional (see also \cite[Theorem 4.5.1]{bocoro98}).
Moreover, the hypersurfaces $\mathcal{S}(\lambda)$ defined as in (\ref{pesto}) are irreducible in $W\subset\R^n$,
since the polynomials (\ref{polinomi}) are assumed to be irreducible.
We then have a family $\{\mathcal{S}(\lambda)\}_{\lambda\in E}$ of smooth and irreducible hypersurfaces in $W$ having the same degree.
Clearly, over an algebraically closed field $K$, the irreducible polynomial $f_\lambda(x)\in K[x_1,\ldots,x_n]$ always defines
an irreducible hypersurface in $K^n$.

For each point $P$ of coordinates $x=(x_1(P),\ldots, x_n(P))$ in the image space $\R^n$, the Hough transform $\mathcal{H}(x)$
of the point $P$, if non-empty, is a hypersurface in the parameter space $\R^t$, defined by the polynomial equation
\be
f_x(\lambda) =
\sum_{|\alpha|=0}^d x_1(P)^{\alpha_1} \cdots  x_n(P)^{\alpha_n} g_{\alpha}( \lambda)=0,
\ee
provided that $(\mathrm{grad}_\lambda\, f_x)(\lambda)\neq 0$ for some $\lambda\in{\mathcal H}(x)$
(see \cite[Theorem 4.5.1]{bocoro98} again).

\begin{remark}
We observe that, in the algebraic case, it would be more natural to take the whole $\R^n$ as the domain of definition of the polynomials
$f_\lambda(x)$ appearing in (\ref{polinomi}).
On the one hand, this would allow considering possible singularities of $\mathcal{S}(\lambda)$, which
indeed characterize the geometry of the algebraic set itself. On the other hand, the presence of singularities, as well
as possible non-pure dimensionality issues (see example \ref{CS} again, and \cite[Section 1]{ribema15}), do not match with the typical
requirements
of smoothness which manifolds are assumed to satisfy in order to develop the classical integration theory on them. In turn, this theory
is crucial for defining and investigating important properties of the generalized Radon transform. This is the reason why we limit
ourselves to considering smooth manifolds $\mathcal{S}(\lambda)$, as done from the very beginning in Subsection \ref{SU}.
\end{remark}

Let us conclude this section with some examples.

\begin{example}\label{HIP}
\normalfont
\textit{Hyperplanes.} As in Section \ref{RT}, consider the parameters $\lambda'=\omega=(\omega_1,\ldots,\omega_n)$ and $\lambda_t=\gamma$,
as well as the corresponding family of hyperplanes in the image space $\A_x^n(\R)$, defined by
${\mathcal P}(\omega,\gamma)=\left\{x\in \R^n : \gamma-\omega\cdot x=0 \right\}$
and having distance  $|\gamma|/|\omega|$ from the origin. The function $f(x;\gamma,\omega)=\gamma-\omega\cdot x$ is clearly $\gamma$-solvable.
The Hough transforms are hyperplanes in the parameter space $\A_{\lambda}^{n+1}(\R)$ of coordinates $\lambda=(\omega_1,\ldots,\omega_n,\gamma)$.
\end{example}

\begin{example}\label{EC}
\normalfont
\textit{Elliptic curves.} In the image plane $\A_{x}^2(\R)$, consider the family of cubic curves
expressed in the canonical Weierstrass form as
$x_2^2=x_1^3 + a x_1 + b$. With respect to the notation adopted for the general setting, here we have the identifications
$n=t=2$, $\lambda=(\lambda_1,\lambda_2)=(a,b)$, $\mathcal{S}(\lambda)=
\mathcal{C}(a,b)=\left\lbrace (x_1,x_2)\in\R^2 :  x_2^2 - x_1^3 - a x_1 - b =0 \right\rbrace$.
Non-singular curves from this family are elliptic curves. The function $f(x_1,x_2;a,b)=x_2^2 - x_1^3 - a x_1 - b$ is clearly
$b$-solvable. The Hough transforms are straight lines in the parameter plane $\A_{\lambda}^2(\R)$ of coordinates $\lambda=(a,b)$.
Slight variants of this family of curves has been successfully used to detect profiles of interest in both astronomical and
medical images (see \cite{bemapi13, macapebe15}).
\end{example}

\begin{example} \label{CS}
\normalfont
{\em Conchoid of Sl\"use.}
In the image plane ${\mathbb A}_x^2(\R)$, consider the family
of rational cubic curves defined by the equation
$a(x_1-a)(x_1^2+x_2^2) = b^2x_1^2$. With respect to the general notation, the identifications are now
$n=t=2$, $\lambda=(\lambda_1,\lambda_2)=(a,b)\in E=\R^+ \times \R^+$, $\mathcal{S}(\lambda)=\mathcal{C}(a,b)=\{(x_1,x_2)\in\R^2 :
a(x_1-a)(x_1^2+x_2^2) = b^2x_1^2 \}$. Note that the function $f(x_1,x_2;a,b)= a(x_1-a)(x_1^2+x_2^2)-b^2x_1^2$
is neither $a$- nor $b$-solvable.
Such a cubic is classically known as {\it conchoid of Sl\"use} of parameters $a$, $b$.
This curve has a double nodal point at the origin $O$, with complex conjugate tangent lines of equation $a^2(x_1^2+x_2^2)+b^2x_1^2=0$,
so that $O$ is an isolated point of the curve.
For any point $P$ of coordinates $x=\left(x_1(P), x_2(P)\right)$ in the image plane ${\mathbb A}_{x}^2(\R)$,
the Hough transform $\mathcal{H}(x)$ is an ellipse of equation
$\left[x_1(P)^2+x_2(P)^2\right]a^2 + x_1(P)^2\, b^2 - x_1(P)\left[x_1(P)^2+x_2(P)^2\right]a=0$ in the parameter plane
${\mathbb A}_{\lambda}^2(\R)$.
Clearly, the Hough transform of $O$ is the whole affine plane ${\mathbb A}_{\lambda}^2(\mathbb R)$. Moreover,
$\mathcal{H}(x)=\mathcal{H}(x')$ whenever $x'=\left(\pm x_1(P), \pm x_2(P)\right)$, which provides a simple example of the non-injectivity
of the Hough transform operator $\mathcal{H}$.
\end{example}

\section{Link between the Radon transform and the Hough transform: the case of discrete images}\label{ziodiscreto}

This section is devoted to proving the following Theorem \ref{equidiscreto}. Roughly speaking, its statement can be summarized as follows:
given a discrete image, i.e., an image formed by a finite number of pixels $P_1,\ldots,P_{\nu}$
having respective grey levels $\mu_1,\ldots,\mu_{\nu}$ (cf. relation (\ref{vinci1})), the corresponding rescaled Hough counter
tends to become the generalized Radon transform of the image itself (cf. equality (\ref{bilzerian}))
as the discretization of the parameter space becomes finer and finer. The precise statement is as follows.

\begin{theorem}\label{equidiscreto}
Let $f:W\times E\rightarrow\R$ be a function satisfying properties \textnormal{(a)--(e)} stated in the previous section,
i.e., \textnormal{(a)} ${\mathcal S}(\lambda):=\{x\in W : f_{ \lambda}(x)= 0\}\neq \emptyset$ $\forall\lambda\in E$;
\textnormal{(b)} $f_{\lambda}\in C^1(W)$ $\forall\lambda\in E$; \textnormal{(c)}
$(\mathrm{grad}_x\, f_{\lambda})(x)\neq 0$ $\forall x\in {\mathcal S}(\lambda)$, $\forall\lambda\in E$; \textnormal{(d)}
$f(x;\lambda_1,\ldots,\lambda_t)= \lambda_t-F(x;\lambda_1,\ldots,\lambda_{t-1})$; \textnormal{(e)} $F_x\in C^1(E')$ $\forall x\in W$.
Moreover, let $\{\lambda^\ast,d\}$ be a discretization of the parameter space, and define
$D:=\max\{d_1,\ldots,d_t\}$, where $d_k$, for $k=1,\ldots,t$, is the sampling distance with respect to the component $\lambda_k$,
as explained in item \textnormal{I} of Subsection \textnormal{\ref{WHC}}.
Finally, let $m$ be a discrete image, $(R_f\, m)(\lambda)$ its generalized Radon transform
and $H(\lambda;\lambda^\ast,d)/d_t$ the corresponding rescaled Hough
counter\footnote{Cf. equalities (\ref{vinci1}), (\ref{bilzerian}) and Definition \ref{Hsinog}, respectively.}, defined
on a bounded and open investigation domain $\mathcal{T}\subset E$. Then
\be\label{minacce}
\lim_{D\rightarrow 0^+}\,\iota_1\left(\frac{H(\lambda;\lambda^\ast,d)}{d_t}\right)=
(R_f\, m)(\lambda)
\ \ \mbox{in}\ \ \mathcal{D}'_1\left(\mathcal{T}\right),
\ee
where $\iota_1:L^1_{\mathrm{loc}}(\mathcal{T})\rightarrow \mathcal{D}'_{1}(\mathcal{T})$
denotes the inclusion map defined as in \textnormal{(\ref{pairint})}.
\end{theorem}

This theorem is an immediate consequence of the following technical Lemma \ref{haendel} and Corollary \ref{gasquet},
together with the subsequent identifications (\ref{rava})--(\ref{fava}).

\begin{lemma}\label{haendel}
Let $\Xi$ be a subset of $\,\R$ such that $\bar{\xi}$ is an accumulation
point for $\Xi$, and let $\xi\in\Xi$ be a parameter.
Moreover, for $t\in\N\setminus\{0,1\}$, let $E'$ be a non-empty open subset of $\,\R^{t-1}$
and, for each $\xi\in\Xi$, let
$U_\xi,V_\xi:E'\rightarrow\R$ be two functions of the variable $\lambda'\in E'$, endowed with the following properties:
\begin{itemize}
 \item[\textnormal{(i)}] both of them are
 elements of the space\footnote{Cf. Definition \ref{pezzi}
 in \ref{ziopiecewise}.} $PC^1(E')$;
 \item[\textnormal{(ii)}] $\exists\, \epsilon_\xi >0$ such that
 $V_\xi(\lambda') - U_\xi(\lambda') > \epsilon_\xi$ $\forall\lambda'\in E'$, $\forall\xi\in\Xi$;
 \item[\textnormal{(iii)}] $\forall\lambda'\in E'$
 $\displaystyle\exists \lim_{\xi\rightarrow\bar{\xi}}U_\xi(\lambda')=
 \lim_{\xi\rightarrow\bar{\xi}}V_\xi(\lambda')=:G(\lambda')\in\R$, with $G\in C^1(E')$;
 \item[\textnormal{(iv)}] the functions $u_\xi:=U_\xi-G$ and
$v_\xi:=V_\xi-G$ are uniformly bounded with respect to the parameter
$\xi$, i.e., there
exists a constant $M\geq 0$ such that
$\,\left|u_\xi(\lambda')\right|\leq M$ and $\,\left|v_\xi(\lambda')\right|\leq M$
$\forall\lambda'\in E'$, $\forall\xi\in\Xi$.
\end{itemize}
Finally, for each $\xi\in\Xi$ let us define:
\begin{itemize}
 \item[\textnormal{(a)}] the set
$C_{\xi}:=\left\{\lambda=(\lambda',\lambda_t)\in E'\times\R : U_\xi(\lambda')
\leq\lambda_t < V_\xi(\lambda')\right\}$;
\item[\textnormal{(b)}] the characteristic function of
$C_{\xi}$, i.e., $\mathbf{1}_{C_{\xi}}:E'\times\R\rightarrow\{0,1\}$;
\item[\textnormal{(c)}] the function $r_\xi:=v_\xi-u_\xi$; 
\item[\textnormal{(d)}] the function $T_\xi:=
\mathbf{1}_{C_{\xi}}/r_\xi\in L^1_{\mathrm{loc}}\left(E'\times\R\right)$ and the corresponding distribution
$\iota_1\left(T_\xi\right)\in\mathcal{D}'_1\left(E'\times\R\right)$;
\item[\textnormal{(e)}] the function defined by $E'\times\R\ni (\lambda',\lambda_t)\mapsto g(\lambda):=\lambda_t-G(\lambda')\in\R$ and the corresponding
Dirac delta $\delta(g)\in\mathcal{D}'_1\left(E'\times\R\right)$.
\end{itemize}
Then, it holds that
$\iota_1\left(T_\xi\right) \rightarrow \delta\left(g\right)$ in
$\mathcal{D}'_1(E'\times\R)$ as $\xi\rightarrow\bar{\xi}$.
\end{lemma}

\proof
According to (\ref{kant}), the thesis of the theorem can be recast as
\be \label{brunetta}
\lim_{\xi\rightarrow \bar{\xi}}\left\langle \iota_1\left(T_\xi\right),\psi\right\rangle=
\left\langle \delta\left(g\right),
\psi\right\rangle\ \ \ \forall \psi\in\mathcal{D}_1\left(E'\times\R\right).
\ee
Note that, by points (iii) and (e), we have $g\in C^1\left(E'\times\R\right)$, with $\partial g/\partial\lambda_t= 1$ identically. Thus,
$\mathrm{grad}\,g(\lambda)\neq 0$ $\forall \lambda\in E'\times\R$ and, in particular, $\delta(g)\in \mathcal{D}'_1(E'\times\R)$
is well-defined, according to definition (\ref{defdelta}).

In order to prove (\ref{brunetta}), we begin by recalling (\ref{pairint}) and the definition of $T_\xi$ in point (d) above,
so that, for all $\psi\in\mathcal{D}_1\left(E'\times\R\right)$ and $\xi\in\Xi$, we have
\be\label{edda}
\langle \iota_1\left(T_\xi\right),\psi\rangle  =
\int_{E'\times\R}\frac{\mathbf{1}_{C_{\xi}}(\lambda)}{r_\xi(\lambda')}\,\psi(\lambda)\,d\lambda
= \int_{C_\xi}\frac{\psi(\lambda)}{r_\xi(\lambda')}\,d\lambda.
\ee
By setting
\begin{align}
A_{\xi} & :=\{(\lambda',\lambda_t)\in E'\times\R :\lambda_t\geq U_\xi(\lambda')\},\label{A}\\
B_{\xi} & :=\{(\lambda',\lambda_t)\in E'\times\R :\lambda_t > V_\xi(\lambda')\},\label{B}
\end{align}
we easily realize that $C_{\xi}=A_{\xi}\Delta B_{\xi}$, where $\Delta$ denotes the symmetric
difference between two sets: accordingly, we have
\be\label{delpo}
\int_{C_\xi}\frac{\psi(\lambda)}{r_\xi(\lambda')}\,d\lambda=
\int_{A_\xi}\frac{\psi(\lambda)}{r_\xi(\lambda')}\,d\lambda-
\int_{B_\xi}\frac{\psi(\lambda)}{r_\xi(\lambda')}\,d\lambda.
\ee
Now, remembering the definition of $u_\xi$ and $v_\xi$ given in assumption (iv),
from (\ref{A})--(\ref{B}) we immediately get
\begin{align}
A_{\xi} &
=\{(\lambda',\lambda_t)\in E'\times\R : \lambda_t - u_\xi(\lambda')- G(\lambda') \geq 0\},\label{Au}\\
B_{\xi} &
=\{(\lambda',\lambda_t)\in E'\times\R : \lambda_t - v_\xi(\lambda')- G(\lambda') > 0\}.\label{Av}
\end{align}
Then, by assumption (i), we can make the following two $\xi$-dependent changes of coordinates
$\lambda=(\lambda',\lambda_t)\mapsto \eta^\xi=\left(\eta', \eta_t^\xi\right)$ almost everywhere
on $A_{\xi}$ and $B_{\xi}$ respectively:
\begin{align}
\eta'=\lambda',\ \ \eta_t^\xi=\lambda_t - u_\xi(\lambda')\ \ \ & \mbox{for}\ (\lambda',\lambda_t)
\in A_{\xi},\label{AJ}\\
\eta'=\lambda',\ \ \eta_t^\xi=\lambda_t - v_\xi(\lambda')\ \ \ & \mbox{for}\ (\lambda',\lambda_t)
\in B_{\xi}.\label{BJ}
\end{align}
An immediate check shows that both the Jacobian matrices of transformations (\ref{AJ}) and (\ref{BJ}) are triangular with $1$ on the
diagonal, so that their determinant is $1$. Moreover, relations (\ref{Au})--(\ref{BJ}) show that, when expressed in the new coordinates
$\eta^\xi$, the integration domains $A_{\xi}$ and $B_{\xi}$ become
\begin{align}
A &
:=\left\{\left(\eta',\eta_t^\xi\right)\in E'\times\R : \eta_t^{\xi} - G(\eta') \geq 0\right\},\label{Aue}\\
B &
:=\left\{\left(\eta',\eta_t^\xi\right)\in E'\times\R : \eta_t^{\xi} - G(\eta') > 0\right\},\label{Ave}
\end{align}
i.e., $A$ and $B$ are independent of $\xi$ and
coincide up to a zero-measure subset of $E'\times\R$.
Accordingly, from (\ref{edda}), (\ref{delpo}) and (\ref{Au})-(\ref{Ave}), we easily find
\be\label{delpoeta}
\langle\iota_1\left(T_\xi\right),\psi\rangle  =
\int_{A}\frac{\psi\big(\eta',\eta_t^\xi+u_\xi(\eta')\big)-\psi\big(\eta',\eta_t^\xi+v_\xi(\eta')\big)}
{r_\xi(\eta')}\,d\eta'd\eta_t^\xi.
\ee

For notational simplicity, we now change the names of the integration variables by setting\footnote{Incidentally, $\lambda'=\eta'$ is also
the first set of equations in the coordinate transformations (\ref{AJ})--(\ref{BJ}), but this has nothing to do with the current
renaming of the integration variables.}
$\lambda':=\eta'$ and $\lambda_t:=\eta_t^\xi$, so that (\ref{delpoeta}) becomes
\be\label{delpoetal}
\langle\iota_1\left(T_\xi\right),\psi\rangle  =
\int_{A}\frac{\psi\left(\lambda',\lambda_t + u_\xi(\lambda')\right)-\psi\left(\lambda',\lambda_t + v_\xi(\lambda')\right)}
{r_\xi(\lambda')}\,d\lambda.
\ee
As far as the integrand function in (\ref{delpoetal}) is concerned, we remember that, by definition (c),
$r_\xi(\lambda')=v_\xi(\lambda')-u_\xi(\lambda')>0$. Thus, by applying Lagrange mean value theorem, there exists
$\widetilde{\lambda}_t^{\xi}\in [\lambda_t +u_\xi(\lambda'),\lambda_t +v_\xi(\lambda')]$ such that
\be \label{Lagrange}
\frac{\psi\left(\lambda',\lambda_t + u_\xi(\lambda')\right)-\psi\left(\lambda',\lambda_t + v_\xi(\lambda')\right)}
{v_\xi(\lambda')-u_\xi(\lambda')}=
-\frac{\partial\psi}{\partial \lambda_t}\big(\lambda',\widetilde{\lambda}_t^{\xi}\big).
\ee

Moreover, by assumptions (iii) and (iv), we have
$\lim_{\xi\rightarrow \bar{\xi}} u_\xi(\eta')=0$ and
$\lim_{\xi\rightarrow \bar{\xi}}v_\xi(\eta')=0$:
as a consequence, $\lim_{\xi\rightarrow \bar{\xi}}\widetilde{\lambda}_t^{\xi}=\lambda_t$. Then, by (\ref{Lagrange}) and the
continuity of $\partial\psi/\partial\lambda_t$, we find
\be \label{Lagrangelim}
\lim_{\xi\rightarrow \bar{\xi}}\frac{\psi\left(\lambda',\lambda_t + u_\xi(\lambda')\right)-
\psi\left(\lambda',\lambda_t + v_\xi(\lambda')\right)}{v_\xi(\lambda')-u_\xi(\lambda')}=
-\frac{\partial\psi}{\partial\lambda_t}(\lambda',\lambda_t).
\ee
Furthermore, remembering assumption (iv) and the fact that $\psi\in C^{1}\left(\R^t\right)$, we deduce the inequality
\be \label{domina}
\left|\frac{\psi\left(\lambda',\lambda_t + u_\xi(\lambda')\right)-
\psi\left(\lambda',\lambda_t + v_\xi(\lambda')\right)}{v_\xi(\lambda')-u_\xi(\lambda')}\right|\leq
\|\psi\|_{C^1}\,\mathbf{1}_{K_M}(\lambda)\ \ \ \forall\xi\in\Xi,
\ee
where $\mathbf{1}_{K_M}(\cdot)$ denotes the characteristic function of the compact subset of $\R^t$ defined as
$K_M:=\{(\lambda',\lambda_t+\bar{\lambda}_t)\in\R^t : (\lambda',\lambda_t)\in\mathrm{supp}\,\psi\ \mbox{and}\ |\bar{\lambda}_t|\leq M\}$.
Since $\mathbf{1}_{K_M}(\cdot)\in L^1(\R^t)$, by (\ref{Lagrangelim})--(\ref{domina}) we can apply Lebesgue dominated convergence theorem in
(\ref{delpoetal}), thus obtaining
\be\label{gilbert}
\lim_{\xi\rightarrow \bar{\xi}}\langle\iota_1\left(T_\xi\right),\psi\rangle  =
-\int_{A}\frac{\partial\psi}{\partial\lambda_t}(\lambda)\,d\lambda.
\ee

Finally, having defined $g(\lambda):=\lambda_t-G(\lambda')$ in (e),
by (\ref{Aue}) and (\ref{defHf1}) the characteristic function of
$A$ can be written as $\R^t\ni\lambda\mapsto \Theta\left(g(\lambda)\right)$. Hence,
from (\ref{pairint}), (\ref{derdistr}) and (\ref{derH1}), we have
\begin{align}
\int_{A}\frac{\partial\psi}{\partial\lambda_t}(\lambda)\,d\lambda & =
\int_{\R^t}\Theta\left(g(\lambda)\right)\frac{\partial\psi}{\partial\lambda_t}(\lambda)\,d\lambda=
\left\langle\iota_0\big(\Theta(g)\big),\frac{\partial\psi}{\partial\lambda_t}\right\rangle\nonumber\\
&=-\left\langle\frac{\partial\, \iota_0\big(\Theta(g)\big)}{\partial \lambda_t},\psi\right\rangle
=-\left\langle\delta(g),\psi\right\rangle.\label{brad}
\end{align}
Now, an immediate comparison between (\ref{gilbert}) and (\ref{brad}) proves equality (\ref{brunetta}), as wanted.~$\square$

\begin{corollary}\label{gasquet}
 For each $j\in\{1,\ldots,J\}$, assume that $U_\xi^j$, $V_\xi^j$, $G^j$, $g^j$, $C_\xi^j$, $T_\xi^j$
 verify the hypotheses satisfied, respectively, by $U_\xi$, $V_\xi$, $G$, $g$, $C_\xi$, $T_\xi$ in
 Lemma~\textnormal{\ref{haendel}}.
 Moreover, let $\beta_j\in\C$ for all $j\in\{1,\ldots,J\}$. Then
 \be\label{marcello}
 \lim_{\xi\rightarrow\bar{\xi}}\,\iota_1\left(\sum_{j=1}^J \beta_j\, T_\xi^j\right) = \sum_{j=1}^{J} \beta_j\,
 \delta(g^j)\ \ \mbox{in}\ \ \mathcal{D}'_1\left(E'\times\R\right).
 \ee
 \end{corollary}
\proof
From Lemma \ref{haendel}, for each $j=1,\ldots,J$, we have that $\iota_1\left(T_\xi^j\right)\rightarrow \delta(g^j)$ in
$\mathcal{D}'_1\left(E'\times\R\right)$ as $\xi\rightarrow \bar{\xi}$. From (\ref{hume}) and (\ref{kant}), it follows that
$\sum_{j=1}^J \beta_j\,\iota_1\left(T_\xi^j\right)\rightarrow\sum_{j=1}^{J} \beta_j\, \delta(g^j)$ in $\mathcal{D}'_1\left(E'\times\R\right)$
as $\xi\rightarrow \bar{\xi}$. Then, limit (\ref{marcello}) immediately follows from the linearity of the inclusion map
$\iota_1:L^1_{\mathrm{loc}}(E'\times\R)\hookrightarrow\mathcal{D}'_1\left(E'\times\R\right)$.~$\square$

\vspace{3mm}
Let us now see how Corollary \ref{gasquet} applies to the rescaled Hough counter. To this end,
we first observe that, by Definition \ref{hegel} in \ref{test},
convergence in $\mathcal{D}'_1\left(E'\times\R\right)$ implies
convergence in $\mathcal{D}'_1\left(\mathcal{T}\right)$,
for any open and bounded investigation domain $\mathcal{T}\subset E\subset E'\times\R$.
Furthermore, keeping into account the notation and
definitions introduced throughout Subsections \ref{WHC}, \ref{solF} in the case of a discrete image
(i.e., an image described as in (\ref{vinci1}),
with $\R^2$ replaced by $\R^n$), we can choose or identify the functions and parameters appearing in the
statement of Corollary \ref{gasquet} as follows:
\begin{align}
& \xi=D=\max\{d_1,\ldots,d_t\},\ \ \bar{\xi}=0,\label{rava}\\
& J=\nu,\ \ \beta_j=\mu_j,\ \ r_{\xi}^j=d_t\ \ \forall j=1,\ldots,\nu,\\
& U_{\xi}^j(\lambda')=-d_t/2 + F\big(x(P_j);c'(\lambda')\big),\ \ V_{\xi}^j(\lambda')=d_t/2 + F\big(x(P_j);c'(\lambda')\big),\\
& G^j(\lambda')=F\big(x(P_j); \lambda'\big),\ \ g^j(\lambda)=f\big(x(P_j);\lambda\big)=\lambda_t-F\big(x(P_j);\lambda'\big),\\
& C_\xi^j=\left\{\lambda=(\lambda',\lambda_t)\in E : -d_t/2 \leq \lambda_t - F\big(x(P_j);c'(\lambda')\big) <d_t/2 \right\},\\
& \mathbf{1}_{C_\xi^j}(\lambda)=p\big(x(P_j),\lambda;\lambda^\ast,d \big),\ \
T_\xi^j(\lambda)=\frac{\mathbf{1}_{C_\xi^j}(\lambda)}{r_\xi^j}=\frac{p\big(x(P_j),\lambda;\lambda^\ast,d\big)}{d_t},\\
& \sum_{j=1}^J \beta_j\, T_\xi^j(\lambda)=\sum_{j=1}^\nu \mu_j\, \frac{p\big(x(P_j),\lambda;\lambda^\ast,d\big)}{d_t}=
\frac{H(\lambda;\lambda^\ast,d)}{d_t}.
\label{fava}
\end{align}
An easy check shows that identifications (\ref{rava})--(\ref{fava}) ensure the fulfilment of the hypotheses required by Corollary
\ref{gasquet}, so that the corresponding form of statement (\ref{marcello}) is now
\be\label{minacce1}
\lim_{D\rightarrow 0^+} \,\iota_1\left(\frac{H(\lambda;\lambda^\ast,d)}{d_t}\right)=
\sum_{j=1}^\nu\mu_j\,\delta\big(\lambda_t-F\left(x(P_j);\lambda'\right)\big)
\ \ \mbox{in}\ \ \mathcal{D}'_1\left(\mathcal{T}\right).
\ee
Finally, it suffices to recall that
the right-hand side of equality (\ref{minacce1}) is just the generalized Radon transform of a discrete image formed by
$\nu$ points $P_1,\ldots,P_\nu$, with corresponding grey levels $\mu_1,\ldots,\mu_\nu$, as shown in relations (\ref{vinci1}) and
(\ref{bilzerian}). Accordingly, relation (\ref{minacce1}) can be equivalently rewritten as (\ref{minacce}), thus proving Theorem
\ref{equidiscreto} and justifying the claims opening this section.

\subsection{The one-dimensional case $t=1$}
For sake of completeness, let us now see how the previous investigation trivializes when $t=1$.
The one-dimensional counterpart of Lemma \ref{haendel} can be formulated as follows.

\begin{lemma}\label{zarelis}
Let $\Xi$ be a subset of $\,\R$ such that $\bar{\xi}$ is an accumulation point for $\Xi$, and let $\xi\in\Xi$ be a parameter.
For each $\xi\in\Xi$, let $U_{\xi}$ and $V_{\xi}$ be two real numbers endowed with the following properties:
\begin{itemize}
 \item[\textnormal{(i)}] $\exists\,\epsilon_\xi>0$ such that $V_\xi-U_\xi>\epsilon_\xi$ $\forall\xi\in\Xi$;
 \item[\textnormal{(ii)}] $\displaystyle\exists\,\lim_{\xi\rightarrow\bar{\xi}} U_\xi=\lim_{\xi\rightarrow\bar{\xi}} V_\xi=:\lambda_0\in\R$.
\end{itemize}
Moreover, for each $\xi\in\Xi$ let us define:
\begin{itemize}

\item[\textnormal{(a)}] the set
$C_{\xi}:=\left\{\lambda\in \R : U_\xi \leq \lambda < V_\xi\right\}$;
\item[\textnormal{(b)}] the characteristic function of
$C_{\xi}$, i.e., $\mathbf{1}_{C_{\xi}}:\R\rightarrow\{0,1\}$;
\item[\textnormal{(c)}] the number $r_\xi:=V_\xi-U_\xi > \epsilon_\xi$;
\item[\textnormal{(d)}] the function $T_\xi:=
\mathbf{1}_{C_{\xi}}/r_\xi\in L^1_{\mathrm{loc}}\left(\R\right)$ and the corresponding distribution
$\iota_0\left(T_\xi\right)\in\mathcal{D}'_0\left(\R\right)$.
\end{itemize}
Then, it holds that
$\iota_0\left(T_\xi\right) \rightarrow \delta_{\lambda_0}$ in
$\mathcal{D}'_0(\R)$ as $\xi\rightarrow\bar{\xi}$.
\end{lemma}

\proof According to (\ref{kant}), the thesis of the theorem can be recast as
\be \label{bruttetta}
\lim_{\xi\rightarrow \bar{\xi}}\left\langle \iota_0 \left(T_\xi\right),\psi\right\rangle=
\left\langle \delta_{\lambda_0},
\psi\right\rangle\ \ \ \forall \psi\in\mathcal{D}_0\left(\R\right).
\ee
By (\ref{pairint}) and definitions (a)--(d), for any $\psi\in\mathcal{D}_0\left(\R\right)$ we have
\be\label{stravaganza}
\left\langle\iota_0\left(T_\xi\right),\psi\right\rangle=\int_{\R}\frac{\mathbf{1}_{C_{\xi}}(\lambda)}{r_\xi}\,\psi(\lambda)\,d\lambda=
\frac{1}{r_\xi}\int_{C_{\xi}}\psi(\lambda)\,d\lambda=\frac{1}{V_\xi - U_\xi}\int_{U_\xi}^{V_{\xi}}\psi(\lambda)\,d\lambda.
\ee
Moreover, by the integral mean value theorem,
\be\label{cetra}
\exists\,\bar{\lambda}(\xi)\in \left[U_\xi,V_\xi\right]\, :\, \frac{1}{V_\xi - U_\xi}\int_{U_\xi}^{V_{\xi}}\psi(\lambda)\,d\lambda
=\psi\left(\bar{\lambda}(\xi)\right).
\ee
By assumption (ii) we then have $\lim_{\xi\rightarrow\bar{\xi}}\bar{\lambda}(\xi)=\lambda_0$. Since $\psi$ is continuous, this implies
that $\lim_{\xi\rightarrow\bar{\xi}}\psi\left(\bar{\lambda}(\xi)\right)=\psi\left(\lambda_0\right)$. The latter limit, together with
(\ref{stravaganza})--(\ref{cetra}), yields
$\lim_{\xi\rightarrow\bar{\xi}}\left\langle\iota_0\left(T_\xi\right),\psi\right\rangle =\psi(\lambda_0)=
\left\langle \delta_{\lambda_0},\psi\right\rangle$,
i.e., limit (\ref{bruttetta}). This concludes the proof.~$\square$

\vspace{3mm}
Obviously, the one-dimensional counterpart of Corollary \ref{gasquet} is as follows.

\begin{corollary}\label{tsonga}
 For each $j\in\{1,\ldots,J\}$, assume that $U_\xi^j$, $V_\xi^j$, $\lambda_0^j$, $T_\xi^j$
 verify the hypotheses satisfied, respectively, by $U_\xi$, $V_\xi$, $\lambda_0$, $T_\xi$ in Lemma \textnormal{\ref{zarelis}}.
 Moreover, let $\beta_j\in\C$ for all $j\in\{1,\ldots,J\}$. Then
 \be\label{gasparini}
 \lim_{\xi\rightarrow\bar{\xi}}\,\iota_0\left(\sum_{j=1}^J \beta_j\, T_\xi^j\right) =
 \sum_{j=1}^{J}\beta_j\, \delta_{\lambda_0^j}  \ \ \mbox{in}\ \ \mathcal{D}'_0\left(\R\right).
 \ee
\end{corollary}

Finally, the one-dimensional counterpart of Theorem \ref{equidiscreto} is obtained from Corollary \ref{tsonga} by making the following
identifications:
\begin{align}
& \xi=D=d_1,\ \ \bar{\xi}=0,\ \ J=\nu,\ \ \beta_j=\mu_j,\ \
r_{\xi}^j=d_1\ \ \forall j=1,\ldots,\nu,\label{rava1}\\
& U_{\xi}^j=-d_1/2 + F\big(x(P_j)\big),\ \ V_{\xi}^j=d_1/2 + F\big(x(P_j)\big),\ \ \lambda_0^j=F\big(x(P_j)\big),\\
& C_\xi^j=\left\{\lambda\in \R : -d_1/2 \leq \lambda - F\big(x(P_j)\big) <d_1/2 \right\},\\
& \mathbf{1}_{C_\xi^j}(\lambda)=p\big(x(P_j),\lambda;\lambda^\ast,d\big),\ \
T_\xi^j(\lambda)=\frac{\mathbf{1}_{C_\xi^j}(\lambda)}{r_\xi^j}=\frac{p\big(x(P_j),\lambda;\lambda^\ast,d\big)}{d_1},\\
& \sum_{j=1}^J \beta_j\, T_\xi^j(\lambda)=\sum_{j=1}^\nu \mu_j\, \frac{p\big(x(P_j),\lambda;\lambda^\ast,d\big)}{d_1}=
\frac{H(\lambda;\lambda^\ast,d)}{d_1}.\label{fava1}
\end{align}
It is understood that, when considering identifications (\ref{rava1})--(\ref{fava1}), convergence in $\mathcal{D}'_0\left(\R\right)$,
involved in relation (\ref{gasparini}), can be replaced by convergence in $\mathcal{D}'_0\left(\mathcal{T}\right)$,
for any open and bounded investigation domain $\mathcal{T}\subset\R$.

\section{Link between the Radon transform and the Hough transform: the case of piecewise continuous images}\label{ziocontinuo}

The results obtained in the previous section for discrete images can be extended to the case of piecewise continuous images, i.e.,
images described by a function $m\in\mathcal{PD}_0(W)$ (cf. Definition \ref{convDkpezzi} in \ref{ziopiecewise}).
The goal of this section is just to prove the $x$-continuum analogous of
Theorem \ref{equidiscreto}, i.e., roughly speaking, to show again that the rescaled Hough counter of a piecewise continuous image,
defined by analogy with
its discrete counterpart (\ref{WHT}), tends to the generalized Radon transform of the image itself (cf. Definition \ref{radongen})
as the discretization of the parameter space becomes finer and finer. The precise statement of this property will be given in
Theorem \ref{muguruza} at the end of the current section.

The first step is to generalize the formulation of Corollary \ref{gasquet}, by
replacing the discrete sum $\sum_{j=1}^J \beta_j\, T_\xi^j$ with the integral $\int_W \beta(x)\, T_\xi(x;\cdot)\, dx$.

\begin{lemma}\label{bach}
Let $\Xi$ be a subset of $\,\R$ such that $\bar{\xi}$ is an accumulation point for $\Xi$, and let
$\xi\in\Xi$ be a parameter.
For $n,t\in\N\setminus\{0,1\}$, let $W$ and $E'$ be non-empty open subsets of $\R^n$ and $\R^{t-1}$, respectively.
Moreover, for each $\xi\in\Xi$, let
$W\times E' \ni(x,\lambda')\mapsto U_\xi(x;\lambda')\in\R$ and
$W\times E'\ni(x,\lambda')\mapsto V_\xi(x;\lambda')\in\R$
be two functions endowed with the following properties:
\begin{itemize}
 \item[\textnormal{(i)}] $\forall x\in W$, both $U_\xi(x;\cdot)$ and $V_\xi(x;\cdot)$
 are elements of the space $PC^1\left(E'\right)$;
 \item[\textnormal{(ii)}] $\exists\,\epsilon_\xi >0$ such that
 $V_\xi(x;\lambda')-U_\xi(x;\lambda')>\epsilon_{\xi}$ $\forall(x,\lambda')\in W\times E'$;
 \item[\textnormal{(iii)}] $\forall (x,\lambda')\in W\times E'$ $\displaystyle \exists\lim_{\xi\rightarrow \bar{\xi}}U_\xi(x;\lambda')=
 \lim_{\xi\rightarrow \bar{\xi}}V_\xi(x;\lambda')=:G(x;\lambda')\in\R$, with $G(x;\cdot)\in C^1(E')$ $\forall x\in W$;
 \item[\textnormal{(iv)}] the functions $u_\xi:=U_\xi-G$ and $v_\xi:=V_\xi-G$ are uniformly bounded with respect to the parameter $\xi$,
 i.e., there exists a constant $M\geq 0$ such that
 $\,|u_\xi(x;\lambda')|\leq M$ and $\,|v_\xi(x;\lambda')|\leq M$ $\forall(x,\lambda')\in W\times E'$, $\forall\xi\in\Xi$.
\end{itemize}
Furthermore, for each $x\in W$ and $\xi\in\Xi$, let us define:
\begin{itemize}
\item[\textnormal{(a)}] the set
$C_{\xi}(x):=\{\lambda=(\lambda',\lambda_t)\in E'\times\R : U_\xi(x;\lambda')\leq\lambda_t < V_\xi(x;\lambda')\}$;
\item[\textnormal{(b)}] the characteristic function of $C_{\xi}(x)$, i.e., $\mathbf{1}_{C_{\xi}(x)}:E'\times\R\rightarrow\{0,1\}$;
\item[\textnormal{(c)}] the function $r_\xi(x;\cdot):=v_\xi(x;\cdot)-u_\xi(x;\cdot)>\epsilon_\xi$;
\item[\textnormal{(d)}] the function $T_\xi(x;\cdot):=\left[\mathbf{1}_{C_{\xi}(x)}(\cdot)/r_\xi(x;\cdot)\right]\in
L^1_{\mathrm{loc}}\left(E'\times\R\right)$;
\item[\textnormal{(e)}] the function $E'\times\R\ni(\lambda',\lambda_t)\mapsto g(x;\lambda):=\lambda_t-G(x;\lambda')\in\R$ and
the corresponding
Dirac delta $\delta\left(g(x;\cdot)\right)\in\mathcal{D}'_1\left(E'\times\R\right)$.
\end{itemize}
Finally, let $\beta\in\mathcal{PD}_0(W)$ a piecewise continuous function, compactly supported in $W$.
Then
\be\label{provenzale}
\int_{W}\beta(x)\,T_\xi(x;\cdot)\,dx\in L^1_{\mathrm{loc}}\left(E'\times\R\right)\ \ \forall\xi\in\Xi
\ee
and
\be\label{scarlatti}
\lim_{\xi\rightarrow \bar{\xi}}\,\iota_1\left(\int_{W}\beta(x)\,T_\xi(x;\cdot)\,dx\right)=
\int_{W}\beta(x)\,\delta(g(x;\cdot))\,dx\ \ \mbox{in}\ \mathcal{D}'_1\left(E'\times\R\right),
\ee
where the integral on the right-hand side of \textnormal{(\ref{scarlatti})} is to be understood in the sense of Definition
\textnormal{\ref{intparam}} in \ref{test}.
\end{lemma}

\proof
Since both $\beta$ and $T_\xi$ are bounded, we have that
$\beta(\cdot)\,T_\xi(\cdot;\cdot)\in L^1_{\mathrm{loc}}\left(W\times\left(E'\times\R\right)\right)$ $\forall\xi\in\Xi$. Accordingly,
statement (\ref{provenzale}) follows from Fubini theorem.

Moreover, according to (\ref{kant}), statement (\ref{scarlatti}) is equivalent to claiming that, for all
$\psi\in\mathcal{D}_1\left(E'\times\R\right)$,
\be \label{fago}
\lim_{\xi\rightarrow \bar{\xi}}\left\langle \iota_1\left(\int_{W}\beta(x)\,T_\xi(x;\cdot)\,dx\right),\psi\right\rangle
=\left\langle\int_{W}\beta(x)\,\delta(g(x;\cdot))\,dx,\psi\right\rangle.
\ee
Note that, by assumption (iii) and definition (e), we have $g(x;\cdot)\in C^1\left(E'\times\R\right)$ $\forall x\in W$, with
$\partial g(x;\cdot)/\partial\lambda_t=1$ identically. Thus, $\mathrm{grad}_\lambda\,g(x;\lambda)\neq 0$
$\forall (x;\lambda)\in W\times\left(E'\times\R\right)$ and, in particular, $\delta(g(x;\cdot))\in \mathcal{D}'_1\left(E'\times\R\right)$
is well-defined $\forall x\in W$, according to definition (\ref{defdelta}).

In order to prove (\ref{fago}), we begin by observing that an immediate application of Lemma \ref{haendel}, i.e., of equality
(\ref{brunetta}), yields, for each $\psi\in\mathcal{D}_1\left(E'\times\R\right)$,
\be\label{giordani}
\lim_{\xi\rightarrow \bar{\xi}}\beta(x)\langle \iota_1\left(T_\xi(x;\cdot)\right),\psi\rangle=
\beta(x)\langle\delta(g(x;\cdot)),\psi\rangle\ \ \forall x\in W.
\ee
Moreover, if we set $M_{\beta}:=\max_{x\in W}|\beta(x)|$, $S_{\beta}:=\mathrm{supp}\,\beta$, and denote by $\mathbf{1}_{S_\beta}:W\rightarrow
\{0,1\}$ the characteristic function of $S_\beta$ and by $\mathcal{L}^n$ the Lebesgue measure on $\R^n$,
from (\ref{delpoeta}) and (\ref{domina}) (with $u_\xi(\lambda')$ and $v_\xi(\lambda')$ replaced by $u_\xi(x;\lambda')$ and
$v_\xi(x;\lambda')$, respectively) we have
\be\label{ferrandini}
\left|\beta(x)\langle \iota_1\left(T_\xi(x;\cdot)\right),\psi\rangle\right|\leq
M_\beta\, \|\psi\|_{C^1}\,\mathcal{L}^n\left(K_M\right)\,\mathbf{1}_{S_\beta}(x).
\ee
By assumption, $S_\beta\subset W$ is compact and then $\mathbf{1}_{S_\beta}\in L^1(W)$. Hence, from (\ref{giordani})--(\ref{ferrandini}) and
Lebesgue dominated convergence theorem, we find, for each $\psi\in\mathcal{D}_1\left(E'\times\R\right)$,
\be\label{radwanska}
\lim_{\xi\rightarrow \bar{\xi}}\int_{W} \beta(x)\langle \iota_1\left(T_\xi(x;\cdot)\right),\psi\rangle\, dx =
\int_{W}\beta(x)\langle\delta(g(x;\cdot)),\psi\rangle \,dx.
\ee

Now, we are going to show that (\ref{radwanska}) is actually thesis (\ref{fago}).
Indeed, as far as the left-hand side of (\ref{radwanska}) is concerned, by (\ref{pairint}) we can write
\be\label{azarenka}
\int_{W} \beta(x)\langle \iota_1\left(T_\xi(x;\cdot)\right),\psi\rangle\, dx=
\int_{W} \beta(x)\left[\int_{E'\times\R}T_{\xi}(x;\lambda)\,\psi(\lambda)\, d\lambda\right]dx.
\ee
From the properties of functions $\beta$, $T_\xi$ and $\psi$ (cf., in particular, definitions (b)--(d)), it follows that their product
is bounded and compactly supported, and then $\beta(\cdot)\,T_\xi(\cdot;\cdot)\,\psi(\cdot) \in
L^1\left(W\times\left(E'\times\R\right)\right)$ $\forall\xi\in\Xi$.
Hence, by Fubini theorem, we can rewrite
(\ref{azarenka}) as
\be\label{williams}
\int_{W} \beta(x)\langle \iota_1\left(T_\xi(x;\cdot)\right),\psi\rangle\, dx=
\int_{E'\times\R} \left[\int_{W}\beta(x)\,T_{\xi}(x;\lambda)\, dx\right]\psi(\lambda)\,d\lambda.
\ee
Moreover, property (\ref{provenzale}) and definition (\ref{pairint}) imply that
relation (\ref{williams}) can be written in the form
\be\label{ivanovic}
\int_{W} \beta(x)\langle \iota_1\left(T_\xi(x;\cdot)\right),\psi\rangle\, dx=
 \left\langle \iota_1\left(\int_{W}\beta(x)\,T_{\xi}(x;\cdot)\, dx\right),\psi\right\rangle.
\ee
An immediate comparison between (\ref{radwanska}) and (\ref{ivanovic}) now shows that the left-hand side of (\ref{radwanska})
coincides with that of statement (\ref{fago}).

As far as the right-hand side of (\ref{radwanska}) is concerned, we remember that, according to definition (e),
$g(x;\lambda)=\lambda_t-G(x;\lambda')$. Thus, by (\ref{hume}), (\ref{defdelta}) and (\ref{viotti}) in the appendix, we have
\be\label{schubert}
\left\langle \beta(x)\,\delta(g(x;\cdot)),\psi\right\rangle=
\beta(x)\left\langle\delta(g(x;\cdot)),\psi\right\rangle = \beta(x)\int_{E'}\psi\left(\lambda', G(x;\lambda')\right)d\lambda'.
\ee
From the properties of the functions $\beta$ and $\psi$, it follows that their product
is bounded and compactly supported: then,
$\beta(\cdot)\,\psi(\cdot,G(\cdot;\cdot))\in L^1\left(W\times E'\right)$.
Hence, by Fubini theorem, we have that $\beta(\cdot)\int_{E'}\psi\left(\lambda',G(\cdot;\lambda')\right)d\lambda'\in L^1(W)$ and
\be \label{luchesi}
\int_{W}\beta(x)\left[\int_{E'}\psi\left(\lambda',G(x;\lambda')\right)d\lambda'\right]dx=
\int_{W\times E'}\beta(x)\,\psi\left(\lambda',G(x;\lambda')\right)\,dx\, d\lambda'.
\ee
Now, it is clear that the right-hand side of (\ref{luchesi}) defines a linear and continuous functional on
$\mathcal{D}_1\left(E'\times\R\right)\ni\psi$, i.e., an element $\ell\in\mathcal{D}'_1\left(E'\times\R\right)$. Indeed, linearity is
obvious, while continuity readily follows from Lebesgue dominated
convergence theorem and the notion of convergence in $\mathcal{D}_1\left(E'\times\R\right)$, as given by Definition \ref{convDk}
in \ref{test}. In agreement with Definition \ref{intparam} in \ref{test}, such functional $\ell$ is denoted by
$\int_{W}\beta(x)\,\delta(g(x;\cdot))\,dx$.
Summing up, from (\ref{schubert}) and (\ref{luchesi}) we find
\be\label{schumann}
\int_{W} \beta(x)\left\langle\delta(g(x;\cdot)),\psi\right\rangle dx=\left\langle\int_{W}\beta(x)\,\delta(g(x;\cdot))\,dx,\,\psi\right\rangle,
\ee
thus showing that the right-hand side of (\ref{radwanska}) coincides with that of the claimed assertion (\ref{fago}).
This concludes the proof.~$\square$

\vspace{3mm}
It is worth noting that, according to our notation, \textit{a priori}
\be\label{tensiun}
\int_{W} \beta(x)\,\delta\left(g(x;\cdot)\right)\,dx \neq \int_{W} \delta\left(g(x;\cdot)\right)\,\beta(x)\,dx.
\ee
Indeed, the left-hand side of (\ref{tensiun}) is the integral of the one-parameter family of distributions
$\beta(x)\,\delta\left(g(x;\cdot)\right)\in\mathcal{D}'_0\left(E'\times\R\right)$ with respect to the parameter $x\in W$, in the sense of
Definition \ref{intparam}.
Instead, according to (\ref{radongeneq}), the right-hand side of (\ref{tensiun}) denotes the generalized Radon transform
$\left(R_g\,\beta\right)$, i.e., the map defined by $\lambda\mapsto \langle\delta\left(g(\cdot;\lambda)\right), \beta \rangle$,
in agreement with the integral notation adopted in (\ref{defdelta}) for the pairing.

However, the following lemma states that, under appropriate hypotheses,
the two sides of (\ref{tensiun}) are equal as elements of $\mathcal{D}'_0\left(E'\times\R\right)$.
\begin{lemma}\label{beethoven}
For $n,t\in \N\setminus\{0,1\}$, let $W$ and $E'$ be non-empty open subsets of $\,\R^n$ and $\,\R^{t-1}$ respectively, and let
$\beta\in\mathcal{PD}_0(W)$. Moreover, let $G:W\times E'\rightarrow\R$ be a function such that
$G(\cdot;\lambda')\in C^1(W)$ $\forall \lambda'\in E'$ and $G(x;\cdot)\in C^1(E')$ $\forall x\in W$.
Finally, assume that the function $g:W\times \left(E'\times\R\right)\rightarrow\R$ defined as
$g(x;\lambda):=\lambda_t-G(x;\lambda')$ satisfies the following two conditions:
\begin{itemize}
\item[\textnormal{(i)}] $\forall\lambda\in E'\times\R$ such that $\mathcal{S}(\lambda):=\{x\in W : g(x;\lambda)=0\}\neq\emptyset$,
it holds that $\mathrm{grad}_x\, g(x;\lambda)\neq 0 \ \forall x\in\mathcal{S}(\lambda)$;
\item[\textnormal{(ii)}] $\displaystyle \int_{W} \delta\left(g(x;\cdot)\right)\,\beta(x)\,dx\in L^1_{\mathrm{loc}}\left(E'\times\R\right)$.
\end{itemize}
Then, as elements of $\mathcal{D}_0'\left(E'\times\R\right)$,
\be\label{muffat}
\int_{W} \beta(x)\,\delta\left(g(x;\cdot)\right)\,dx =\iota_0\left(\int_{W} \delta\left(g(x;\cdot)\right)\,\beta(x)\,dx\right).
\ee
\end{lemma}

\proof
Proving equality (\ref{muffat}) amounts to proving that, for all $\psi\in \mathcal{D}_0\left(E'\times\R\right)$,
\be\label{muffatpsi}
\left\langle \int_{W} \beta(x)\,\delta\left(g(x;\cdot)\right)\,dx,\,\psi\right\rangle=
\left\langle\iota_0\left(\int_{W} \delta\left(g(x;\cdot)\right)\,\beta(x)\,dx\right),\,\psi \right\rangle.
\ee

Now, from (\ref{schubert})--(\ref{schumann}), we immediately find
\be\label{chopin}
\left\langle \int_{W} \beta(x)\,\delta\left(g(x;\cdot)\right)\,dx,\,\psi\right\rangle=
\int_{W\times E'}\beta(x)\,\psi\left(\lambda',G(x;\lambda')\right)\,dx\, d\lambda'.
\ee
As already observed just below (\ref{schubert}), $\beta(\cdot)\,\psi\left(\cdot,G(\cdot;\cdot)\right)\in L^1\left(W\times E'\right)$. Thus, by
Fubini theorem, we can rewrite (\ref{chopin}) as
\be\label{mendelssohn}
\left\langle \int_{W} \beta(x)\,\delta\left(g(x;\cdot)\right)\,dx,\,\psi\right\rangle=
\int_{E'}\left[\int_W \beta(x)\,\psi\left(\lambda',G(x;\lambda')\right)\,dx\right] d\lambda'.
\ee

Now, for each $\lambda'\in E'$, we introduce the auxiliary functions $\psi_{\lambda'}:\R\rightarrow\R$ and $G_{\lambda'}:W\rightarrow\R$,
defined by the conditions
\be\label{liszt}
\psi_{\lambda'}(\lambda_t):=\psi(\lambda',\lambda_t)=\psi(\lambda),\  G_{\lambda'}(x):=G(x;\lambda')
\ \ \forall (x,\lambda',\lambda_t)\in W\times E'\times\R.
\ee
From the expression of $g(x;\lambda):=\lambda_t-G(x;\lambda')$ and definitions (\ref{liszt}), it follows that
\be\label{zitromax}
G^{-1}_{\lambda'}(\lambda_t)=\{x\in W : G_{\lambda'}(x)=\lambda_t\}=\{x\in W : g(x;\lambda)=0\}=\mathcal{S}(\lambda),
\ee
as well as $\mathrm{grad}\,G_{\lambda'}(x)=\mathrm{grad}_x\,g(x;\lambda)$ for all $(x,\lambda)\in W\times\left(E'\times\R\right)$.

Then, by assumption (i), for each $\lambda'\in E'$ we can apply the coarea
formula\footnote{Cf. also the short discussion just below (\ref{coarea}) itself. The identifications to be done in (\ref{coarea}) to obtain
(\ref{brahms}) are the following ones: $W=A$, $\beta(x)\,\psi_{\lambda'}\left(G_{\lambda'}(x)\right)=g(x)$, $G_{\lambda'}=\Psi$, $\lambda_t=s$.}
(\ref{coarea}) to the internal integral on the right-hand side of (\ref{mendelssohn}),
thus obtaining
\be\label{brahms}
\int_W \beta(x)\,\psi_{\lambda'}\left(G_{\lambda'}(x)\right)\,dx =
\int_{G_{\lambda'}(W)}\left[\int_{\mathcal{S}(\lambda)}\frac{\beta(x)\,\psi_{\lambda'}(\lambda_t)}
{|\mathrm{grad}_x\,g(x;\lambda)|}\,d\sigma(x)\right]d\lambda_t.
\ee
Taking into account
(\ref{liszt}), (\ref{defdelta}), (\ref{desigma}) and observing that, by (\ref{zitromax}), if
$\lambda_t\notin G_{\lambda'}(W)$ then $\mathcal{S}(\lambda)=\emptyset$,
we can rewrite (\ref{brahms}) as
\be\label{busoni}
\int_W \beta(x)\,\psi\left(\lambda',G_{\lambda'}(x)\right)\,dx=
\int_{\R}\psi(\lambda',\lambda_t)\left[\int_{W}\delta\left(g(x;\lambda)\right)\beta(x)\,dx\right]d\lambda_t.
\ee
Next, by substituting (\ref{busoni}) into (\ref{mendelssohn}), we find
\begin{align}\label{breva}
& \left\langle \int_{W} \beta(x)\,\delta\left(g(x;\cdot)\right)\,dx,\,\psi\right\rangle \\
&\ \ \ \ \ \ \ \ \ \ \ \ \ \ \ \ \ \ \ \ \ \ \
=\int_{E'}\left\{\int_{\R}\psi(\lambda)
\left[\int_{W}\delta\left(g(x;\lambda)\right)\beta(x)\,dx\right]d\lambda_t\right\}d\lambda'.\nonumber
\end{align}
From hypothesis (ii) and the fact that $\psi\in\mathcal{D}_0\left(E'\times\R\right)$, it follows that the map defined by
\be
E'\times\R\ni\lambda\mapsto \psi(\lambda) \int_{W}\delta\left(g(x;\lambda)\right)\beta(x)\,dx\in\R
\ee
is in $L^1_{\mathrm{loc}}\left(E'\times\R\right)$. Thus, by Fubini theorem, relation (\ref{breva}) can be written as
\be\label{polaramin}
\left\langle \int_{W} \beta(x)\,\delta\left(g(x;\cdot)\right)\,dx,\,\psi\right\rangle =
\int_{E'\times\R}\psi(\lambda)\left[\int_{W}\delta\left(g(x;\lambda)\right)\beta(x)\,dx\right]d\lambda.
\ee
Finally, the same hypothesis (ii) and definition (\ref{pairint}) easily allow recognizing that (\ref{polaramin})
coincides with equality (\ref{muffatpsi}), as wanted.~$\square$

\begin{remark}\label{luxilon}
If $E$ is an open subset of $E'\times\R$ and hypothesis (ii) in Lemma \ref{beethoven} is reformulated as (ii$'$)
$\int_{W} \delta\left(g(x;\cdot)\right)\,\beta(x)\,dx\in L^1_{\mathrm{loc}}(E)$, then thesis (\ref{muffat})
still holds true, provided that it be regarded as an equality between elements of $\mathcal{D}'_0(E)$. The proof is almost
identical to that of Lemma \ref{beethoven} itself. However, as an assumption on $g$, property (ii), or (ii$'$), in Lemma \ref{beethoven}
is rather implicit and, in principle, its fulfilment depends not only on
$g$, but also on $\beta$. Accordingly, it is important to establish sufficient and more explicit conditions on $g$ only,
ensuring that such a property holds true for all $\beta\in\mathcal{PD}_0(W)$.
This task has already been accomplished by Theorem \ref{monteverdi}.
\end{remark}

We can now come back to our main problem, i.e., the link between the Radon transform and the Hough transform. To this end, we first need to
formulate an appropriate definition of the weighted Hough accumulator $H(\lambda;\lambda^\ast,d)$ for a piecewise continuous image
$m\in\mathcal{PD}_0(W)$. Taking inspiration from (\ref{WHT}), it is natural to define\footnote{See also footnote no. \ref{zioambiguo}.}
\be\label{Hcontinuo}
H(\lambda;\lambda^\ast,d):=\int_W m(x)\,p(x,\lambda;\lambda^\ast,d)\,dx,
\ee
where $p(x,\lambda;\lambda^\ast,d)$ is the Hough transform kernel (\ref{nuovokernel}).
Then, we can state the main result of this section as follows.
\begin{theorem}\label{muguruza}
For $n,t\in \N\setminus\{0,1\}$, let $W$ and $E$ be non-empty open subsets of $\,\R^n$ and $\,\R^{t}$, respectively.
Moreover, let $f:W\times E\rightarrow\R$ be a function satisfying the following properties:
\textnormal{(a)} ${\mathcal S}(\lambda):=\{x\in W : f_{ \lambda}(x)= 0\}\neq \emptyset$ $\forall\lambda\in E$;
\textnormal{(b)} $f\in C^1(W\times E)$;  \textnormal{(c)} $(\mathrm{grad}_x\, f_{\lambda})(x)\neq 0$ $\forall x\in  {\mathcal S}(\lambda)$, $\forall\lambda\in E$;
 \textnormal{(d)} $f$ is $\lambda_t$-solvable, i.e., $f(x;\lambda)= \lambda_t-F(x;\lambda')$, with $\lambda=(\lambda',\lambda_t)$.
Moreover, let $\{\lambda^\ast,d\}$ be a discretization of the parameter space, and define
$D:=\max\{d_1,\ldots,d_t\}$, where $d_k$, for $k=1,\ldots,t$, is the sampling distance with respect to the component $\lambda_k$,
as explained in item \textnormal{I} of Subsection \textnormal{\ref{WHC}}.
Finally, let $m\in\mathcal{PD}_0(W)$ be a piecewise continuous and compactly supported image,
$(R_f\, m)(\lambda)$ its generalized Radon transform
and $H(\lambda;\lambda^\ast,d)/d_t$ the corresponding rescaled Hough
counter\footnote{Cf. definitions (\ref{hume}), (\ref{radongeneq}) and (\ref{Hcontinuo}), respectively.}, defined
on a bounded and open investigation domain $\mathcal{T}\subset E$. Then
\be\label{minacce2}
\lim_{D\rightarrow 0^+}\,\iota_1\left(\frac{H(\lambda;\lambda^\ast,d)}{d_t}\right)=
(R_f\, m)(\lambda)
\ \ \mbox{in}\ \ \mathcal{D}'_1\left(\mathcal{T}\right),
\ee
where $\iota_1:L^1_{\mathrm{loc}}(\mathcal{T})\rightarrow \mathcal{D}'_{1}(\mathcal{T})$
denotes the inclusion map defined as in \textnormal{(\ref{pairint})}.
\end{theorem}

\proof First, we choose or identify the functions and parameters appearing in the statements of Lemmas \ref{bach} and
\ref{beethoven} as follows:
\begin{align}
& \xi=D=\max\{d_1,\ldots,d_t\},\ \ \bar{\xi}=0,\ \ \beta=m,\ \ r_{\xi}(x;\cdot)=d_t\ \forall x\in W,\label{rava2} \\
& U_{\xi}(x;\lambda')=-d_t/2 + F\big(x;c'(\lambda')\big),\ \ V_{\xi}(x;\lambda')=d_t/2 + F\big(x;c'(\lambda')\big),\\
& G(x;\lambda')=F(x; \lambda'),\ \ g(x;\lambda)=f(x;\lambda)=\lambda_t-F(x;\lambda'),\\
& C_\xi(x)=\left\{\lambda=(\lambda',\lambda_t)\in E : -d_t/2 \leq \lambda_t - F\big(x;c'(\lambda')\big) <d_t/2 \right\},\\
& \mathbf{1}_{C_\xi(x)}(\lambda)=p(x,\lambda;\lambda^\ast,d),\ \
T_\xi(x;\lambda)=\frac{\mathbf{1}_{C_\xi(x)}(\lambda)}{r_\xi(x;\lambda)}=\frac{p(x,\lambda;\lambda^\ast,d)}{d_t},\\
& \int_W \beta(x)\, T_\xi(x;\lambda)\,dx=\int_W m(x)\,\frac{p(x,\lambda;\lambda^\ast,d)}{d_t}\,dx=\frac{H(\lambda;\lambda^\ast,d)}{d_t}.
\label{fava2}
\end{align}
An easy check shows that identifications (\ref{rava2})--(\ref{fava2}) ensure the fulfilment of the hypotheses required by
Lemmas \ref{bach}, \ref{beethoven} and Remark \ref{luxilon}. Accordingly, we find that
\be\label{minacce3}
\lim_{D\rightarrow 0^+}\,\iota_1\left(\frac{H(\cdot;\lambda^\ast,d)}{d_t}\right)=
\int_W \delta(f(x;\cdot))\,m(x)\,dx
\ \ \mbox{in}\ \ \mathcal{D}'_1\left(\mathcal{T}\right).
\ee
Finally, it suffices to recall that the right-hand side of equality (\ref{minacce3}) is just the generalized Radon transform of
a piecewise continuous image $m\in\mathcal{PD}_0(W)$, as shown in definition (\ref{radongeneq}). Accordingly, relation (\ref{minacce3})
can be equivalently rewritten as (\ref{minacce2}), thus proving Theorem \ref{muguruza} and justifying the claims opening this
section.~$\square$

\subsection{The one-dimensional case $t=1$}

For sake of completeness, let us now sketch how the previous investigation trivializes when $t=1$.
The one-dimensional counterpart of Lemma \ref{bach} can be formulated as follows.

\begin{lemma}\label{bach1}
Let $\Xi$ be a subset of $\,\R$ such that $\bar{\xi}$ is an accumulation point for $\Xi$, and let
$\xi\in\Xi$ be a parameter. For $n\in\N\setminus\{0,1\}$, let $W$ be a non-empty open subset of $\R^n$.
Moreover, for each $\xi\in\Xi$, let $U_\xi:W\rightarrow\R$ and $V_\xi:W\rightarrow\R$
be two functions endowed with the following properties:
\begin{itemize}
 \item[\textnormal{(i)}] $\exists\,\epsilon_\xi >0$ such that
 $V_\xi(x)-U_\xi(x)>\epsilon_{\xi}$ $\forall x\in W$;
 \item[\textnormal{(ii)}] $\forall x\in W$ $\displaystyle \exists\lim_{\xi\rightarrow \bar{\xi}}U_\xi(x)=
 \lim_{\xi\rightarrow \bar{\xi}}V_\xi(x)=:G(x)\in\R$.
 \end{itemize}
Furthermore, for each $x\in W$ and $\xi\in\Xi$, let us define:
\begin{itemize}
\item[\textnormal{(a)}] the set
$C_{\xi}(x):=\{\lambda\in\R : U_\xi(x)\leq\lambda < V_\xi(x)\}$;
\item[\textnormal{(b)}] the characteristic function of $C_{\xi}(x)$, i.e., $\mathbf{1}_{C_{\xi}(x)}:\R\rightarrow\{0,1\}$;
\item[\textnormal{(c)}] the number $r_\xi(x):=V_\xi(x)-U_\xi(x)>\epsilon_\xi$;
\item[\textnormal{(d)}] the function $T_\xi(x;\cdot):=\left[\mathbf{1}_{C_{\xi}(x)}(\cdot)/r_\xi(x)\right]\in
L^1_{\mathrm{loc}}\left(\R\right)$.
\end{itemize}
Finally, let $\beta\in\mathcal{PD}_0(W)$. Then
\be\label{provenzale1}
\int_{W}\beta(x)\,T_\xi(x;\cdot)\,dx\in L^1_{\mathrm{loc}}\left(\R\right)\ \ \forall\xi\in\Xi
\ee
and
\be\label{scarlatti1}
\lim_{\xi\rightarrow \bar{\xi}}\,\iota_0\left(\int_{W}\beta(x)\,T_\xi(x;\cdot)\,dx\right)=
\int_{W}\beta(x)\,\delta\big(\cdot- G(x)\big)\,dx\ \ \mbox{in}\ \mathcal{D}'_0\left(\R\right),
\ee
where $\delta\big(\cdot - G(x)\big)=\delta_{G(x)}\in \mathcal{D}'_0\left(\R\right)$ is the Dirac delta centred at the point $G(x)$,
and the integral on the right-hand side of \textnormal{(\ref{scarlatti1})} is to be understood in the sense of
Definition \textnormal{\ref{intparam}} in \ref{test}.
\end{lemma}

\proof
The result follows by adapting and simplifying the proof of Lemma \ref{bach}.~$\square$

\vspace{3mm}
The one-dimensional counterpart of Lemma \ref{beethoven} can be formulated as follows.
\begin{lemma}\label{beethoven1}
For $n\in\N\setminus\{0,1\}$, let $W$ be a non-empty open subset of $\,\R^n$, and let
$\beta\in\mathcal{PD}_0(W)$. Moreover, let $G\in C^1(W)$ be a real-valued function satisfying the following two conditions:
\begin{itemize}
\item[\textnormal{(i)}] $\forall\lambda\in\R$ such that $\mathcal{S}(\lambda):=\{x\in W : \lambda - G(x)=0\}\neq\emptyset$, it holds that
$\mathrm{grad}\, G(x)\neq 0$ $\forall x\in\mathcal{S}(\lambda)$;
\item[\textnormal{(ii)}] $\int_{W} \delta\big(\cdot - G(x)\big)\,\beta(x)\,dx\in L^1_{\mathrm{loc}}\left(\R\right)$.
\end{itemize}
Then, as elements of $\mathcal{D}_0'\left(\R\right)$,
\be\label{muffat1}
\int_{W} \beta(x)\,\delta\big(\cdot - G(x)\big)\,dx =\iota_0\left(\int_{W} \delta\big(\cdot - G(x)\big)\,\beta(x)\,dx\right).
\ee
\end{lemma}

\proof
The result follows by adapting and simplifying the proof of Lemma \ref{beethoven}.~$\square$

\vspace{3mm}
Note that Theorem \ref{monteverdi}, mentioned in Remark \ref{luxilon}, already comprise the case $t=1$. Then, the
one-dimensional counterpart of Theorem \ref{muguruza} is obtained from Lemmas \ref{bach1}, \ref{beethoven1} and Remark \ref{luxilon}
by making the following identifications:
\begin{align}
& \xi=D=d_1,\ \ \bar{\xi}=0,\ \ \beta=m,\ \ r_{\xi}(x;\cdot)=d_1\ \forall x\in W,\label{rava3} \\
& U_{\xi}(x)=-d_1/2 + F(x),\ \ V_{\xi}(x)=d_1/2 + F(x),\\
& G(x)=F(x),\ \ g(x;\lambda)=f(x;\lambda)=\lambda-F(x),\\
& C_\xi(x)=\left\{\lambda\in \R : -d_1/2 \leq \lambda - F(x) <d_1/2 \right\},\\
& \mathbf{1}_{C_\xi(x)}(\lambda)=p(x,\lambda;\lambda^\ast,d),\ \
T_\xi(x;\lambda)=\frac{\mathbf{1}_{C_\xi(x)}(\lambda)}{r_\xi(x;\lambda)}=\frac{p(x,\lambda;\lambda^\ast,d)}{d_1},\\
& \int_W \beta(x)\, T_\xi(x;\lambda)\,dx=\int_W m(x)\,\frac{p(x,\lambda;\lambda^\ast,d)}{d_1}\,dx=\frac{H(\lambda;\lambda^\ast,d)}{d_1},
\label{fava3}
\end{align}
so that relation (\ref{minacce2}) now reads
\be\label{minacce4}
\lim_{D\rightarrow 0^+}\,\iota_0\left(\frac{H(\lambda;\lambda^\ast,d)}{d_1}\right)=
(R_f\, m)(\lambda)
\ \ \mbox{in}\ \ \mathcal{D}'_0\left(\mathcal{T}\right),
\ee
for any open and bounded investigation domain $\mathcal{T}\subset\R$.

\section{Applications and numerical examples}\label{zionumerico}

The investigation performed in the previous sections has highlighted a complex relationship between the Radon transform and the Hough
transform. That is, according to Theorem \ref{equidiscreto} and Theorem \ref{muguruza}, the rescaled (weighted) Hough counter
of either a discrete or a piecewise continuous image tends to the generalized Radon transform of the image itself
as the discretization of the parameter space becomes finer and finer.

Although this is an interesting result in itself, an issue naturally arises concerning its possible numerical applications.
Here we just outline a new inversion technique for visualizing an object from a very noisy Radon sinogram, by regarding
the latter as an approximation of a rescaled Hough sinogram (cf. Definition \ref{Hsinog}). This possibility is suggested, in particular,
by limit (\ref{minacce2}), and it is worth
investigating, since there are cases (e.g., Positron Emission Tomography) in which a high level of noise affects the Radon sinogram,
thus preventing the traditional (i.e., Radon-based) inversion techniques from providing satisfactory reconstruction of the unknown object.

\begin{figure}[t]
\begin{center}
\includegraphics[scale=0.4]{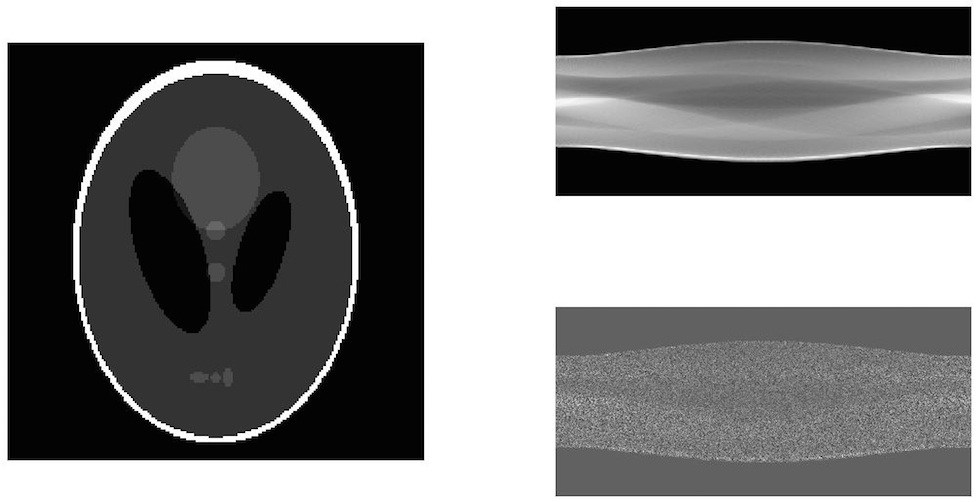}\\
(a)\ \ \ \ \ \ \ \ \ \ \ \ \ \ \ \ \ \ \ \ \ \ \ \ \ \ \ \ \ \ \ \ \ \ \ \ \ \ \ \ \ \ \ \ \ \ \ \ (b)
\caption{(a) The Shepp--Logan phantom. (b) Upper plot: the Radon noise-free sinogram of the Shepp--Logan phantom, computed for $629$ values of
$\vartheta\in [0,\pi)$ (horizontal axis) and $287$ values of $\gamma\in[-\sqrt{2},\sqrt{2}]$ (vertical axis). Lower plot: the noisy sinogram of the
Shepp--Logan phantom, obtained from the upper one by corrupting it with additive Gaussian noise at a level $\ell = 100\%$, according to formula (\ref{eq:sinonoise}).}
\label{phantom}
\end{center}
\end{figure}

Then, we consider the well-known Shepp--Logan phantom \cite{shlo74}, shown in panel (a) of Figure \ref{phantom}, as the piecewise constant
image\footnote{This image is contained in a square of sides ranging from $-1$ to $1$ and is formed by pixels with values ranging from $0$
to $1$.} to be recovered from a very noisy Radon sinogram. To this end, we first compute the exact Radon transform with respect to the
family of straight lines of equation
\be\label{retta2}
\gamma-x_1\cos\vartheta-x_2\sin\vartheta=0,
\ee
which is of the form $f(x;\lambda)=0$, under the identifications
$x=(x_1,x_2)$, $\lambda=(\lambda_1,\lambda_2)=(\vartheta,\gamma)$ and $f(x;\lambda)= \lambda_2 -x_1\cos\lambda_1-x_2\sin\lambda_1$ (cf. the notation adopted in Remark \ref{def-sinogramma}). As explained in Subsection \ref{zioquadrato}, such computation can be easily performed by means
of an appropriate implementation\footnote{In particular, owing to the different forms of the straight line equation adopted in Subsection
\ref{zioquadrato} and in the current section (compare (\ref{retta1}) with (\ref{retta2})), in formula (\ref{radonquadrato}) the
following substitutions should be made: $\gamma\mapsto\gamma/\sin\vartheta$, $\omega_1\mapsto\cot\vartheta$.}
of formula (\ref{radonquadrato}), by considering $I$ discretized values $\vartheta_i$
(with $i=1,\ldots,I$) of $\vartheta\in [0, \pi)$ and $J$ discretized values $\gamma_j$ (with $j=1,\ldots,J$) of
$\gamma\in \left[-\sqrt{2},\sqrt{2}\right]$.
The corresponding noise-free sinogram, obtained for $I=629$ and $J=287$, is represented in the upper part of panel (b) in Figure
\ref{phantom}.

The noise-free sinogram is then corrupted by additive Gaussian noise by using the formula
\be\label{eq:sinonoise}
S_n(\vartheta_i,\gamma_j)=S_t(\vartheta_i,\gamma_j)+\ell\,\varepsilon\, S_t(\vartheta_i,\gamma_j),
\ee
where
\begin{itemize}
\item  $\vartheta_i$, for $i = 1\dots I$, are discretized values of $\vartheta$;
\item $\gamma_j$, for $j = 1\dots J$, are discretized values of $\gamma$;
\item $S_t(\vartheta_i,\gamma_j)$ is the true value of the sinogram at the point  $(\vartheta_i,\gamma_j)$;
\item $S_n(\vartheta_i,\gamma_j)$ is the noisy value of the sinogram at the point $(\vartheta_i,\gamma_j)$;
\item $\varepsilon$ is a realization of a normal Gaussian random variable;
\item $\ell$ is the percentage noise level ($\ell = 100\%$).
\end{itemize}
The resulting noisy Radon sinogram is shown in the lower part of panel (b) in Figure \ref{phantom}.

\begin{figure}[t]
\begin{center}
\includegraphics[scale=0.4]{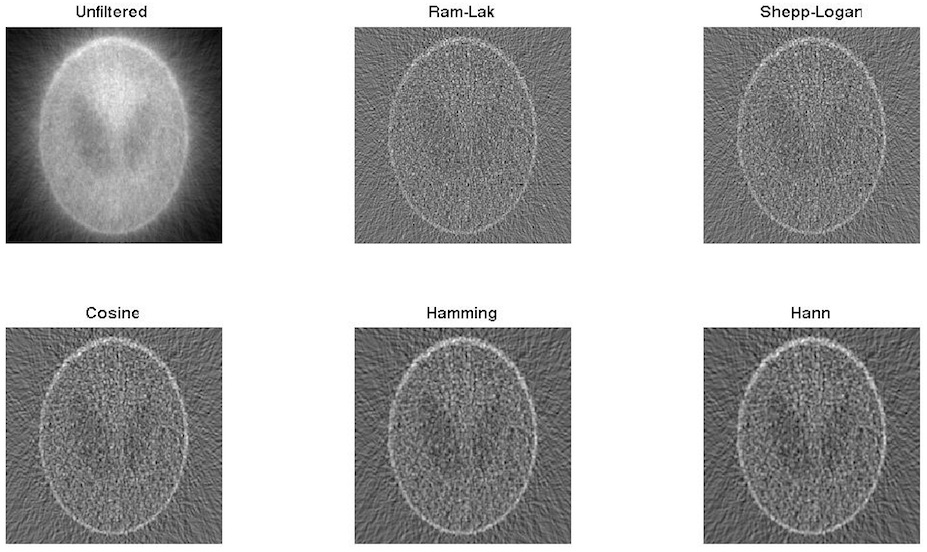}\\
\caption{Reconstructions of the Shepp--Logan phantom, obtained from the noisy sinogram (shown in the bottom part of Figure
\ref{phantom}, panel (b)) by using the unfiltered back-projection, and the FBP algorithm with five different choices for the filtering function. Except for the case of the unfiltered back-projection, the internal structure of the phantom is almost completely lost.}
\label{radonA}
\end{center}
\end{figure}

Usually, the inversion of the Radon transform is numerically performed by using the filtered back-projection (FBP) algorithm, where the presence of a ramp filter (Ram--Lak filter) in the frequency domain attenuates the blurring effect of a crude unfiltered back-projection and where, at the same time, a second filtering function multiplying the ramp filter allows the attenuation of high frequency noise that can be present in the sinogram. Common choices for this filtering function are \cite{bebo98, nawu01} the Shepp--Logan filter (a sinc function); the Cosine filter (a cosine function); the Hamming window; the Hann window. We have then applied both the unfiltered and the filtered back-projection to recover the Shepp--Logan image from its noisy Radon sinogram, using all the cited filters in the case of the FBP algorithm. The corresponding results are shown in Figure \ref{radonA}. It is clear that, independently of the particular filter adopted, the FBP algorithm fails to recover the internal structure of the phantom, while the unfiltered back-projection can at least visualize its main features.

\begin{figure}[t]
\begin{center}
\includegraphics[scale=0.4]{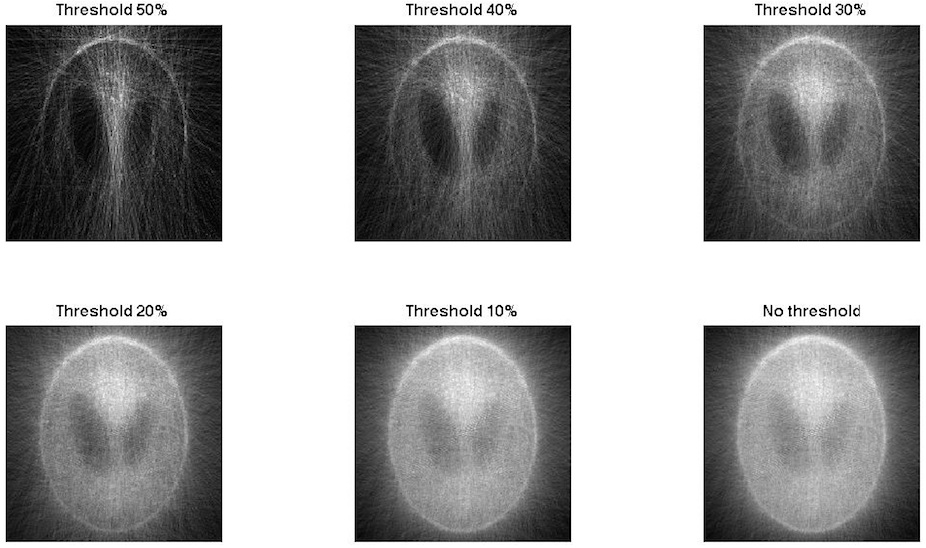}\\
\caption{Visualizations of the Shepp--Logan phantom, obtained by drawing straight lines identified by pairs of parameters corresponding to cells in the Hough counter (i.e., the noisy Radon sinogram shown in the bottom part of Figure \ref{phantom}, panel (b)) whose values are higher than a fixed percentage of the maximum value. Five different thresholds are chosen, while ``no threshold'' means that all the pairs of parameters related to non-empty cells are used to identify straight lines in the image space.}
\label{hough}
\end{center}
\end{figure}

Let us now see what kind of visualization we can obtain by regarding the noisy Radon sinogram as an approximation of a rescaled Hough
sinogram, as suggested by limit (\ref{minacce2}). Each pixel of the noisy Radon sinogram is regarded as a cell of centre
$(\vartheta_i,\gamma_j)$ in the parameter space, and the value of the pixel, multiplied by the sampling distance\footnote{Cf. definition (\ref{IN}), for $t=2$.} $d_2$ with respect to the component $\lambda_2=\gamma$, is regarded as the number of straight
lines characterized by parameters $(\vartheta_i,\gamma_j)$ and to be considered in the image space. Note that this number need not be an
integer. In fact, more precisely, all the pixel values in the image space $\A^2_{(x_1,x_2)}(\R)$
are initialized to zero and then, for any pixel centred at $(\vartheta_i,\gamma_j)$ in the parameter space $\A^2_{(\vartheta,\gamma)}(\R)$ and having value $S_n(\vartheta_i,\gamma_j)$, we trace back in $\A^2_{(x_1,x_2)}(\R)$ a straight line of equation $\gamma_j-x_1\cos\vartheta_i-x_2\sin\vartheta_i=0$, and the value of each pixel crossed by this straight line is increased by $d_2 S_n(\vartheta_i,\gamma_j)$. The resulting visualization is shown in the bottom-right panel of Figure \ref{hough}.

It is also interesting to implement the above procedure by taking into account only the ``principal'' pixels, i.e., the pixels whose
values are larger than a certain threshold. Various thresholds are considered, as five different percentages of the maximum value of the pixels
in the noisy Radon sinogram (multiplied by $d_2$). The corresponding visualizations are shown in the first five panels of Figure \ref{hough}.

It is worth observing that, unlike Figure \ref{radonA}, the pixel values in the panels of Figure \ref{hough} are not
related, in principle, to the true values of the Shepp--Logan phantom. However, a visual comparison between Figure \ref{radonA} and
\ref{hough} suggests that, for appropriate values of the threshold, our new (Hough-based) approach can provide visualizations
that are more informative and accurate than those provided by the (Radon-based) FBP algorithm.

This is confirmed by a quantitative and objective analysis, performed as follows. In order to verify the existence of an optimal threshold,
and to compare the quality of the visualizations with those obtained by filtered/unfiltered back-projection with different filters (and
shown in Figure \ref{radonA}), we have
\begin{itemize}
\item rescaled the grey levels of all the visualizations in the range $[0,1]$, as in the original Shepp--Logan phantom;
\item masked the pixels of the background in order to compare just the values of the pixels inside the phantom;
\item defined and computed the visualization error as the Frobenius norm of the matrices given by the differences between each visualization and the original Shepp--Logan phantom.
\end{itemize}

\begin{figure}[!htbp]
\begin{center}
\includegraphics[scale=0.3]{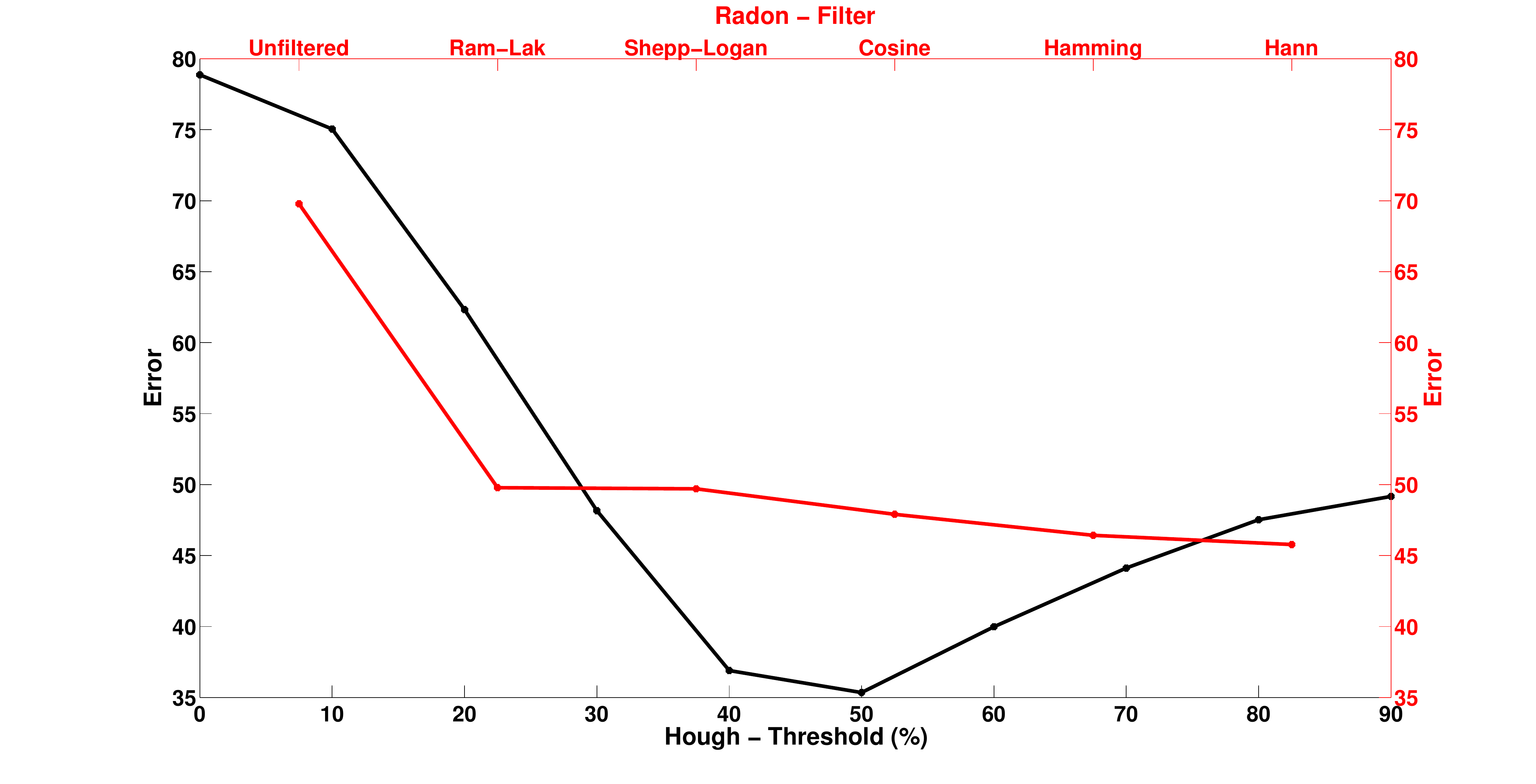}\\
\caption{Visualization errors committed in the inversion of a noisy Radon sinogram by using unfiltered/filtered back projection with different filtering functions (red plot and axes) and by using the Hough-based procedure (black plot and axes).}
\label{frobenius}
\end{center}
\end{figure}

The results of this analysis are summarized in the graph with multiple $x$- and $y$-axes shown in Figure \ref{frobenius}, where the black axes and plot refer to the errors computed from the Hough visualizations, while the red ones refer to back-projection rescaled reconstructions. The black plot clearly shows that there exists an optimal threshold value minimizing the Hough error function. Moreover, for a rather large range of threshold values, the visualizations obtained from Hough inversion seem to be more accurate than those obtained from usual back-projection inversion.

The black plot in Figure \ref{frobenius} also seems to suggest that the threshold may play the role of a regularization parameter in the
Hough-based inversion algorithm. This is one of several issues to be investigated in order to make this algorithm a feasible and accepted
alternative to the classical FBP, at least when the latter does not provide satisfactory results.

\section{Conclusions and future perspectives}\label{zioprospettico}

This paper provides for the first time a rigorous description of the formal equivalence between the Radon transform, introduced in harmonic analysis and at the basis of the mathematical theory of X-ray tomography, and the Hough transform, utilized in image processing for pattern recognition. Specifically, the main theoretical result of the paper is concerned with the forward problem associated to the Radon transform, i.e., the proof that the rescaled Hough counter of either a linear combination of Dirac deltas or a piecewise constant function tends to the Radon sinogram as the discretization step in the parameter space vanishes. Moreover, we briefly discussed how the Hough-Radon equivalence may have impacts on the inverse problem associated to image reconstruction in the case of modalities in which the data formation process is modeled by the Radon transform. Indeed, application perspectives of this paper are concerned with the possibility to invert a Radon sinogram by regarding it as a rescaled Hough sinogram. There is no doubt that the FBP algorithm is a very powerful tool for the inversion of Radon sinograms in X-ray computed tomography, but its performance can degrade in presence of very noisy sinograms, as in the case of Positron Emission Tomography (PET) imaging, where well known physical effects limit the spatial resolution. Exploiting the identification between Radon and Hough sinograms proved in this paper, for a simple and synthetic example we have here shown that improvements can be obtained by using the Hough procedure with an optimal threshold to invert a Radon sinogram. Further studies should be aimed at testing this computational method in realistic PET conditions.

\section*{Acknowledgements}

Our collegues and friends Filippo De Mari, Ernesto De Vito and Nicola Pinamonti (Universit\`{a} di Genova)
are kindly acknowledged for their valuable suggestions.



\newpage

\appendix

\section{Some elements of distributions theory}\label{distrapp}

The aim of this appendix is to introduce the notation adopted in the paper and to recall some definitions and results in distribution theory,
since we state them in a formulation that is sometimes different from the standard one.

\subsection{Test functions and functionals acting on them}\label{test}

Let $W$ be a non-empty open subset of $\R^n$, with $n\in\N\setminus\{0\}$.
For $k\in\N$ or $k=\infty$, let $C^k_C(W)$ be the vector space of $k$-times continuously differentiable functions $\phi:W\rightarrow\C$
whose supports $\mathrm{supp}\,\phi$ are compact subsets of $W$. To denote partial derivatives of such functions, we use the multi-index
$\alpha=(\alpha_1,\ldots,\alpha_n)\in\N^n$, whereby the
order of the differential operator is $|\alpha|:=\sum_{i=1}^{n}\alpha_i$. According to this notation, we have
$\partial^{\alpha}:=\partial_1^{\alpha_1}\cdots\partial_n^{\alpha_n}=\frac{\partial^{|\alpha|}}{\partial x_1^{\alpha_1}
\cdots\partial x_n^{\alpha_n}}$,
where in the right-hand side the components of $x=(x_1,\ldots,x_n)\in\R^n$ appear explicitly.
\begin{definition}\label{convDk}
For $W$ and $k$ as above, the space of test functions $\mathcal{D}_k(W)$ is the space $C^k_C(W)$ endowed with the following notion
of convergence of sequences. Given a sequence $\left\lbrace\phi_j\right\rbrace_{j\in\N}$
in $C^k_C(W)$, we say that $\phi_j\rightarrow 0$ in $\mathcal{D}_k(W)$ as $j\rightarrow\infty$
if and only if there exists a compact subset $K$ of $W$ such that $\mathrm{supp}\,\phi_j\subset K$ for all $j\in\N$ and
$\lim_{j\rightarrow\infty}\sup_{x\in K}|\partial^\alpha \phi_j(x)|=0$ for all $\alpha$ such that $|\alpha|\leq k$.
The convergence of a sequence $\left\lbrace\phi_j\right\rbrace_{j\in\N}\subset C^k_C(W)$ to a non-zero function $\phi\in C^k_C(W)$ is trivially defined by the condition $\left(\phi_j-\phi\right)\rightarrow 0$ in $\mathcal{D}_k(W)$.
\end{definition}

\begin{definition}\label{schelling}
For $W$ and $k$ as above, let $\ell:\mathcal{D}_k(W)\rightarrow\C$ be a functional
endowed with the two following properties:
\begin{itemize}
 \item[\textnormal{(i)}] linearity, i.e., $\langle\ell,a\phi_1 + b\phi_2\rangle = a\langle\ell,\phi_1\rangle + b \langle\ell,\phi_2\rangle$
 $\forall a, b\in\C$, $\forall\phi_1, \phi_2\in \mathcal{D}_k(W)$, having adopted the pairing notation $\langle\ell,\phi\rangle$
 to denote the action of $\ell$ on $\phi$, i.e., $\ell(\phi)$;
 \item[\textnormal{(ii)}] continuity, i.e., if
 $\phi_j\rightarrow\phi$ in $\mathcal{D}_k(W)$ as $j\rightarrow\infty$, then
 $\langle\ell,\phi_j\rangle\rightarrow \langle\ell,\phi\rangle$ in $\C$ as $j\rightarrow\infty$.
\end{itemize}
Such a functional $\ell$ is called a \textnormal{distribution\footnote{Usually, the term ``distribution'' is referred to the case $k=\infty$,
whereby $\mathcal{D}_\infty(W)$ is simply denoted by $\mathcal{D}(W)$. However, it is well known that distributions can be defined on
spaces of test functions that are larger than $\mathcal{D}(W)$: see, e.g., \cite[pp. 14--15]{zem87} or \cite{deki09}.}
(on $\mathcal{D}_k(W)$)}. The
set of all such functionals is denoted by $\mathcal{D}'_k(W)$ and is made a vector space by setting
\be\label{hume}
\langle a\ell_1+b\ell_2,\phi \rangle:=a\langle\ell_1,\phi\rangle + b \langle\ell_2,\phi\rangle\
\forall a,b\in\C,\ \forall\ell_1,\ell_2\in\mathcal{D}'_k(W),\ \forall\phi\in \mathcal{D}_k(W).
\ee
\end{definition}
\begin{definition}\label{hegel}
Let $\Xi$ be a subset of $\,\R$ such that $\bar{\xi}$ is an accumulation
point\footnote{Typically, $\bar{\xi}=0$ or, mainly when $\Xi=\N$, $\bar{\xi}=+\infty$.} for $\Xi$. For each $\xi\in\Xi$, let
$\ell_\xi\in \mathcal{D}'_k(W)$, with $W$, $k$ as above, and let
$\ell\in \mathcal{D}'_k(W)$. Then we say that $\ell_\xi\rightarrow\ell$
in $\mathcal{D}'_k(W)$ as $\xi\rightarrow\bar{\xi}$
if and only if
\be \label{kant}
\lim_{\xi\rightarrow\bar{\xi}}\,\langle\ell_\xi,\phi\rangle=\langle\ell,\phi\rangle\ \ \ \forall\phi\in\mathcal{D}_k(W).
\ee
\end{definition}

We observe that, by Lebesgue dominated convergence theorem, any locally integrable function $u\in L^1_{\mathrm{loc}}(W)$
uniquely defines a corresponding distribution $\iota_k (u)\in\mathcal{D}'_k(W)$ (for any $k\in\N$ or $k=\infty$) by means of the rule
\be \label{pairint}
\langle\iota_k (u),\phi\rangle:=\int_{W}u(x)\phi(x)\,dx\ \ \ \ \forall\phi\in\mathcal{D}_k(W).
\ee
It can be proved that the linear map $\iota_k:L^1_{\mathrm{loc}}(W)\rightarrow \mathcal{D}'_k(W)$ is injective: see, e.g.,
\cite[pp. 64--66]{mclean}. In particular, we can regard $L^1_{\mathrm{loc}}(W)$ as a subspace of $\mathcal{D}'_k(W)$, and then denote
$\iota_k$ as an inclusion map, i.e., $\iota_k:L^1_{\mathrm{loc}}(W)\hookrightarrow \mathcal{D}'_k(W)$.

Moreover, we recall that any distribution $\ell\in\mathcal{D}'_k(W)$ admits (distributional) partial derivatives
$\partial^\alpha\ell\in\mathcal{D}'_{k+|\alpha|}(W)$ of any order $|\alpha|\geq 0$,
according to the definition
\be \label{derdistr}
\langle\partial^\alpha\ell,\phi\rangle:=(-1)^{|\alpha|}\,\langle\ell,\partial^\alpha\phi\rangle
\ \ \ \ \forall\alpha\in\N^n,\ \ \forall\phi\in
\mathcal{D}_{k+|\alpha|}(W).
\ee

\begin{definition}\label{intparam}
For $p\in\N\setminus\{0\}$ and $k\in\N$ or $k=\infty$, let $Y$ be a Lebesgue-measurable subset of $\,\R^p$ and,
for each $y\in Y$, let $\ell_y\in\mathcal{D}'_k(W)$. Moreover, assume that
\begin{itemize}
 \item[\textnormal{(i)}] the map defined by $Y\ni y\mapsto \langle\ell_y,\phi\rangle\in\C$ belongs to
 $L^1(Y)$ $\forall\phi\in\mathcal{D}_k(W)$;
 \item[\textnormal{(ii)}] there exists $\ell\in\mathcal{D}'_k(W)$ such that
 $\int_Y \langle\ell_y,\phi\rangle\, dy = \langle\ell,\phi\rangle$
 $\forall\phi\in\mathcal{D}_k(W)$.
\end{itemize}
Then $\ell$ is said to be the integral of $\ell_y$ with respect to $y$, and
the equality $\ell= \int_Y \ell_y\,dy$ is written to summarize properties \textnormal{(i)} and \textnormal{(ii)}.
\end{definition}

\subsection{The Dirac delta of a function}\label{mementodelta}

First, we consider the one-dimensional case. Let $W$ be a non-empty open subset of $\R$ and let
$f:W\rightarrow\R$, with $f\in C^1(W)$. We denote by
\be\label{defS}
\mathcal{S}:=\{x\in W\, :\, f(x)=0\}
\ee
the locus of its zeros. If $\mathcal{S}=\emptyset$, we define the distribution $\delta(f)$ as coinciding with the zero
of the vector space $\mathcal{D}'_0(W)$. If $\mathcal{S}\neq\emptyset$, we assume that
1) $\mathcal{S}$ is at most countable: in particular, we can write
$\mathcal{S}=\{x_0(i)\in W : i\in I\}$, being the set $I$ of indices finite or countable; 2) the only possible
accumulation points for $\mathcal{S}$ are $-\infty$ and $+\infty$; 3) each zero $x_0(i)$ of $f$ is simple, i.e.,
$f'\big(x_0(i)\big)\neq 0$ $\forall i\in I$. Then, we define the distribution $\delta(f)\in\mathcal{D}'_0(W)$ as
\be \label{deltaf1}
\left\langle\delta(f),\phi\right\rangle:=\sum_{i\in I}\frac{\phi\big(x_0(i)\big)}{\left|f'\big(x_0(i)\big)\right|}
\ \ \ \ \ \forall\phi\in\mathcal{D}_0(W),
\ee
or, equivalently,
\be\label{ulisse}
\delta(f(x)):=\sum_{i\in I}\frac{1}{\left|f'\big(x_0(i)\big)\right|}\,\delta\big(x-x_0(i)\big),
\ee
being, in general, $\delta(\cdot - x_0)=\delta_{x_0}\in\mathcal{D}'_0(W)$ the Dirac delta centred at a point
$x_0\in \R$. For further details, see, e.g., \cite[pp. 184--185]{gesh1}.

Then, we turn to the $n$-dimensional case, with $n\geq 2$. Let $W$ be a non-empty open subset of $\R^n$,
with $n\in\N\setminus\{0,1\}$. Again, given
$f:W\rightarrow\R$, with $f\in C^1(W)$,
let $\mathcal{S}\subset W$ be the locus of its zeros, according to the same definition (\ref{defS}), provided that $x$ is now understood as
$x=(x_1,\ldots,x_n)\in\R^n$. As before, if $\mathcal{S}=\emptyset$, the distribution $\delta(f)$ coincides, by definition, with the zero of
$\mathcal{D}'_0(W)$.
If $\mathcal{S}\neq\emptyset$, we assume that
$\mathrm{grad}\, f(x)\neq 0$ $\forall x\in \mathcal{S}$:
such condition, together with $f\in C^1(W)$, allows defining on $\mathcal{S}$ a nowhere-vanishing differential form of maximum degree
\cite[chap.~III]{gesh1}, which amounts to prove (see, e.g., \cite[prop. 4.2]{wa83}) that $\mathcal{S}$ is an orientable manifold.
Summing up, the above assumptions on $f$ imply that $\mathcal{S}$
is a smooth, closed, orientable and $(n-1)$-dimensional submanifold
of $W$, described by the Cartesian equation $f(x)=0$.

Now, as explained, e.g., in
\cite[chap.~III]{gesh1}
and in
\cite[chap.~8]{jones}
(to which the reader is referred for a detailed treatment), $\delta(f)\in\mathcal{D}'_0(W)$ is defined as follows.

First, let us consider the particular case where $\mathcal{S}$ is the zero-locus of an $x_i$-solvable function, i.e.,
$\mathcal{S}$ is the graph of a function $\mathsf{F}$ of the coordinates
$x_1,\ldots,x_{i-1},x_{i+1},\ldots,x_n$; for notational simplicity, we assume $i=n$. Since $f\in C^1(W)$ by hypothesis, a sufficient
condition for $x_n$-solvability is that $\partial_n f(x)\neq 0$ $\forall x \in\mathcal{S}$. Indeed, in such case, if we define
the open subset $W'$ of $\R^{n-1}$ as $W':=\{x'\in \R^{n-1} : \exists\, x_n\in\R : (x',x_n)\in W\}$, by the implicit function theorem
there exists a function $\mathsf{F}\in C^1\left(W'\right)$ such that
\be\label{dini}
\mathcal{S}=\{x\in W\, :\, x_n=\mathsf{F}(x')\}.
\ee
We shall often denote the action of the functional $\delta(f)$ on the test function $\phi\in\mathcal{D}_0(W)$ by the integral notation
$\int_{W}\delta\left(f(x)\right)\phi(x)\,dx$, as synonymous with $\langle\delta(f),\phi\rangle$. This action is defined as
\be\label{defdelta}
\int_{W}\delta\left(f(x)\right)\phi(x)\,dx:=\int_{W'}\phi\left(x',\mathsf{F}(x')\right)
\frac{\sqrt{1 + \left|\mathrm{grad}\,\mathsf{F}(x')\right|^2}}
{\left|\mathrm{grad}\,f\left(x',\mathsf{F}(x')\right)\right|}\,dx',
\ee
where the right-hand side is a proper Lebesgue integral.

Note that, in general, the right-hand side of (\ref{defdelta}) does not coincide with the usual surface integral
$\int_{\mathcal{S}}\phi(x)\,d\sigma(x)$ of $\phi$ on $\mathcal{S}$, since
\be\label{intsup}
\int_{\mathcal{S}}\phi(x)\,d\sigma(x):=\int_{W'}\phi\left(x',\mathsf{F}(x')\right) \sqrt{1 + \left|\mathrm{grad}\,\mathsf{F}(x')\right|^2}\, dx',
\ee
being
\be\label{desigma}
d\sigma(x)=\sqrt{1 + \left|\mathrm{grad}\,\mathsf{F}(x')\right|^2}\, dx'
\ee
the Euclidean surface element on $\mathcal{S}$.
In particular, the right-hand side of (\ref{intsup}) does not depend on the particular function $f$ chosen to describe $\mathcal{S}$,
while that of (\ref{defdelta}) does. However, the two integrals coincide when
\be\label{grad1}
\left|\mathrm{grad}\,f(x)\right|=1\ \ \forall x\in\mathcal{S}.
\ee
Moreover, whenever $f$ is of the form $f(x)=x_n-\mathsf{F}(x')$, it holds that
\be\label{viotti}
\frac{\sqrt{1 + \left|\mathrm{grad}\,\mathsf{F}(x')\right|^2}}{\left|\mathrm{grad}\,f(x',\mathsf{F}(x'))\right|}=1\ \ \ \forall x'\in W',
\ee
which is a useful identity when inserted into the integral on the right-hand side of (\ref{defdelta}).

The simplification due to the $x_n$-solvability of $f$ is that $(x_1,\ldots,x_{n-1})$ become global coordinates on the whole $\mathcal{S}$.
When this condition does not hold, the local solvability of $f$ is anyway ensured, owing to the implicit function theorem and
the basic assumption that $\mathrm{grad}\,f(x)\neq 0$ $\forall x\in\mathcal{S}$. Then,
the local expressions representing on each coordinate chart the integrand function in (\ref{defdelta}) (or in (\ref{intsup})) can be
glued together by using the partition of unity, and the integral can be computed as a finite sum of integrals on each chart covering
$\mathcal{S}\cap S_{\phi}$, being $S_{\phi}:=\mathrm{supp}\,\phi$. For sake of notational simplicity, even in this
more general case we shall maintain expressions (\ref{defdelta})--({\ref{desigma}}), being understood that $(x_1,\ldots,x_{n-1})$ is only a symbol for one of the coordinate charts covering $\mathcal{S}\cap S_{\phi}$.

Finally, we recall the following result.
Let $\Theta(f):W\rightarrow\{0,1\}$ be the characteristic function of the region $\{x\in W : f(x)\geq 0\}$, i.e.,
\be\label{defHf1}
\Theta(f(x)):=\left\{
\begin{array}{ll}
0 & \mbox{for}\ x\in W\, :\, f(x)<0,\\
1 & \mbox{for}\ x\in W\, :\, f(x)\geq 0.
\end{array}
\right.
\ee
Obviously, $\Theta(f)=\Theta\circ f$, where $\Theta$ is the Heaviside function, i.e., $\Theta(t)=0$ for $t<0$ and $\Theta(t)=1$ for $t\geq 0$.
Note that $\Theta(f)\in L^1_{\mathrm{loc}}(W)$: then, recalling the inclusion map
$\iota_k:L^1_{\mathrm{loc}}(W)\hookrightarrow \mathcal{D}'_k(W)$ defined by (\ref{pairint}),
$\Theta(f)$ can also be regarded as an element of $\mathcal{D}_k'(W)$, i.e., as $\iota_k\big(\Theta(f)\big)$,
which, according to (\ref{derdistr}), admits partial derivatives
with respect to $x_i$, for all $i=1,\ldots,n$.
In particular, it can be proved that (see \cite{gegr5}, p. 224), for any $k\in\N$ or $k=\infty$,
\be\label{derH1}
\frac{\partial\,\iota_k\big(\Theta(f)\big)}{\partial x_i}=\frac{\partial f}{\partial x_i}\,\delta(f)\in\mathcal{D}'_{k+1}(W)
\ \ \ \ \ \ \forall i\in\{1,\ldots,n\}.
\ee

\subsection{Piecewise continuous test functions}\label{ziopiecewise}

As definition (\ref{defdelta}) clearly shows, the fact that $\phi$ belongs to $\mathcal{D}_0(W)$
is a sufficient but not necessary condition for the existence of
$\left\langle\delta(f),\phi\right\rangle$. In fact, we can enlarge the spaces $\mathcal{D}_k(W)$ of test functions in order to
include functions that are piecewise $C^k(W)$. To this end, we introduce the following definitions.
\begin{definition}\label{pezzi}
Given an open subset $W$ of $\,\R^n$, with $n\in\N\setminus\{0\}$,
we say that a function $\phi:W\rightarrow\C$ is \textit{piecewise} $C^k(W)$ (for $k\in\N$ or $k=\infty$) if and only if
there exist a finite or countable set $Q$ of indices $q$ and a corresponding family $\left\lbrace W_q\right\rbrace_{q\in Q}$
of open subsets $W_q$ of $W$ such that
\begin{itemize}
\item[\textnormal{(i)}] $\mathcal{L}^{n}\left(\partial W_q\right)=0$ $\,\forall q\in Q$, where $\mathcal{L}^{n}$ denotes the
 Lebesgue measure on $\R^{n}$;
 \item[\textnormal{(ii)}] $W_r\cap W_s=\emptyset$ $\forall r,s\in Q$ such that $r\neq s$, and $\overline{W}=\bigcup_{q\in Q} \overline{W}_q$;
 \item[\textnormal{(iii)}] $\phi\in C^k\left(\overline{W_q\cap W}\right)$ $\forall q\in Q$: this notation means that, for each $q\in Q$,
 $\phi\in C^k\left(W_q\cap W\right)$ and both $\phi$ and its partial derivatives up to the $k$-th order are bounded on the open set $W_q\cap W$;
 \item[\textnormal{(iv)}] for any bounded subset $X$ of $\,W$, there exists $Q'\subset Q$ such that $\# Q'\in\N$ and
 $X\subset \bigcup_{q\in Q'} W_q$.
\end{itemize}
The vector space of all such functions will be denoted by $PC^k(W)$. The vector subspace of all such functions whose supports are
compact subsets of $W$ will be denoted by $PC^k_C(W)$.
\end{definition}

\begin{definition}\label{convDkpezzi}
For $W$ and $k$ as above, the space of test functions $\mathcal{PD}_k(W)$ is the space $PC^k_{\mathrm{comp}}(W)$
endowed with the following notion of convergence of sequences.
Given a sequence $\left\lbrace\phi_j\right\rbrace_{j\in\N}$ in $PC^k_{\mathrm{comp}}(W)$,
we say that $\phi_j\rightarrow 0$ in $\mathcal{PD}_k(W)$ as $j\rightarrow\infty$
if and only if the following three conditions are fulfilled:
\begin{itemize}
\item[\textnormal{(i)}] there exists a compact subset $K$ of $W$ such that $\mathrm{supp}\,\phi_j\subset K$ for all $j\in\N$;
\item[\textnormal{(ii)}] $\displaystyle\lim_{j\rightarrow\infty}\sup_{x\in K}|\phi_j(x)|=0$;
\item[\textnormal{(iii)}] $\exists Z\subset K$ such that $\mathcal{L}^n(Z)=0$ and
$\displaystyle\lim_{j\rightarrow\infty}\sup_{x\in K\setminus Z}|\partial^\alpha \phi_j(x)|=0$ $\ \forall\alpha\,:\,|\alpha|\leq k$.
\end{itemize}
The convergence of a sequence $\left\lbrace\phi_j\right\rbrace_{j\in\N}\subset PC^k_C(W)$ to a non-zero function $\phi\in PC^k_C(W)$ is trivially defined by the condition $\left(\phi_j-\phi\right)\rightarrow 0$ in $\mathcal{PD}_k(W)$.
\end{definition}

All the remaining part of \ref{test}, except definition (\ref{derdistr}), can now be trivially paralleled,
by replacing $\mathcal{D}_k(W)$ with $\mathcal{PD}_k(W)$. For example, Definition \ref{schelling} allows introducing the vector space
$\mathcal{PD}'_k(W)$ of linear and continuous functionals on $\mathcal{PD}_k(W)$
(note that $\mathcal{D}_k(W)\subset\mathcal{PD}_k(W)$
implies $\mathcal{D}'_k(W)\supset\mathcal{PD}'_k(W)$); rule (\ref{pairint}) defines the inclusion map
$\widetilde{\iota}_k:L^1_{\mathrm{loc}}(W)\hookrightarrow\mathcal{PD}'_k(W)$; Definition \ref{intparam} establishes the concept
of integral of a $y$-parametrized family of distributions $\ell_y$ with respect to the parameter $y\in Y\subset\R^p$.

Of course, defining an
appropriate concept of derivation in $\mathcal{PD}'_k(W)$ is not easy, but since we do not need to
consider such an operation in our paper, we shall not address this issue (we refer, e.g., to \cite{deki09}
for a theory of distributions acting on discontinuous test functions).

Finally, it is not difficult to realize that, by taking $\phi\in\mathcal{PD}_0(W)$, definition (\ref{defdelta}) can still be adopted in order to make $\delta(f)$ an element of $\mathcal{PD}'_0(W)$.

\subsection{A link between $\delta(f)$ and the coarea formula}\label{ziocoarea}

Interestingly, definition (\ref{defdelta}) is related to the \textit{coarea formula}: as we are going to see, this relationship is a
consequence of the following Theorem \ref{devito}, which is a corollary of the result known as the ``coarea formula'' and is proved
in \cite[pp. 118--119]{evga92}.

\begin{theorem}\label{devito}
Let $A$ be an open subset of $\R^n$, with $n\in\N\setminus\{0\}$, and let $\Psi:A\rightarrow\R$ be a Lipschitz function
such that $\mathrm{ess} \inf |\mathrm{grad}\, \Psi|>0$.
Moreover, let $g:A\rightarrow\R$ be such that $g\in L^1(A)$.
Then, it holds that
\be \label{coarea}
\int_{A}g(x)\,dx=\int_{\Psi(A)}\left(\int_{\Psi^{-1}(s)}\frac{g(x)}{|\mathrm{grad}\, \Psi(x)|}\,d\sigma(x)\right) ds,
\ee
where $d\sigma(x)$ is to be understood as either the Euclidean surface element on $\Psi^{-1}(s)$ if $n\geq 2$, or
the counting measure on the (at most countable) set $\Psi^{-1}(s)$ if $n=1$.
\end{theorem}

Now, for $n\geq 2$, it is clear that definition (\ref{defdelta}) coincides with the internal integral on the right-hand side of
(\ref{coarea}) under the
identifications $A=W$, $\Psi=f$, $g=\phi\in\mathcal{D}_0(W)$ and $s=0$. The only apparent mismatches are that in (\ref{defdelta})
the function $f$ is not required to be Lipschitz and, rather than
$\mathrm{ess} \inf |\mathrm{grad}\, f|>0$, only the weaker
condition $|\mathrm{grad}\, f(x)|\neq 0$ $\forall x\in f^{-1}(0)=\mathcal{S}$ is assumed.
However, in (\ref{defdelta}) the integral on $\mathcal{S}$ can be regarded as performed on the compact set
$K_{\phi}:=S_\phi\cap\mathcal{S}$, being $S_\phi:=\mathrm{supp}\,\phi$.
Since $S_\phi\subset W$ is compact and $f\in C^1(W)$, there exists an open set $V$ such that
$S_\phi\subset V \subset W$ and the restriction $f\restriction_V$ is Lipschitz.
Then, by Kirszbraun theorem (see \cite{federer69}, p. 201),
there exists a Lipschitz function
$\bar{f}:\R^n\rightarrow\R$ with the property $\bar{f}\restriction_V=f\restriction_V$,
so that we can replace $f$ with $\bar{f}\restriction_V$ in (\ref{defdelta}). Moreover,
since $\bar{f}\restriction_V\in C^1(V)$, we have that
$\left|\mathrm{grad}\, \bar{f}\restriction_{V}\right|\in C^0(V)$. Thus, condition
$|\mathrm{grad}\, f(x)|\neq 0$ $\forall x\in \mathcal{S}$ in (\ref{defdelta}) implies that
$\mathrm{ess} \inf_{x\in K_{\phi}} \left|\mathrm{grad}\, \bar{f}\restriction_{V}(x)\right|=
\inf_{x\in K_{\phi}} \left|\mathrm{grad}\, \bar{f}\restriction_{V}(x)\right|=\min_{x\in K_{\phi}} \left|\mathrm{grad}\, \bar{f}\restriction_{V}(x)\right|>0$,
since the continuous function $\left|\mathrm{grad}\, \bar{f}\restriction_{V}\right|$ admits maximum and minimum values on the compact set $K_{\phi}$.

The previous argument can also be adapted to the case $n=1$. Indeed, maintaining the above identifications, it is easy to realize that the internal integral on the right-hand side of (\ref{coarea}), computed with respect to the counting measure on
$\mathcal{S}=\{x_0(i)\in W : i\in I\}$, coincides with definition (\ref{deltaf1}).



\bibliographystyle{plain}
\bibliography{Hough_Radon_ACHA_2016}







\end{document}